\documentclass[11pt]{article}
\usepackage[margin=1in]{geometry}
\usepackage{amsthm}

\usepackage[numbers,sort&compress]{natbib}

\usepackage{braket}
\usepackage[utf8]{inputenc}
\usepackage[english]{babel}

\usepackage{tikz}
\usetikzlibrary{shapes.geometric}

\usepackage{graphicx}
\usepackage{stmaryrd}
\usepackage[normalem]{ulem}

\usepackage{amsfonts}
\usepackage{amssymb}
\usepackage{amsmath}

\usepackage{mathtools}
\usepackage{enumitem}

\usepackage{xcolor}

\usepackage{mathrsfs}
\usepackage{overpic}

\usepackage{mathtools,booktabs}

\DeclarePairedDelimiter\floor{\lfloor}{\rfloor}

\usepackage{varwidth}

\newcommand{\rank}{{\rm rank}}

\graphicspath{{images/}}

\newcommand{\R}[0]{\mathbb{R}}
\newcommand{\C}[0]{\mathbb{C}}

\usepackage{aliascnt}

\newtheorem{theorem}{Theorem}[section]

\newaliascnt{lemma}{theorem}
\newtheorem{lemma}[lemma]{Lemma}
\aliascntresetthe{lemma}

\newaliascnt{corollary}{theorem}

\aliascntresetthe{corollary}

\newaliascnt{proposition}{theorem}

\aliascntresetthe{proposition}

\newaliascnt{definition}{theorem}

\aliascntresetthe{definition}

\newaliascnt{assumption}{theorem}
\newtheorem{assumption}[assumption]{Assumption}
\aliascntresetthe{assumption}

\newaliascnt{remark}{theorem}
\newtheorem{remark}[remark]{Remark}
\aliascntresetthe{remark}

\newaliascnt{example}{theorem}

\aliascntresetthe{example}

\usepackage{algorithm}
\numberwithin{algorithm}{section}

\usepackage[colorlinks=true,linkcolor=blue,citecolor=blue,urlcolor=blue]{hyperref}
\usepackage[capitalise,nameinlink]{cleveref}
\crefname{equation}{}{}
\creflabelformat{equation}{#2\textup{(#1)}#3}
\crefname{section}{section}{sections}
\crefname{subsection}{subsection}{subsections}
\crefname{algorithm}{Algorithm}{Algorithms}
\Crefname{algorithm}{Algorithm}{Algorithms}
\crefname{theorem}{Theorem}{Theorems}
\Crefname{theorem}{Theorem}{Theorems}
\crefname{lemma}{Lemma}{Lemmas}
\Crefname{lemma}{Lemma}{Lemmas}
\crefname{corollary}{Corollary}{Corollaries}
\Crefname{corollary}{Corollary}{Corollaries}
\crefname{proposition}{Proposition}{Propositions}
\Crefname{proposition}{Proposition}{Propositions}
\crefname{definition}{Definition}{Definitions}
\Crefname{definition}{Definition}{Definitions}
\crefname{assumption}{Assumption}{Assumptions}
\Crefname{assumption}{Assumption}{Assumptions}
\crefname{remark}{Remark}{Remarks}
\Crefname{remark}{Remark}{Remarks}

\usepackage[noend]{algpseudocode}
\algrenewcommand\algorithmicrequire{\textbf{Input:}}
\algrenewcommand\algorithmicensure{\textbf{Output:}}

\makeatletter
\let\arxivorigappendix\appendix
\renewcommand{\appendix}{
  \arxivorigappendix
  \numberwithin{equation}{section}
  \numberwithin{figure}{section}
  \numberwithin{table}{section}
}
\makeatother

\title{An Efficient Projection Enhanced Interpolation Based Tensor Train Decomposition\thanks{D. Hayes and J.-M. Qiu acknowledge support provided by Department of Energy DE-SC0023164 and the Patient-Centered Outcomes Research Institute (PCORI) under award ME-2024C3-42660. J.-M. Qiu also acknowledges support provided by Air Force Office of Scientific Research FA9550-24-1-0254 and the National Institute of Mental Health of the National Institutes of Health under award R01MH142790. T. Shi acknowledges support provided by the Director, Office of Science, Office of Advanced Scientific Computing Research, of the U.S. Department of Energy under Contract No. DE-AC02-05CH11231. This research used resources of the National Energy Research Scientific Computing Center (NERSC), a Department of Energy User Facility using NERSC award under the Contract No. DE-AC02-05CH11231. This research was supported in part through the use of DARWIN computing system at the University of Delaware supported by NSF under Grant Number 1919839.}}
\author{Daniel Hayes\thanks{Department of Mathematical Sciences, University of Delaware, Newark, DE 19716. \texttt{dphayes@udel.edu}}
\and Jing-Mei Qiu\thanks{Department of Mathematical Sciences, University of Delaware, Newark, DE 19716. \texttt{jingqiu@udel.edu}}
\and Tianyi Shi\thanks{Scalable Solvers Group, Lawrence Berkeley National Laboratory, Berkeley, CA 94720. \texttt{tianyishi@lbl.gov}}}
\date{}
\begin{document}
\maketitle
\begin{abstract}
The tensor-train (TT) format is a data-sparse tensor representation commonly used in high dimensional data approximations. In order to represent data with interpretability in data science, researchers develop data-centric skeletonized low rank approximations. However, these methods might still suffer from accuracy degradation, nonrobustness, and high computation costs. In this paper, given existing skeletonized TT approximations, we propose a family of projection enhanced
interpolation based algorithms to further improve approximation accuracy while keeping low computational complexity. We do this as a postprocessing step to existing interpolative decompositions, via oversampling data not in skeletons to include more information and selecting subsets of pivots for faster projections. We provide a priori error bounds for the projection enhanced TT format and for the proposed dimension parallel algorithm, in which the contributions of the skeleton, the oversampling, and the oblique projections appear as separate terms. We illustrate the performances of our proposed methods with extensive numerical experiments. These include up to 10D synthetic datasets such as tensors generated from kernel functions, and 4D and 6D tensors constructed from Maxwellian distribution functions that arise in kinetic theory. Our results demonstrate accuracy improvements between 2 and 100 times on tested tensors over original skeletonized TT approximations, while using a limited amount of computational resources.
\end{abstract}


\section{Introduction}
\label{sec:introduction}
A wide range of applications involve multidimensional data as observations or solutions. These datasets are often referred to as tensors~\cite{kolda2009tensor,ballard2025tensor}, and a general tensor $\mathcal{X} \in \C^{n_1\times\cdots \times n_d}$ requires a storage cost of $\prod_{j=1}^d n_j$ degrees of freedom. This scales exponentially with the dimension $d$, and is often referred to as ``the curse of dimensionality" (CoD). Tensors are higher order analogues of matrices, so a common way to alleviate CoD is to represent tensors with low rank approximations, similar to data-sparse matrix approximations. The Eckart-Young theorem~\cite{eckart1936approximation} bounds singular values of a matrix, thus guaranteeing optimal low rank approximation through the singular value decomposition (SVD), but it does not hold for multi-linear algebra. As a result, researchers develop a handful of low rank tensor formats, including canonical polyadic (CP)~\cite{hitchcock1927expression}, Tucker~\cite{de2000multilinear,tucker1966some} and hierarchical Tucker~\cite{grasedyck2010hierarchical}, tensor-train (TT)~\cite{oseledets2011tensor} and quantized TT~\cite{dolgov2012fast}, and tensor networks~\cite{evenbly2011tensor} with more general network topologies. In particular, because of the close connections with matrices for algorithmic designs and linear scaling with respect to storage cost, the TT format, also known as the matrix product state (MPS) in tensor networks and quantum physics, is widely used in applications such as molecular simulations~\cite{savostyanov2014exact}, high-order correlation functions~\cite{kressner2015low}, partial differential equations~\cite{guo2023local, dektor2025collocation, einkemmer2025review}, constrained optimization~\cite{dolgov2017low,benner2020low}, and machine learning~\cite{vandereycken2022ttml,novikov2020tensor}.

Although SVD provides best low rank matrix approximations, it is computationally expensive. In practice, for low rank approximations with accuracy bounds, researchers use other sub-optimal algorithms, including randomized SVD~\cite{halko2011finding}, matrix sketching via randomized range finder~\cite{halko2011finding}, rank-revealing QR~\cite{chan1987rank}, and so on. These methods are all based on operations to project a full rank matrix onto lower-dimensional subspaces, and the matrices to represent the projector have orthonormal (ON) columns. These methods, including SVD, are often referred to as orthogonal factorizations or projection-based factorizations. By applying these projection-based factorizations on tensor unfoldings, one can build low rank TT approximations with deterministic and randomized TTSVD~\cite{oseledets2011tensor,che2019randomized} and randomized TT sketching~\cite{shi2023parallel, kressner2023streaming}.

In some application areas such as data analysis, researchers sometimes want to represent data as a linear combination of a selected subset of data for better interpretability. This leads to skeletonized or interpolation-based matrix and tensor low rank approximations, where the low rank factors are constructed with a subset of the data structure. In the matrix regime, such algorithms include the interpolative decomposition (ID)~\cite{cheng2005compression,gu1996efficient}, CUR~\cite{mahoney2009cur, hamm2020perspectives}, adaptive cross approximation (ACA)~\cite{zhao2005adaptive}, and others. For tensors, and specifically the TT format, the commonly used ones are cross (TT-cross)~\cite{oseledets2010tt,dektor2025collocation,savostyanov2011fast,dolgov2014alternating} and adaptive cross (TTACA)~\cite{oseledets2010tt,dolgov2020parallel,shi2026distributed}. Unlike projection-based algorithms, interpolation-based algorithms do not need to access all elements of the data. They provide sub-optimal low rank approximations, but avoid CoD especially for high dimensional tensors. The key procedure in the skeletonized algorithms is to find the ``pivots", or ``skeletons" of the data, i.e., important indices to represent a given dimension of the matrix or tensor. There are various criteria and methods for the pivot selection, including maximum volume~\cite{goreinov1997theory,goreinov2010find}, statistical leverage scores~\cite{mahoney2009cur}, discrete empirical interpolation method (DEIM) from model order reduction~\cite{sorensen2016deim,chaturantabut2010nonlinear}, greedy approaches based on approximation differences~\cite{davis1997adaptive,dong2023simpler,dolgov2020parallel}, among others. For the TT format, these criteria lead to different cross algorithms. The original TT-cross~\cite{oseledets2010tt} uses alternating maximum-volume sweeps with prescribed ranks, the DMRG-type cross of~\cite{savostyanov2011fast} and the AMEn-type cross implemented in the TT-Toolbox~\cite{dolgov2014alternating,OseledetsTTToolbox} adapt the ranks during the sweeps, and the greedy cross of~\cite{savostyanov2014quasioptimality,dolgov2020parallel,shi2026distributed} adds one pivot at a time where the current error is largest. In practice for matrices, we observe that pivots selected via different approaches share a large overlap, and can all lead to reasonable matrix approximations. Therefore, in this work, we focus on methods that only need to traverse through linearly-scaled number of tensor entries with respect to the mode size, such as the greedy method. These approaches are thus termed to have a linearly-scaled complexity, and indicate that they are suitable for large datasets.

In spite of computational efficiency, interpolation-based approximations might suffer from instability and non-robustness from noisy data \cite{peherstorfer2020stability}. One remedy is to use the skeletonized decomposition as an initial guess for constructing orthogonal factors, from which one can construct a projection-based approximation. Such approximation becomes a post-processing step of an existing interpolation-based decomposition for robustness and accuracy. It is the focus of this manuscript, and we call such a method ``projection enhanced interpolation-based decomposition" (PEID). The baseline routine of PEID for matrices can be briefly described as the following: we find bases with ON columns for the skeletonized low rank factors using QR-preconditioned SVD~\cite{drmavc2017algorithm}, and then we form a low rank approximation using projections of the orthogonal factors (see~\cref{sec:CUR_oblique_sub} for more details). We extend this method to work on the TT format for tensors, with modifications so that the algorithms can be efficiently and accurately applied to high-dimensional data.

PEID methods work well when a linearly-scaled ACA is used as the interpolation-based initial guess, so that one can afford computations from projections to improve accuracy. However, PEID faces some potential challenges:
\begin{enumerate}[leftmargin=*,noitemsep]
\item The accuracy of the approximation depends heavily on the quality of the chosen skeleton in representing full data. This is particularly true for linearly-scaled ACA, where heuristics or greedy criterion are used to ensure efficiency in interpolation, so some important information might be hiding in data not selected.
\item The overall complexity of PEID is bounded by projections, which includes finding basis with ON columns through QR, SVD, etc., and conducting projections via matrix-matrix multiplication. This becomes particularly troublesome for high dimensional tensors, as these operations involve matrices related to tensor unfoldings. Therefore, using full matrices directly result in ``full-complexity" methods, whose algorithmic complexity is dominated by high-degree polynomials of the mode size, and it incurs CoD.
\end{enumerate}
In this manuscript, we propose two strategies in addressing the above mentioned challenges, with details described in~\cref{sec:CUR_oblique_sub}.
\begin{itemize}[leftmargin=*,noitemsep]
\item We oversample indices from unselected skeletons~\cite{cortinovis2024sublinear} to improve interpolation accuracy. This is similar to the oversampled CUR~\cite{park2025accuracy} and Gappy proper orthogonal decomposition (GappyPOD)~\cite{peherstorfer2020stability} from DEIM~\cite{chaturantabut2010nonlinear}. Since we use low-cost linearly-scaled ACA, we simply use uniform oversampling with little computational overhead. In comparison, in~\cite{peherstorfer2020stability}, DEIM with uniform oversampling is referred to as GappyPOD$+$R, and oversampling generated via singular value information is GappyPOD$+$E. Although GappyPOD$+$E is the winner technique in terms of accuracy improvement, it is significantly more expensive to compute.
\item We select submatrices of the orthogonal factors to perform oblique projections~\cite{engquist2007fast}, in order to reduce the complexity of full orthogonal projections. As ACA skeletons are good representations of the dataset, we propose to use them for selection of effective submatrices. See right panel of~\cref{fig:examples24} for a comparison of oblique projection with pivots and full orthogonal projection in terms of accuracy for matrix approximations.
\end{itemize}

In this paper, we propose a family of PEID algorithms for high-dimensional tensors in the TT format. Such algorithms oversample from indices uncovered by skeletonized pivots, and select most informative submatrices for oblique projections across all tensor unfoldings. The proposed algorithms can be thought of as a low-cost post-processing step for skeletonized TT approximation to further improve accuracy, regardless of whether the skeleton comes from a greedy, DMRG-type, or AMEn-type cross.
Our proposed methods are designed with similar oversampling ideas, but they differ on whether the oversampled indices are nested (see~\cref{sec:cross_algs} for the definition of nestedness). To be specific, \textbf{TT-PEID-Par} (\cref{alg:TT_update_parallel}) and \textbf{TT-PEID-Par2} (\cref{alg:TT_update_parallel_2sided}) require the oversampled indices to be nested, so building TT cores can be conducted in a dimension parallel manner, but the choices for oversampling are significantly limited. Comparatively, \textbf{TT-PEID-Seq} (\cref{alg:TT_update_sequential}) and \textbf{TT-PEID-Seq2} (\cref{alg:TT_update_sequential_2sided}) do not have this nestedness requirement, leading to sequential algorithm execution and broader oversampling choices. Two-sided versions (\textbf{TT-PEID-Par2} and \textbf{TT-PEID-Seq2}) can be viewed as extensions of matrix PEID, where orthogonal bases are generated for both row and column spaces, such as CUR or two-sided ID. In contrast, one-sided TT algorithms (\textbf{TT-PEID-Par} and \textbf{TT-PEID-Seq}) are higher-order analogues of matrix PEID that uses ID as the interpolation-based approximation.
In addition, \textbf{TT-PEID-Par} and \textbf{TT-PEID-Par2} require further that the skeleton pivots are nested so we can reuse the pivot rows of consecutive unfoldings without any reconstruction, but \textbf{TT-PEID-Seq} and \textbf{TT-PEID-Seq2} do not. We present error analysis of the underlying projection enhanced format,~\cref{thm:PEID_format_main,thm:PEID_format_simplified}, which requires neither skeleton pivots nestedness nor oversampling sets nestedness.
Finally, if we have sufficient computational power, we design \textbf{TT-PEID-R} (\cref{alg:TT_update_rounding}) that uses two one-sided runs to get a better approximation via rounding.

The remainder of the manuscript is organized as follows.~\Cref{sec:preliminary} reviews some necessary tensor notations, the TT format, skeletonized approximations, and details of matrix PEID with accuracy and stability improvements. In~\cref{sec:algorithm}, we introduce new PEID algorithms for the TT format. Finally, we demonstrate the performance on a variety of synthetic test tensors in~\cref{sec:numerical_tests}.

\section{Tensor notations and basic algorithms for matrix and tensor decompositions}
\label{sec:preliminary}
In this section, we review some tensor notations, the TT format for low rank tensor approximations, skeletonized matrix and tensor decompositions, and the PEID method with improvements for matrix approximations.

\subsection{Tensor notations} \label{sec:notation}
We use lower case letters for vectors, capital letters for matrices, and calligraphic capital letters for tensors. The number of entries for one dimension of a tensor is referred to as the mode size. We use MATLAB-style symbol ``:" to represent all the indices in one specific dimension. For example, if $\mathcal{Y}$ is a 3D tensor, then $\mathcal{Y}(:,j,:)$ denotes the $j$th lateral slice of $\mathcal{Y}$, and $\mathcal{Y}(:,:,k)$ denotes its $k$th mode-3 fiber. In addition, we use calligraphic capital letters for index sets in submatrix and subtensor selection. For example, $A(\mathcal{L},:)$ is a submatrix of $A$ with rows of $A$ selected from the set $\mathcal{L}$. 

We use the MATLAB command ``reshape" for the transformation between matrices and tensors. The new structure is constructed according to the desired multi indexing, without changing the element ordering. For example, if $\mathcal{Y} \in \C^{n_1 \times n_2 \times n_3}$, then $Z = {\rm reshape}(\mathcal{Y},n_1n_2,n_3)$ is a matricized output of size $n_1n_2 \times n_3$. Each column $Z(:,j)$ of $Z$ constitutes of elements from the corresponding frontal slice $\mathcal{Y}(:,:,j)$. In reverse, we have $\mathcal{Y} = {\rm reshape}(Z,n_1,n_2,n_3)$. For a $d$-dimensional tensor, there are $d-1$ ways to flatten the tensor into matrices via ``reshaping", and we refer to the result matrices as the unfoldings of the tensor. The majority of algorithms related to the TT tensor format relies on computations on the tensor unfoldings. Specifically, for a tensor $\mathcal{X}\in\C^{n_1\times\cdots\times n_d}$, we denote the $k$th unfolding as
\[X_k={\rm reshape}\left(\mathcal{X},\prod_{s=1}^k n_s,\prod_{s=k+1}^d n_s\right), \quad 1 \le k \le d-1.\]
Note that this is different from the tensor ``matricizations" $X_{(k)} \in \C^{n_k \times \prod_{j \neq k}n_j}$ used often in the Tucker format, where the row indices come directly from the $k$th tensor dimension, and the ordering of the entries changes between $\mathcal{X}$ and $X_{(k)}$.

\subsection{Tensor-train format} 
\label{sec:TT_format}
The TT format of a tensor $\mathcal{X}\in\C^{n_1\times \cdots \times n_d}$ is represented as
\[
\mathcal{X}_{i_1,\ldots,i_d} = \mathcal{G}_1(:,i_1,:)\mathcal{G}_2(:,i_2,:) \cdots \mathcal{G}_d(:,i_d,:), \qquad 1\leq i_k \leq n_k,
\]
where $\mathcal{G}_k \in \C^{s_{k-1} \times n_k \times s_k}$ for $1 \le k \le d$ are called the TT cores, and the vector $\pmb{s} = (s_0,\ldots,s_d)$ is referred to as the size of the TT cores. For the formula to hold, we require $s_0 = s_d = 1$.~\cref{fig:TT} illustrates one TT format with TT core size $\pmb{s}$ of $\mathcal{X}$. It is easy to see that this TT representation has a storage cost of $\sum_{k=1}^d s_{k-1}s_k n_k$, which is linear with respect to both $d$ and $n_k$. The storage cost is optimal if $s_k$ is as small as possible for each $1 \le k \le d-1$, and thus we define the TT rank as a vector $\pmb{r}= (r_0,\ldots,r_d)$ containing entry-by-entry smallest values of the TT core size. It is shown in~\cite{oseledets2011tensor} that if the TT representation of $\mathcal{X}$ is exact, then
\[ r_k = {\rm rank}(X_k), \quad 1 \le k \le d-1,\]
where $\rank(X_k)$ is the rank of the $k$th unfolding of $\mathcal{X}$. In practice, we often only use the numerical rank $r_k(\epsilon)$ of each $X_k$, which is the number of singular values of $X_k$ that are larger than $\epsilon$. As a result, we aim for $\pmb{s}$ that
\begin{equation} \label{eq:TT_trivial}
r_k(\epsilon) \le s_k \le {\rm rank}(X_k), \quad 1 \le k \le d-1.
\end{equation}
There are generally two types of algorithms regarding the choice of $\pmb{s}$: (1) rank-adaptive scheme, where $\pmb{s}$ is computed alongside the factorization steps using a given accuracy threshold; and (2) fixed-rank scheme, where $\pmb{s}$ is prior information known to the users and serve as an input to the algorithms. Algorithms in category (1) are suitable for theoretical analysis, such as TTSVD~\cite{oseledets2011tensor}, while algorithms in category (2) often excel in practical performance with fixed $\pmb{s}$. Some examples of these algorithms include TT sketching and randomized TTSVD~\cite{che2019randomized,shi2023parallel}, and TT cross~\cite{oseledets2010tt,dolgov2020parallel,shi2026distributed}.
\begin{figure}[!htbp]
\centering
\begin{tikzpicture}
\filldraw[black] (0,-0.5) node {$\mathcal{X}_{i_1,\ldots,i_d}$};
\filldraw[black] (1,-0.5) node {$=$};
\filldraw[color=black,fill=gray!20] (1.5,0) rectangle (4,-.5);
\filldraw[black] (2.8,-0.25) node {$\mathcal{G}_1(i_1,:)$};
\filldraw[black] (2.8,0.3) node {$s_0 \! \times \! s_1$};
\filldraw[color=black,fill=gray!20] (4.2,0) rectangle (7.1,-2.5);
\filldraw[black] (5.7,-1.3) node {$\mathcal{G}_2(:,i_2,:)$};
\filldraw[black] (5.7,0.3) node {$s_1 \! \times \! s_2$};
\filldraw[black] (7.5,-1) node {$\cdots$};
\filldraw[color=black,fill=gray!20] (7.9,0) rectangle (10.3,-2.3);
 \filldraw[black] (9.1,-1) node {$\mathcal{G}_{d-1}(:,i_{d-1},:)$};
\filldraw[black] (9.1,0.3) node {$s_{d-2} \! \times \! s_{d-1}$};
\filldraw[color=black,fill=gray!20] (10.5,0) rectangle (11,-2.4);
\filldraw[black] (10.75,-1) node {\rotatebox{270}{$\mathcal{G}_{d}(:,i_{d})$}};
\filldraw[black] (10.75,0.3) node {$s_{d-1} \! \times \! s_d$};
\end{tikzpicture}
\caption{The TT format with TT core size $\pmb{s} = (s_0,\ldots,s_d)$. Each entry of a tensor is represented by the product of $d$ matrices, where the $k$th matrix in the ``train" is selected based on the value of $i_k$.}
\label{fig:TT}
\end{figure}

For completeness, we briefly outline TT-Sketching~\cite{che2019randomized} (see~\cref{alg:TTsketching}) and Parallel-TT-Sketching~\cite{shi2023parallel} (see~\cref{alg:parallel_TTsketching}) with fixed ranks and no oversampling parameter, so that it is easier for our new proposed algorithm discussion in~\cref{sec:algorithm}.
In short, TT-Sketching finds TT cores by sequentially peeling off dimensions, while Parallel-TT-Sketching constructs TT cores by exploring the connections between column spaces of consecutive tensor unfoldings. Note that~\cref{alg:parallel_TTsketching} is parallel as respective iterations in line 1 and 3 can be carried out simultaneously.

\subsection{Skeletonized matrix and tensor decomposition}
\label{sec:cross_algs}
For matrices, if $\mathcal{I}$ and $\mathcal{J}$ are two index sets to denote the ``skeletons" of rows and columns respectively, then one can approximate a matrix $A$ by $A \approx \tilde{A} = CUR$, where $C = A(:,\mathcal{J})$ and $R = A(\mathcal{I},:)$. In this paper, we construct $U$ by $U = C^\dagger AR^\dagger$ since it is easier to do projection enhancement.
For the TT format, we need multiple pivot sets $\mathcal{I}_{\le k}$ and $\mathcal{J}_{>k}$ with $1 \le k \le d-1$. For easier understanding, $\mathcal{I}_{\le k}$ and $\mathcal{J}_{>k}$ can be perceived as the pivot sets for the $k$th unfolding of the target tensor $\mathcal{X}$. Using these pivot sets, the TT cores can be obtained via
\begin{align}
    \mathcal{G}_1 &= X_1(:,\mathcal{J}_{>1}), \ \mathcal{G}_d = Y_{d-1}X_{d-1}(\mathcal{I}_{\le d-1},:), \nonumber \\
    \mathcal{G}_{k} &= {\rm reshape}\left(Y_{k-1}X_{k-1}(\mathcal{I}_{\le k-1},\mathbb{I}_k\otimes\mathcal{J}_{>k}),r,n,r\right), \ 2 \le k \le d-1, \label{eq:ttcross_cores}
\end{align}
where 
\begin{equation} \label{eq:ttcross_PEID_build}
Y_k = X_k(:,\mathcal{J}_{>k})^\dagger X_k X_k(\mathcal{I}_{\le k},:)^\dagger, \quad 1\le k \le d-1,
\end{equation}
$\mathbb{I}_k$ denotes all indices for dimension $k$, and we assume $\mathcal{X}$ has uniform mode size $n$ and all index sets have cardinality $r$. One can easily tell that computing $Y_k$ naively with ``full-complexity" algorithms is very expensive because of the size of $X_k$. In addition, note that the TT approximation assembled through~\cref{eq:ttcross_cores,eq:ttcross_PEID_build} is generally not an orthogonal projection of the original tensor. Connecting cores across successive separations introduces oblique reconstruction operators that extend information from the selected pivot rows, and only the final right hand factor remains an orthogonal row projection.

Nevertheless,~\cref{eq:ttcross_cores,eq:ttcross_PEID_build} is a natural target for error analysis. For the $k$th unfolding, we write
\[
C_k = X_k(:,\mathcal{J}_{>k}),\qquad R_k = X_k(\mathcal{I}_{\le k},:),\qquad 1\le k\le d-1,
\]
and the following result presents a simplified version of an error bound for this TT approximation. The complete theorem and its proof are in~\cref{app:format_analysis}.

\begin{theorem}[Informal version of~\cref{thm:PEID_format_main}]\label{thm:PEID_format_simplified}
For $1\le k\le d-1$, assume that $R_k$ has full row rank, and that the pivot pair $(\mathcal{I}_{\le k},\mathcal{J}_{>k})$ is produced by $r$ steps of ACA with complete pivoting applied to $X_k$, where $\rank(X_k)>r$. Then the TT approximation $\tilde{\mathcal{X}}$ defined by \cref{eq:ttcross_cores,eq:ttcross_PEID_build} satisfies
\[
\left\|\mathcal{X}-\tilde{\mathcal{X}}\right\|_F\;\le\;\sum_{k=1}^{d-1}c_k\,\sigma_{r+1}(X_k),
\]
where $\sigma_{r+1}(X_k)$ is the $(r+1)$st singular value of the $k$th unfolding, and each constant $c_k$ depends on $n$, $d$, $r$, and the conditioning of the pivot rows at the separations preceding $k$.
\end{theorem}

Since $\sigma_{r+1}(X_k)\le\|\mathcal{X}-\mathcal{Y}\|_F$ for every TT tensor $\mathcal{Y}$ with TT ranks at most $r$~\cite{oseledets2011tensor},~\cref{thm:PEID_format_simplified} expresses quasi-optimality with respect to the best TT approximation, up to constants $c_k$, which include amplification factors from row interpolative decomposition and the (pessimistic, but in general unimprovable~\cite{cortinovis2020maximum}) complete-pivoting constant.

Some interpolation-based decomposition algorithms for the TT format choose the pivot sets to be nested:
\begin{align}
    \mathcal{I}_{\le k+1} &\subseteq \mathcal{I}_{\le k} \otimes \mathbb{I}_{k+1}, \nonumber \\ 
    \mathcal{J}_{> k} &\subseteq \mathcal{J}_{> k+1} \otimes \mathbb{I}_{k+1}, \quad 1 \le k \le d-2. \label{eq:tt_nestedness}
\end{align}
This ensures that the original tensor and the approximation match exactly at the interpolated skeletons~\cite{savostyanov2014quasioptimality}. Notably,~\cref{thm:PEID_format_simplified} does not require the pivot sets to be nested, because the format~\cref{eq:ttcross_cores,eq:ttcross_PEID_build} is well defined for arbitrary pivot sets.~\Cref{alg:TT_update_parallel,alg:TT_update_parallel_2sided} require the pivot sets to be nested to guarantee that the cores are glued through sampled rows without touching the full unfoldings (see the recursion~\cref{eq:recursion_app} in~\cref{app:alg3_analysis}), while~\cref{alg:TT_update_sequential,alg:TT_update_sequential_2sided} do not have this restriction. If skeleton pivots are not nested, the same algorithms apply after the selected rows of each unfolding are reconstructed from the previously computed cores. Since our motivation is to improve upon linearly-scaled interpolation-based decompositions such as TTACA, whose greedy pivot selection produces nested sets by construction, we restrict the presentation in this manuscript to the nested case.

Alternatively, one can use $U = A(\mathcal{I},\mathcal{J})^{-1}$ in the matrix CUR approximation. For the TT format, we can adjust~\cref{eq:ttcross_cores,eq:ttcross_PEID_build} accordingly. For completeness,~\cref{app:ttaca_baseline} contains an error bound for this TT approximation, with the assumption that the pivot indices are nested.

\subsection{A matrix prototype for projection enhanced interpolation-based approximations}
\label{sec:CUR_oblique_sub}
\begin{figure}[!htbp]
    \centering
    \includegraphics[width=0.45\linewidth]{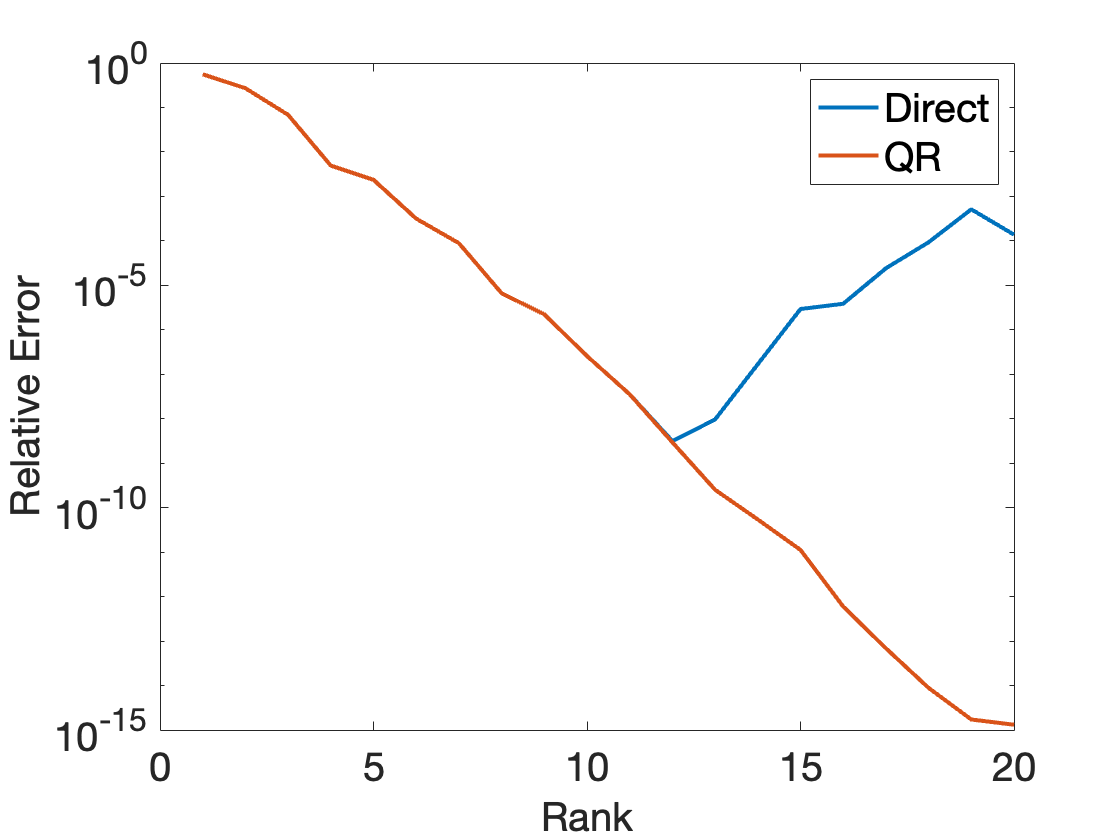}\includegraphics[width=0.45\linewidth]{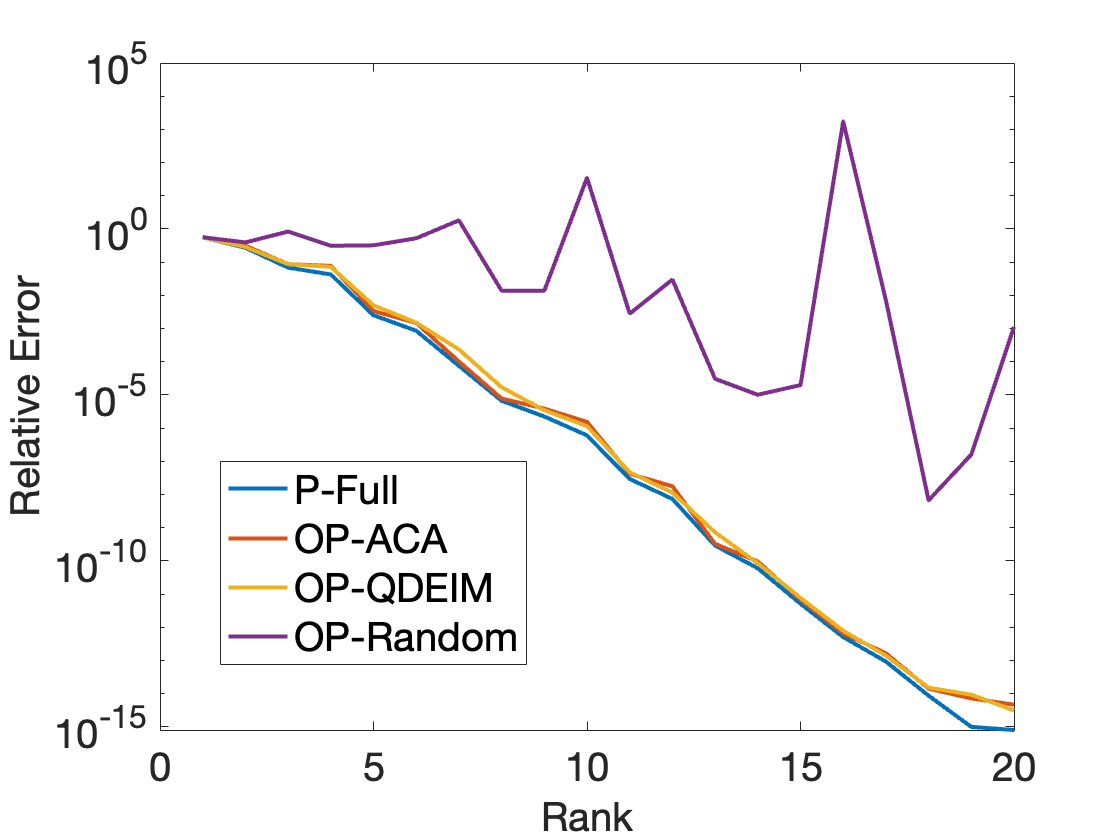}
    \caption{Various ways to build CUR approximations of a $100\times 100$ Hilbert matrix. Left: Demonstrating the stability of using $C$ and $R$ for the construction of $U = C^\dagger AR^\dagger$ (labeled ``Direct"), versus using $Q_C$ and $Q_R$ for $U = Q_C^TAQ_R$ (labeled ``QR"). Right: Skeletonized approximation errors when using varied $U$ construction. Using~\cref{eq:matrix_HIP_ON} is termed ``P-Full",~\cref{eq:matrix_HIP_oblique_sp} is termed ``OP-ACA",~\cref{eq:matrix_HIP_oblique} where $\mathcal{S},\mathcal{T}$ come from applying QDEIM to $Q_C$ and $Q_R$ is termed ``OP-QDEIM", and a random subset of 50 indices for $\mathcal{S},\mathcal{T}$ is termed ``OP-Random".}
    \label{fig:examples24}
\end{figure}

Let $C=A(:,\mathcal{J})$ and $R=A(\mathcal{I},:)$ be the selected columns and rows of $A$. In exact arithmetic,
\[
C\,(C^\dagger A R^\dagger)\,R \;=\; P_C\,A\,P_R,
\qquad P_C:=CC^\dagger,\quad P_R:=R^\dagger R,
\]
is the two-sided orthogonal projection of $A$ onto the selected column and row spaces. The difficulty with $U=C^\dagger AR^\dagger$ is numerical rather than algebraic. When $A$ has low numerical rank relative to the working tolerance, the columns of $C$ and the rows of $R$ are nearly linearly dependent, the pseudo-inverses are formed from ill-conditioned factors, and the computed $U$ depends on roundoff and on the numerical rank thresholds used inside the pseudo-inverse. The blue curve in~\cref{fig:examples24} (left) illustrates this with the Hilbert matrix, where the direct construction stagnates well above machine precision. Perturbations of $A$ by noise have the same effect, since they can turn an exactly dependent set of columns into a nearly dependent one.

Since the only role of $C^\dagger$ and $R^\dagger$ is to represent the column and row spaces of $A$, we can replace them by ON bases. We compute QR factorizations $C=Q_CR_C$ and $R^T=Q_RR_R$ and form
\begin{equation}\label{eq:matrix_HIP_ON}
U = Q_C^TAQ_R,
\end{equation}
which gives the approximation $A\approx Q_CUQ_R^T$. The red curve in~\cref{fig:examples24} (left) shows the performance of this orthogonal construction. Indeed, without truncation, this is the same projection as $P_CAP_R$, so the two curves represent the same approximation in exact arithmetic, and their separation is a pure finite-precision effect. In addition, with rank truncation in the QR step, the construction projects onto the retained subspaces, leading to insensitivity to near dependence and noise. We refer to this ON formulation as the projection enhancement. For the TT format, we substitute $X_k(:,\mathcal{J}_{>k})$ and $X_k(\mathcal{I}_{\le k},:)$ in~\cref{eq:ttcross_cores,eq:ttcross_PEID_build} by such bases. Due to the two challenges posed in~\cref{sec:introduction}, efficient and accurate schemes are needed for the TT format. In this section, we
briefly summarize two methodologies for matrix CUR and PEID. 

In~\cite{engquist2007fast}, the authors propose to select submatrices of $Q_C$ and $Q_R$ to perform projections. In other words, to avoid using the full matrix $A$ in building $U = Q_C^TAQ_R$, one can select index sets $\mathcal{S}$ and $\mathcal{T}$ with some sampling routines,
and calculate $U$ by
\begin{equation} \label{eq:matrix_HIP_oblique}
U = Q_C(\mathcal{S},:)^\dagger A(\mathcal{S},\mathcal{T})Q_R(:,\mathcal{T})^\dagger. 
\end{equation} 
With certain matrices and sampling strategies~\cite{goreinov1997theory,chiu2013sublinear,cortinovis2024sublinear}, this method is shown to have good approximation with high probability, but one can adopt this method as heuristics to general matrices. With $\mathcal{I}$ and $\mathcal{J}$ selected in the skeletonized approximation, they are ideal candidates for performing oblique projections. Therefore, we choose to use $\mathcal{S} = \mathcal{I}$ and $\mathcal{T} = \mathcal{J}$, and then in the matrix case,~\cref{eq:matrix_HIP_oblique} becomes

\begin{equation} \label{eq:matrix_HIP_oblique_sp}
U = Q_C(\mathcal{I},:)^\dagger A(\mathcal{I},\mathcal{J})Q_R(:,\mathcal{J})^\dagger. 
\end{equation}
This allows us to find $U$ efficiently and effectively,
in particular when extending the algorithm to TT decomposition of high dimensional tensors. For example, for the Hilbert matrix,~\cref{fig:examples24} (Right) compares the approximation accuracy of building $U$ with full data~\cref{eq:matrix_HIP_ON}, partial data in~\cref{eq:matrix_HIP_oblique_sp} using ACA pivots, partial data in~\cref{eq:matrix_HIP_oblique} with $\mathcal{S}$ and $\mathcal{T}$ as QDEIM pivots \cite{drmac2016new}, and $\mathcal{S}$ and $\mathcal{T}$ as uniformly random index selections. We can see that using well-selected submatrices from prior knowledge for oblique projections provide reasonably accurate approximations.

Another recent advancement~\cite{cortinovis2024sublinear,park2025accuracy} avoids missing important information in unselected data by oversampling. In short, the skeletonized approximation selects interpolation matrices
\begin{equation} \label{eq:matrix_cur_over}
    C = A(:,\mathcal{J}\cup\mathcal{L}), \quad R = A(\mathcal{I}\cup\mathcal{K},:),
\end{equation}
where $\mathcal{K}$ and $\mathcal{L}$ are index sets with small cardinality for the rows and columns of $A$ respectively, and $\mathcal{K}\cap\mathcal{I}=\emptyset$, $\mathcal{L}\cap\mathcal{J}=\emptyset$. There are various ways to oversample, with uniform oversampling the simplest and cheapest to implement yet still effective.

\cref{fig:Matrix_HIP_Oversample} gives a schematic of the approximation of a given matrix, using an oversampled, normalized, and oblique projection construction for all factors. In the left side of this figure, we assume that we have already selected candidates for $\mathcal{I}$ (purple) and $\mathcal{J}$ (pink). Then, we make an oversampling selection of additional $\mathcal{K}$ (blue) and $\mathcal{L}$ (orange). Following this, we apply a QR on selected columns and rows to obtain $Q_C, Q_R$.
Finally we subselect index sets $\mathcal{S},\mathcal{T}$, e.g. $\mathcal{I}\cup\mathcal{K}$ and $\mathcal{J}\cup\mathcal{L}$ from ACA, to construct the internal component~\cref{eq:matrix_HIP_oblique} via oblique projections.

As a final remark, we emphasize that both techniques, oversampling and selecting submatrices for oblique projections, apply to matrix ID as well, where only column or row space of the matrix is interpolated and projected. This is essential to develop one-sided TT algorithms in the next section.

\begin{figure}[!htbp]
    \begin{center}
        \begin{tikzpicture}[scale=0.75]
            \draw[thick] (0+0.5,0) rectangle (3+0.5,2); 
            \draw[thick,fill,red,opacity=0.2] (0.1+0.5,0) rectangle (0.2+0.5,2);
            \draw[thick,fill,red,opacity=0.2] (0.4+0.5,0) rectangle (0.5+0.5,2);
            \draw[thick,fill,red,opacity=0.2] (1.3+0.5,0) rectangle (1.4+0.5,2);
            \draw[thick,fill,red,opacity=0.2] (2.2+0.5,0) rectangle (2.3+0.5,2);
            \draw[thick,fill,blue,opacity=0.2] (0+0.5,0.2) rectangle (3+0.5,0.3);
            \draw[thick,fill,blue,opacity=0.2] (0+0.5,0.4) rectangle (3+0.5,0.5);
            \draw[thick,fill,blue,opacity=0.2] (0+0.5,1.3) rectangle (3+0.5,1.4);
            \draw[thick,fill,blue,opacity=0.2] (0+0.5,1.7) rectangle (3+0.5,1.8);
            \draw[->] (3.5+0.5,1)--(4.5,2);
            \draw[->] (3.5+0.5,1)--(4.5,0);
            \draw[thick] (0+5,0+1.5) rectangle (3+5,2+1.5); 
            \draw[thick,fill,red,opacity=0.2] (0.1+5,0+1.5) rectangle (0.2+5,2+1.5);
            \draw[thick,fill,red,opacity=0.2] (0.4+5,0+1.5) rectangle (0.5+5,2+1.5);
            \draw[thick,fill,red,opacity=0.2] (1.3+5,0+1.5) rectangle (1.4+5,2+1.5);
            \draw[thick,fill,red,opacity=0.2] (2.2+5,0+1.5) rectangle (2.3+5,2+1.5);
            \draw[thick,fill,orange,opacity=0.75] (0.7+5,0+1.5) rectangle (0.8+5,2+1.5);
            \draw[thick,fill,orange,opacity=0.75] (1.7+5,0+1.5) rectangle (1.8+5,2+1.5);
            \draw[thick,fill,orange,opacity=0.75] (2.7+5,0+1.5) rectangle (2.8+5,2+1.5);
            \draw[thick] (0+5,0-1.5) rectangle (3+5,2-1.5); 
            \draw[thick,fill,blue,opacity=0.2] (0+5,0.2-1.5) rectangle (3+5,0.3-1.5);
            \draw[thick,fill,blue,opacity=0.2] (0+5,0.4-1.5) rectangle (3+5,0.5-1.5);
            \draw[thick,fill,blue,opacity=0.2] (0+5,1.3-1.5) rectangle (3+5,1.4-1.5);
            \draw[thick,fill,blue,opacity=0.2] (0+5,1.7-1.5) rectangle (3+5,1.8-1.5);
            \draw[thick,fill,cyan,opacity=0.75] (0+5,0.6-1.5) rectangle (3+5,0.7-1.5);
            \draw[thick,fill,cyan,opacity=0.75] (0+5,1.0-1.5) rectangle (3+5,1.1-1.5);
            \draw[thick,fill,cyan,opacity=0.75] (0+5,1.5-1.5) rectangle (3+5,1.6-1.5);
            \node[] at (6.5,1) {Oversample};
            \draw[->] (8.5,2.5) -- (9,2.5);
            \draw[thick,fill,red,opacity=0.2] (9.5,1.5) rectangle (9.7,3.5);
            \draw[thick,fill,orange,opacity=0.5] (9.7,1.5) rectangle (9.8,3.5);
            \draw[thick,fill,red,opacity=0.2] (9.8,1.5) rectangle (9.9,3.5);
            \draw[thick,fill,orange,opacity=0.5] (9.9,1.5) rectangle (10,3.5);
            \draw[thick,fill,red,opacity=0.2] (10,1.5) rectangle (10.1,3.5);
            \draw[thick,fill,orange,opacity=0.5] (10.1,1.5) rectangle (10.2,3.5);
            \node[] at (9.8,4) {$A(:,\mathcal{J}\cup\mathcal{L})$};

            \draw[->] (8.5,-0.5)--(9,-0.5);
            \draw[thick,fill,blue,opacity=0.2] (9.5,-1.5) rectangle (9.6,0.5);
            \draw[thick,fill,cyan,opacity=0.5] (9.6,-1.5) rectangle (9.7,0.5);
            \draw[thick,fill,blue,opacity=0.2] (9.7,-1.5) rectangle (9.8,0.5);
            \draw[thick,fill,cyan,opacity=0.5] (9.8,-1.5) rectangle (10.0,0.5);
            \draw[thick,fill,blue,opacity=0.2] (10,-1.5) rectangle (10.2,0.5);
            \node[] at (9.8,-2) {$A(\mathcal{I}\cup\mathcal{K},:)^T$};

            \draw[->] (10.5,2.5) -- (11,2.5); 
            \draw[thick] (11.5,1.5) rectangle (12.1,3.5);
            \node[] at (11.8,4) {$Q_C$};
            
            \draw[thick,fill,blue,opacity=0.2] (11.5,0.2+1.5) rectangle (12.1,0.3+1.5);
            \draw[thick,fill,blue,opacity=0.2] (11.5,0.4+1.5) rectangle (12.1,0.5+1.5);
            \draw[thick,fill,blue,opacity=0.2] (11.5,1.3+1.5) rectangle (12.1,1.4+1.5);
            \draw[thick,fill,blue,opacity=0.2] (11.5,1.7+1.5) rectangle (12.1,1.8+1.5);
            \draw[thick,fill,cyan,opacity=0.75] (11.5,0.6+1.5) rectangle (12.1,0.7+1.5);
            \draw[thick,fill,cyan,opacity=0.75] (11.5,1.0+1.5) rectangle (12.1,1.1+1.5);
            \draw[thick,fill,cyan,opacity=0.75](11.5,1.5+1.5) rectangle (12.1,1.6+1.5);
            \draw[->] (12.5,2.5)--(13,2);

            \node[] at (11.8,1) {Subselect};
            
            \draw[->] (10.5,-0.5)--(11,-0.5);
            \draw[thick] (11.5,-1.5) rectangle (12.1,0.5);
            \node[] at (11.8,-2) {$Q_R$};
            \draw[thick,fill,red,opacity=0.2]     (11.5,1.7-1.5) rectangle (12.1,1.8-1.5);
            \draw[thick,fill,red,opacity=0.2]     (11.5,0.3-1.5) rectangle (12.1,0.4-1.5);
            \draw[thick,fill,red,opacity=0.2]     (11.5,1.4-1.5) rectangle (12.1,1.5-1.5);
            \draw[thick,fill,red,opacity=0.2]     (11.5,0.9-1.5) rectangle (12.1,1.0-1.5);
            \draw[thick,fill,orange,opacity=0.75] (11.5,0.6-1.5) rectangle (12.1,0.7-1.5);
            \draw[thick,fill,orange,opacity=0.75] (11.5,1.1-1.5) rectangle (12.1,1.2-1.5);
            \draw[thick,fill,orange,opacity=0.75] (11.5,0.1-1.5) rectangle (12.1,0.2-1.5);
            \draw[->] (12.5,-2.5+2)--(13,-2+2);

            \node[] at (14.2,1.2) {\footnotesize{$\dagger$}};

            \draw[thick,fill,teal,opacity=0.5] (0+13.5,0.6) rectangle (13.5+0.6,1); 
            \draw[thick,fill,teal,opacity=0.5] (14.3,0.6) rectangle (14.7,0.7);
            \draw[thick,fill,teal,opacity=0.5] (14.3,0.7) rectangle (14.7,0.8);
            \draw[thick,fill,teal,opacity=0.5] (14.3,0.8) rectangle (14.7,0.9);
            \draw[thick,fill,teal,opacity=0.5] (14.3,0.9) rectangle (14.7,1);
            \draw[thick,fill,red,opacity=0.25] (14.3,0.6) rectangle (14.4,1); 
            \draw[thick,fill,red,opacity=0.25] (14.4,0.6) rectangle (14.5,1); 
            \draw[thick,fill,red,opacity=0.25] (14.5,0.6) rectangle (14.6,1); 
            \draw[thick,fill,red,opacity=0.25] (14.6,0.6) rectangle (14.7,1); 
            \draw[thick,fill,red,opacity=0.5] (14.9,0.4) rectangle (15.3,1);
            \node[] at (15.4,1.2) {\footnotesize{$\dagger$}};
            \node[] at (14.5,0) {$U$};
        \end{tikzpicture}
    \end{center}
    \caption{Diagram of the working pieces behind the matrix level approximation using oversampled columns and rows, and submatrix selection for oblique projections.}
    \label{fig:Matrix_HIP_Oversample}
\end{figure}

\section{Projection enhanced interpolation-based decomposition for tensor-train format}
\label{sec:algorithm}
In this section, we target to efficiently conduct PEID for the TT format. Given an existing skeletonized TT approximation, we propose algorithms that enhance the approximation accuracy while maintaining low costs of interpolation-based methods. We use the same tricks of repeatedly oversampling unselected indices~\cite{cortinovis2024sublinear,park2025accuracy} and conducting oblique projections with submatrices~\cite{engquist2007fast} for different tensor dimensions and unfoldings. For notational simplicity, throughout the rest of this section, we assume that $\mathcal{X}$ has uniform mode size $n$, and the number of pivots $r$, is uniform across all dimensions. We also assume that the pivot sets $(\mathcal{I}_{\le k},\mathcal{J}_{>k})$ are known and are nested in the sense of~\cref{eq:tt_nestedness} from an existing skeletonized TT decomposition. In order to conduct fair comparisons with the original skeletonized approximation in~\cref{sec:numerical_tests}, the resulting TT formats from our algorithms still have universal TT core size $r$, but it is straightforward to use other numbers.

\subsection{Dimension parallel approach}
\label{sec:algorithm_parallel}
Since the pivot sets $(\mathcal{I}_{\le k},\mathcal{J}_{>k})$ can be treated independently for different dimension $1 \le k \le d-1$ of a $d$-dimensional tensor $\mathcal{X}$, we start with an ``embarrassingly" dimension parallel implementation (see~\cref{alg:TT_update_parallel}). 

\begin{algorithm}
\caption{TT-PEID-Par: TT construction of a given tensor with prior information on existing skeletonized TT pivots with dimension parallelism.}
\begin{algorithmic}[1]
\label{alg:TT_update_parallel}
\Require {Tensor $\mathcal{X}$, pivot sets $(\mathcal{I}_{\le j},\mathcal{J}_{>j})$, and oversampling parameter $p$.}
\Ensure {The TT format of $\mathcal{X}$ in the form of TT cores $\mathcal{T}_1,\dots,\mathcal{T}_d$.}
\State Set an empty index set $\mathcal{K}_{\le 0}$.
\For {$1 \le j \le d-1$}
\State Construct $\mathcal{K}_{\le j}$ of size $p$ sampled uniformly from $\mathcal{K}_{\le j-1} \otimes \mathbb{I}_{j} \ \backslash \ \mathcal{I}_{\le j}$.
\EndFor
\For {$1 \le j \le d-1$}
\State Construct $\mathcal{L}_{>j}$ of size $p$ sampled uniformly from $\mathbb{I}_{j+1} \otimes \cdots \otimes \mathbb{I}_d \ \backslash \ \mathcal{J}_{> j}$.
\State Find $U_j$ with $r$ ON columns from $X_j\left((\mathcal{I}_{\le j-1} \cup \mathcal{K}_{\le j-1}) \otimes \mathbb{I}_j, \mathcal{J}_{>j} \cup \mathcal{L}_{>j}\right)$.
\State Set $\tilde{U}_j = {\rm reshape}(U_j, r+p, nr)$ for $j > 1$.
\EndFor
\State Set $\mathcal{T}_1 = U_1$.
\For {$2 \le j \le d-1$}
\State Construct $\mathcal{T}_j = {\rm reshape}\left(\left[U_{j-1}(\mathcal{I}_{\le j-1}\cup\mathcal{K}_{\le j-1},:)\right]^\dagger \tilde{U}_j, r, n, r\right)$.
\EndFor
\State Construct $\mathcal{T}_d = \left[U_{d-1}(\mathcal{I}_{\le d-1}\cup\mathcal{K}_{\le d-1},:)\right]^\dagger X_{d-1}(\mathcal{I}_{\le d-1}\cup\mathcal{K}_{\le d-1},:)$.
\end{algorithmic}
\end{algorithm}

\begin{figure}[!htbp]
    \begin{center}
        \begin{tikzpicture}[scale=0.65]
            \draw[thick] (0,0) rectangle (6,2);
            \draw[thick,fill,red,opacity=0.5] (0,0.2) rectangle (6,0.3);
            \draw[thick,fill,red,opacity=0.5] (0,0.7) rectangle (6,0.8);
            \draw[thick,fill,red,opacity=0.5] (0,1.4) rectangle (6,1.5);
            \node[] at (3,-0.5) {$X_1$, red rows: $\mathcal{K}_{\leq 1}$};
            \draw[->] (6.5,1)--(7,1);
            \draw[thick] (7.5,-2) rectangle (11.5,2); 
            \draw[thick,fill,red,opacity = 0.5] (7.5,-1.5) rectangle (11.5,-1);
            \draw[thick,fill,red,opacity = 0.5] (7.5,-0.5) rectangle (11.5,0.0);
            \draw[thick,fill,red,opacity = 0.5] (7.5,0.8) rectangle (11.5,1.3);
            \draw[thick,fill,blue,opacity = 0.5] (7.5,-1.3) rectangle (11.5,-1.2);
            \draw[thick,fill,blue,opacity = 0.5] (7.5,-0.4) rectangle (11.5,-0.3);
            \draw[thick,fill,blue,opacity = 0.5] (7.5,-0.2) rectangle (11.5,-0.1);
            \node[] at (9.5,-2.5) {$X_2$, red region: $\mathcal{K}_{\leq 1}$};
            \node[] at (9.5,-3) {blue rows: $\mathcal{K}_{\leq 2}$};
            \draw[->] (12,1)--(12.5,1);
            \draw[thick] (13,-4) rectangle (15,2);
            \draw[thick,fill,blue,opacity = 0.5] (13,-0.5) rectangle (15,0);
            \draw[thick,fill,blue,opacity = 0.5] (13,-1.5) rectangle (15,-1);
            \draw[thick,fill,blue,opacity = 0.5] (13,-3) rectangle (15,-2.5);
            \draw[thick,fill,yellow,opacity = 0.75] (13,-0.2) rectangle (15,-0.1);
            \draw[thick,fill,yellow,opacity = 0.75] (13,-1.4) rectangle (15,-1.3);
            \draw[thick,fill,yellow,opacity = 0.75] (13,-2.7) rectangle (15,-2.6);
            \node[] at (14,-4.5) {$X_3$, blue regions: $\mathcal{K}_{\leq 2}$};
            \node[] at (14,-5) {yellow rows: $\mathcal{K}_{\leq 3}$};
            \draw[->] (15.5,1)--(16,1);
            \node[] at (16.5,1) {$\dots$};
        \end{tikzpicture}
    \end{center}
    \caption{Visual representation of the nestedness of the index selection for the oversampling sets $\mathcal{K}_{\leq j}$ in~\cref{alg:TT_update_parallel}. $\mathcal{K}_{\leq 1}$ can be freely selected, but the candidate region for $\mathcal{K}_{\leq 2}$ is the region corresponding to all locations in $X_2$ where the first index value is a member of $\mathcal{K}_{\leq 1}$, i.e. blocks of rows in $X_2$. This behavior continues for $\mathcal{K}_{\leq 3}$ and further.}
    \label{fig:TT_update_parallel_sets}
\end{figure}

Overall,~\cref{alg:TT_update_parallel} follows the same logic as~\cref{alg:parallel_TTsketching}. One can easily spot that both algorithms construct low rank TT formats by accomplishing the same goals: (1) line 6 in~\cref{alg:TT_update_parallel} and line 2 in~\cref{alg:parallel_TTsketching} build an ON basis for the column space of (a submatrix of) the $j$th unfolding of $\mathcal{X}$, and (2) line 10 in~\cref{alg:TT_update_parallel} and line 4 in~\cref{alg:parallel_TTsketching} construct TT cores via the column space basis of consecutive unfoldings.

\Cref{alg:TT_update_parallel} can be better understood by considering the oversampled sets $\mathcal{L}_{>j}$ and $\mathcal{K}_{\le j}$ separately. 
\begin{itemize}[leftmargin=*,noitemsep]
    \item The skeletonized TT pivot set $\mathcal{J}_{>j}$ allows us to approximate the column space of the $j$th unfolding of $\mathcal{X}$ with that of $X_j(:,\mathcal{J}_{>j})$. $\mathcal{L}_{>j}$ leads to a better capture using $X_j(:,\mathcal{J}_{>j}\cup\mathcal{L}_{>j})$. This type of oversampling is also used frequently in randomized linear algebra algorithms~\cite{halko2011finding} to ensure high accuracy.

    \item TT core construction in~\cref{alg:parallel_TTsketching} requires an expensive matrix multiplication for high dimensional tensors. For the $j$th TT core, the complexity is $\mathcal{O}(n^{j+1}r^2)$. Instead, the skeletonized TT pivot set $\mathcal{I}_{\le j}$ can be used to select submatrices of $U_j$ and $U_{j+1}$ for oblique projections. To include more uncovered information, we introduce oversampling sets $\mathcal{K}_{\le j}$. As any $U_j$ with $2 \le j \le d$ is used to compute two TT cores, we need to ensure that indices selected by both $\mathcal{I}_{\le j-1}\cup\mathcal{K}_{\le j-1}$ and $\mathcal{I}_{\le j}\cup\mathcal{K}_{\le j}$ can match with those from $U_{j-1}$ and $U_{j+1}$. It is easy to see that the overall nestedness property $\mathcal{I}_{\le j}\cup\mathcal{K}_{\le j} \subseteq (\mathcal{I}_{\le j-1}\cup\mathcal{K}_{\le j-1})\otimes \mathbb{I}_{j}$ is sufficient to meet this requirement. Since we assume pivot sets $\mathcal{I}_{\le j}$ are nested, which can be guaranteed with strategies such as the greedy pivot selection in~\cite{dolgov2020parallel,shi2026distributed}, we need $\mathcal{K}_{\le j}$ to be nested as well, and this is enforced with line 3 of~\cref{alg:TT_update_parallel}. Such a procedure is shown in~\cref{fig:TT_update_parallel_sets}, where the nestedness is enforced by restricting the choices of $\mathcal{K}_{\leq j}$ with regions of $\mathcal{K}_{\leq j-1}$. In a full unfolding, $\mathcal{K}_{\leq j-1}$ produces blocks of rows in $X_j$ as candidates for $\mathcal{K}_{\leq j}$.
\end{itemize}

It is worth noting that one can choose to use any algorithm for line 6 of~\cref{alg:TT_update_parallel}. In particular, for fixed rank algorithms, one of the most efficient option is randomized sketching, while for rank adaptive schemes, rank-revealing QR, or tall-and-skinny QR followed by SVD, are two suitable choices. The matrix accessed in line 6 of~\cref{alg:TT_update_parallel} has $(r+p)n$ rows and $r+p$ columns, so forming it requires $(r+p)^2n$ entry evaluations, and computing $U_j$ by a thin QR factorization, a rank-revealing QR factorization, or a randomized range finder costs $\mathcal{O}\bigl((r+p)^3n\bigr)$ arithmetic operations. Line 10 applies the pseudo-inverse of the $(r+p)\times r$ matrix $U_{j-1}(\mathcal{I}_{\le j-1}\cup \mathcal{K}_{\le j-1},:)$ to the $(r+p)\times nr$ matrix $\tilde U_j$, which costs $\mathcal{O}\bigl((r+p)r^2 + (r+p)r\cdot nr\bigr)=\mathcal{O}\bigl((r+p)r^2n\bigr)$. The total cost of~\cref{alg:TT_update_parallel} is therefore $\mathcal{O}\bigl(dn(r+p)^3\bigr)$ operations and $\mathcal{O}\bigl(dn(r+p)^2\bigr)$ tensor entry evaluations. For fixed $r$ and $p$ this is linear in both $d$ and $n$, and the latter two loops of the algorithm distribute over computing units with $\mathcal{O}\bigl(n(r+p)^3\bigr)$ work each.

The regime $r\ll n$ is an assumption on the data, and under it the linear dependence on $n$ is the dominant feature of the cost. However, when the TT ranks are large, the polynomial dependence on $r$ dominates the runtime, so we keep both the rank $r$ and the oversampling parameter $p$ explicit in the cost estimates. Since PEID is a post-processing step, the relevant practical comparison is TTACA followed by PEID against TTACA run with a tighter tolerance, as both routes spend additional work to reduce the error.~\Cref{sec:numerical_tests} makes this comparison in terms of error versus total runtime.

We conclude this subsection with an a priori error bound for \cref{alg:TT_update_parallel}, where the three algorithmic ingredients---(a) skeleton pivots, (b) oversampling, and (c) submatrix-based oblique projections---contribute separate, interpretable terms. For $1\le k\le d-1$, let $\Delta_k:=\|X_k-(X_k)_r\|_F$ denote the best rank-$r$ approximation error of the $k$th unfolding (the quantity that also governs the TTSVD error~\cite{oseledets2011tensor}), let
\[
A_k := X_k\bigl((\mathcal{I}_{\le k-1}\cup\mathcal{K}_{\le k-1})\otimes\mathbb{I}_k,\,:\bigr)
\]
be the row-restricted unfolding accessed in line 6 of \cref{alg:TT_update_parallel}, and set $\theta_k:=\|A_k(:,\mathcal{J}_{>k}\cup\mathcal{L}_{>k})^\dagger A_k\|_2$. The proof of the following result, together with the intermediate bounds from which it is assembled, is given in~\cref{app:alg3_analysis}.

\begin{theorem}\label{thm:alg3_main}
Let~\cref{ass:incoh_app,ass:incohU_app} of~\cref{app:alg3_analysis} hold (incoherence of the dominant singular subspaces of the unfoldings, with constant $\mu$, and of the computed bases $U_j$, with constant $\mu_U$; the latter is checkable at run time), and let $p\ge c_0\max\{\mu,\mu_U\}\,r\log(r/\delta)$ for an absolute constant $c_0$ and a failure probability parameter $\delta\in(0,1)$. Then, with probability at least $1-3(d-1)\delta$ over the draws of the oversampling sets $\{\mathcal{K}_{\le j},\mathcal{L}_{>j}\}_{j=1}^{d-1}$, the output $\tilde{\mathcal{X}}$ of \cref{alg:TT_update_parallel} satisfies
\begin{equation}\label{eq:alg3_main_bound}
\left\|\mathcal{X}-\tilde{\mathcal{X}}\right\|_F\;\le\;
\sum_{k=1}^{d-1}\Bigl(2\sqrt{\tfrac{(r+p)n}{p}}\Bigr)^{\!k}
\Bigl[\,\min\bigl\{\beta_k,\;(1+c_1\delta^{-1/2})\Delta_k\bigr\}\;+\;\theta_k\Delta_k\Bigr],
\end{equation}
where $c_1\le2\sqrt{2}$ is an absolute constant, and $\beta_k$ is any available bound on the Frobenius error of the skeleton cross approximation restricted to the rows of $A_k$; under the complete-pivoting assumption of \cref{thm:PEID_format_main}, one may take $\beta_k=\sqrt{(r+p)\,n^{d-k+1}}\,4^r\rho_r\,\sigma_{r+1}(X_k)$. In particular, with $\theta:=\max_k\theta_k$,
\[
\left\|\mathcal{X}-\tilde{\mathcal{X}}\right\|_F\;\le\;(d-1)\Bigl(2\sqrt{\tfrac{(r+p)n}{p}}\Bigr)^{\!d-1}\bigl(1+c_1\delta^{-1/2}+\theta\bigr)\,\max_{1\le k\le d-1}\Delta_k .
\]
\end{theorem}

The two arguments of the $\min$ in \cref{eq:alg3_main_bound} quantify the interplay of ingredients (a) and (b): the deterministic branch $\beta_k$ is inherited from the skeleton alone, while the probabilistic branch shows that column oversampling replaces, with high probability, the exponentially pessimistic skeleton constant by an absolute constant times the \emph{optimal} error $\Delta_k$. The prefactor $\bigl(2\sqrt{(r+p)n/p}\bigr)^k$ quantifies ingredient (c): it bounds the amplification incurred by the oblique projections, decreases as the row oversampling $p$ grows, and can be replaced by the computable quantities $\prod_{j\le k}\|[U_j((\mathcal{I}_{\le j}\cup\mathcal{K}_{\le j}),:)^\dagger]\|_2$, which are certifiable a posteriori at negligible cost (see \cref{eq:master_app} and the proof of \cref{thm:alg3_main} in~\cref{app:alg3_analysis}).~\Cref{app:alg3_analysis} also proves two monotonicity statements: enlarging either oversampling set can only decrease the corresponding term in \cref{eq:alg3_main_bound}, on every realization.

\Cref{thm:alg3_main} also gives guidance on the choice of the oversampling parameter. The condition $p\ge c_0\max\{\mu,\mu_U\}\,r\log(r/\delta)$ is sufficient but conservative, and the experiments in~\cref{sec:numerical_tests} use values of $p$ far below this threshold and still observe substantial improvement. Since the incoherence constants are not known in advance, we recommend choosing $p$ adaptively with two inexpensive indicators: (1) the oblique factor $\eta_k=\|[U_k(\mathcal{I}_{\le k}\cup\mathcal{K}_{\le k},:)]^\dagger\|_2$ from~\cref{eq:eta_def}, which costs $\mathcal{O}\bigl((r+p)r^2\bigr)$ to evaluate and enters the a posteriori bound~\cref{eq:alg3_app_bound} directly; (2) the relative error on a fixed set of randomly sampled entries, which is the quantity we report in~\cref{sec:numerical_tests}. Tensor entries evaluated for a smaller $p$ can be reused when $p$ grows.

\subsection{Dimension sequential approach}
\label{sec:algorithm_sequential}
With dimension parallelism, the oversampling index sets need to be nested, which poses a major limitation on the amount of information we can oversample from the unselected parts of the skeletonized approximation. From the discussion of the use of $\mathcal{K}_{\le j}$ in~\cref{sec:algorithm_parallel}, we know these oversampled sets are used mainly to transform orthogonal projections to oblique projections in TT core construction, and they need to be nested since dimension parallel algorithms explore the connection of column space among consecutive tensor unfoldings. However, if the TT cores are generated sequentially, $\mathcal{K}_{\le j}$ can be independent from each other, and in this way, we are able to broaden our scope of oversampling. We describe this sequential procedure to update the TT format in~\cref{alg:TT_update_sequential}.

\begin{algorithm}
\caption{TT-PEID-Seq: TT construction of a given tensor with prior information on existing skeletonized TT pivots.}
\begin{algorithmic}[1]
\label{alg:TT_update_sequential}
\Require {Tensor $\mathcal{X}$, pivot sets $(\mathcal{I}_{\le j},\mathcal{J}_{>j})$, and oversampling parameter $p$.}
\Ensure {The TT format of $\mathcal{X}$ in the form of TT cores $\mathcal{T}_1,\dots,\mathcal{T}_d$.}
\State Construct an index set $\mathcal{L}_{>1}$ of size $p$ sampled uniformly from $\mathbb{I}_{2} \otimes \cdots \otimes \mathbb{I}_d \ \backslash \ \mathcal{J}_{>1}$.
\State Set $Y = X_1\left(:, \mathcal{J}_{>1} \cup \mathcal{L}_{>1}\right)$.
\For {$1 \le j \le d-1$}
\State Find $U_j$ with $r$ ON columns from $Y$.
\State Set $\mathcal{T}_1 = U_1$ or $\mathcal{T}_j = {\rm reshape}(U_j,r,n,r)$ if $j > 1$.
\State Construct $\mathcal{L}_{>j+1}$ of size $p$ sampled uniformly from $\mathbb{I}_{j+2} \otimes \cdots \otimes \mathbb{I}_d  \ \backslash \ \mathcal{J}_{> j+1}$.
\State Construct $\mathcal{K}_{\le j}$ of size $p$ sampled uniformly from $\mathbb{I}_{1} \otimes \cdots \otimes \mathbb{I}_{j} \ \backslash \ \mathcal{I}_{\le j}$.
\State Set an empty matrix $F_{\le j}$ of size $(r+p)$-by-$r$.
\For {$1 \le \ell \le r+p$}
\State Set $m$ as the $\ell$th element of $\mathcal{I}_{\le j} \cup \mathcal{K}_{\le j}$.
\State Find $m_k$ with $1 \le k \le j$ such that $m = \sum_{k=1}^j m_k n^{j-k}$.
\State Construct the $\ell$th row of $F_{\le j}$ by $F_{\le j}(\ell,:) = \mathcal{T}_1(m_1,:)\prod_{k=2}^j \mathcal{T}_k(:,m_k,:)$.
\EndFor
\If {$j < d-1$}
\State Set $Y = {\rm reshape}\left(F_{\le j}^\dagger X_j\left(\mathcal{I}_{\le j}\cup\mathcal{K}_{\le j},\mathbb{I}_{j+1} \otimes (\mathcal{J}_{> j+1}\cup\mathcal{L}_{> j+1})\right),rn,r+p\right)$.
\Else
\State Set $\mathcal{T}_d = F_{\le d-1}^\dagger X_{d-1}\left(\mathcal{I}_{\le d-1}\cup\mathcal{K}_{\le d-1},:\right)$.
\EndIf
\EndFor
\end{algorithmic}
\end{algorithm}

\begin{figure}[!htbp]
    \begin{center}
        \begin{tikzpicture}[scale = 0.89]
            \draw[] (0,0) rectangle (3,2);
            \draw[thick,fill,red,opacity=0.25] (0,0.2) rectangle (3,0.3);
            \draw[thick,fill,red,opacity=0.5] (0,1.2) rectangle (3,1.3);
            \draw[thick,fill,blue,opacity=0.25] (0,0.4) rectangle (3,0.5);
            \draw[thick,fill,blue,opacity=0.5] (0,0.7) rectangle (3,0.8);
            \draw[thick,fill,blue,opacity=0.75] (0,1.6) rectangle (3,1.7);
            \node[] at (1.5,-0.5) {\tiny{$X_j$}};
            \node[] at (1.5,-1) {\tiny{Red: $\mathcal{I}_{\leq j}$. Blue: $\mathcal{K}_{\leq j}$}};
            \node[] at (3.5,1) {$\Rightarrow$};

            \begin{scope}[shift={(0,2.5)}]
            \draw[] (4,0) rectangle (5,2);
            \draw[thick,fill,blue,opacity=0.75] (4,1.6) rectangle (5,1.7); 
            
            \draw[] (5.5,0) rectangle (6.5,2);
            \draw[] (5.5,2)--(6,2.25)--(7,2.25)--(6.5,2);
            \draw[] (7,2.25)--(7,0.25)--(6.5,0);
            \draw[thick,fill,blue,opacity=0.75] (5.75,0.125) rectangle (6.75,2.125);

            \node[] at (7.5,1) {...};
            \draw[] (5.5+2.5,0) rectangle (6.5+2.5,2);
            \draw[] (5.5+2.5,2)--(6+2.5,2.25)--(7+2.5,2.25)--(6.5+2.5,2);
            \draw[] (7+2.5,2.25)--(7+2.5,0.25)--(6.5+2.5,0);
            \draw[thick,fill,blue,opacity=0.75] (5.75+2.5,0.125) rectangle (6.75+2.5,2.125);
            \node[] at (10,1) {$=$};
            \draw[thick,fill,blue,opacity=0.75] (10.5,1) rectangle (11.5,1.1);
            \node[] at (11,0.75) {\tiny{$F_{\leq j}(1,:)$}};
            \node[] at (11,0.4) {\tiny{Size: $1\times r$}};

            \draw[dotted,thick] (3.75,-0.1)--(3.75,2.35)--(11.9,2.35)--(11.9,-0.1);
            \end{scope}

            \begin{scope}[shift={(0,0)}]
            \draw[] (4,0) rectangle (5,2);
            \draw[thick,fill,red,opacity=0.5] (4,1.2) rectangle (5,1.3); 
            
            \draw[] (5.5,0) rectangle (6.5,2);
            \draw[] (5.5,2)--(6,2.25)--(7,2.25)--(6.5,2);
            \draw[] (7,2.25)--(7,0.25)--(6.5,0);
            \draw[thick,fill,red,opacity=0.5] (5.75,0.125) rectangle (6.75,2.125);

            \node[] at (7.5,1) {...};
            \draw[] (5.5+2.5,0) rectangle (6.5+2.5,2);
            \draw[] (5.5+2.5,2)--(6+2.5,2.25)--(7+2.5,2.25)--(6.5+2.5,2);
            \draw[] (7+2.5,2.25)--(7+2.5,0.25)--(6.5+2.5,0);
            \draw[thick,fill,red,opacity=0.5] (5.75+2.5,0.125) rectangle (6.75+2.5,2.125);
            \node[] at (10,1) {$=$};
            \draw[thick,fill,red,opacity=0.5] (10.5,1) rectangle (11.5,1.1);
            \node[] at (11,0.75) {\tiny{$F_{\leq j}(2,:)$}};
            \node[] at (11,0.4) {\tiny{Size: $1\times r$}};
            \draw[dotted,thick] (3.75,-0.1) rectangle (11.9,2.35);
            \end{scope}

            \begin{scope}[shift={(0,-3)}]
            \draw[] (4,0) rectangle (5,2);
            \draw[thick,fill,red,opacity=0.25] (4,0.4) rectangle (5,0.5); 
            
            \draw[] (5.5,0) rectangle (6.5,2);
            \draw[] (5.5,2)--(6,2.25)--(7,2.25)--(6.5,2);
            \draw[] (7,2.25)--(7,0.25)--(6.5,0);
            \draw[thick,fill,red,opacity=0.25] (5.75,0.125) rectangle (6.75,2.125);

            \node[] at (6.5,2.75) {$\vdots$};
            \node[] at (7.5,1) {...};
            \draw[] (5.5+2.5,0) rectangle (6.5+2.5,2);
            \draw[] (5.5+2.5,2)--(6+2.5,2.25)--(7+2.5,2.25)--(6.5+2.5,2);
            \draw[] (7+2.5,2.25)--(7+2.5,0.25)--(6.5+2.5,0);
            \draw[thick,fill,red,opacity=0.25] (5.75+2.5,0.125) rectangle (6.75+2.5,2.125);
            \node[] at (10,1) {$=$};
            \draw[thick,fill,red,opacity=0.25] (10.5,1) rectangle (11.5,1.1);
            \node[] at (11,0.75) {\tiny{$F_{\leq j}(r+p,:)$}};
            \node[] at (11,0.4) {\tiny{Size: $1\times r$}};
            \draw[dotted,thick] (3.75,-0.1) rectangle (11.9,2.35);
            \end{scope}
            \draw[->] (12,3)--(12.5,2);
            \draw[->] (12,1)--(12.5,1);
            \draw[->] (12,-1.5)--(12.5,-0.5);
            \draw[thick,fill,blue,opacity = 0.75] (13,1.2) rectangle (14,1.3);
            \draw[thick,fill,red,opacity = 0.5] (13,1.1) rectangle (14,1.2);
            \draw[thick,fill,blue,opacity = 0.5] (13,1) rectangle (14,1.1);
            \draw[thick,fill,blue,opacity = 0.25] (13,0.9) rectangle (14,1);
            \draw[thick,fill,red,opacity = 0.25] (13,0.8) rectangle (14,0.9);
            \node[] at (13.5,0.5) {\tiny{$F_{\leq j}$}};
            \node[] at (13.5,0.2) {\tiny{Size: $(r+p)\times r$}};
        \end{tikzpicture}
    \end{center}
    \caption{Construction of $F_{\leq j}$ inside of~\cref{alg:TT_update_sequential}. The left image corresponds to the locations $\mathcal{I}_{\leq j}\cup \mathcal{K}_{\leq j}$ in the $j$-th unfolding $X_j$. For each index in $\mathcal{I}_{\leq j}\cup \mathcal{K}_{\leq j}$, we separate it out into $j$ indices. Then, we perform a sequence of vector matrix products, where we use the $j$ indices to take the vector from the first core $\mathcal{T}_1$, and the remaining matrices from the already computed cores $\mathcal{T}_2,\dots,\mathcal{T}_{j}$. This multiplication forms one row of $F_{\leq j}$.
    }
    \label{fig:TT_update_sequential_F_const}
\end{figure}

\Cref{alg:TT_update_sequential} is an extension of~\cref{alg:TTsketching} or other dimension sequential TT factorization algorithms. These algorithms usually constitute of two steps for the $j$th TT core construction: (1) form an ON basis of the column space of a matrix that contains unused data and use it as the $j$th TT core (line 4 of~\cref{alg:TTsketching}), and (2) organize the remaining data into a new matrix to be used for recovering the $(j+1)$st TT core (line 6 of~\cref{alg:TTsketching}). Similar to understanding~\cref{alg:TT_update_parallel} in~\cref{sec:algorithm_parallel}, we discuss the uniformly oversampled sets $\mathcal{K}_{\le j}$ and $\mathcal{L}_{> j+1}$ separately:
\begin{itemize}[leftmargin=*,noitemsep]
    \item One major computational burden of~\cref{alg:TTsketching} is caused by the number of columns of $Y$ in line 4 and 6, which is $\mathcal{O}(n^{d-j})$ when dealing with dimension $j$. Therefore, we select a subset of the columns of $Y$ to use in the computation, where the columns are indexed by the skeletonized and oversampled sets.
    
    \item The oversampled indices $\mathcal{K}_{\le j}$ are used together with $\mathcal{I}_{\le j}$ for subselection to convert orthogonal projections to oblique projections. However, in line 14 of~\cref{alg:TT_update_sequential}, if we use the same idea in~\cref{alg:TTsketching} to update $Y$ with $Y$ from the previous iteration, then the row subsets are already restricted by all previous row selections. In other words, $\mathcal{K}_{\le j}$ are nested just as in~\cref{alg:TT_update_parallel}. Therefore, in order to achieve independent selections of $\mathcal{K}_{\le j}$ for the $j$th dimension for a new $Y$, we have to use the previous TT cores $\mathcal{T}_1,\cdots,\mathcal{T}_j$ to compress the first $j$ dimensions of $\mathcal{X}$. The naive way is to construct a matrix $G_{\le j}$ of size $n^j \times r$ using Kronecker products, project it onto $X_j$, and subselect the rows from $G_{\le j}^\dagger X_j$, but forming $G_{\le j}$ and performing this least squares solve are both expensive due to the number of rows in $G_{\le j}$. Instead, we only form $F_{\le j}$, a subset of rows of $G_{\le j}$, by uniformly sampling row indices in the range of $1$ and $n^j$ not covered in $\mathcal{I}_{\le j}$. This ensures the oversampling to be flexible and independent across dimensions. A visual demonstration of the construction of $F_{\leq j}$ is shown in~\cref{fig:TT_update_sequential_F_const}. Here we see that row locations in the left most portion of the figure are mapped to tuples, which are used to take slices of previously formed cores $\mathcal{T}_1,\dots,\mathcal{T}_j$. These slices are then multiplied as a sequence of vector-matrix products, depicted by the dotted boxes in the center of the figure, to produce a single row of $F_{\leq j}$. After constructing $F_{\leq j}$, we can build the selected rows of $G_{\le j}$ using the previously constructed TT cores as in lines 10-12 in~\cref{alg:TT_update_sequential}. In particular, line 11 can be understood as converting linear indices to subscript tuples of a multidimensional array. If one is familiar with MATLAB notations, then this is the ``\textit{ind2sub}" command.
\end{itemize}

Similar to~\cref{alg:TT_update_parallel}, one is free to choose any algorithm for line 4 of~\cref{alg:TT_update_sequential}. From construction, the size of $Y$ is $n \times (r+p)$ for $j=1$ and $(nr)\times(r+p)$ otherwise, so line 4 normally takes $\mathcal{O}(nr(r+p)^2)$ to run. In addition, for dimension $j$, forming one row of $F_{\le j}$ involves the multiplication of one row vector of size $r$ and $(j-1)$ matrices of size $r\times r$, which results in a complexity of $\mathcal{O}((j-1)r^2)$, and thus building $F_{\le j}$ takes $\mathcal{O}(jr^2(r+p))$ and the $Y$ update in line 14 takes $\mathcal{O}((r+p)^2rn)$. In practice, since $n$ is often a lot larger than $d$, line 14 dominates the runtime complexity of each iteration, so the overall complexity of~\cref{alg:TT_update_sequential} is $\mathcal{O}(dnr(r+p)^2)$. Finally, since $F_{\le j}$ is calculated row by row, we can only use matrix-vector multiplications instead of matrix-matrix multiplications. This indicates that in practice when we implement with basic linear algebra subprograms (BLAS)~\cite{dongarra1990set}, we can only use the slower BLAS2 routines instead of BLAS3. On the contrary,~\cref{alg:TT_update_parallel} is coded with BLAS3, and can run much faster practically, even in a serial computing environment.

A natural question is whether \cref{alg:TT_update_sequential} admits an a priori
error bound analogous to \cref{thm:alg3_main}. The analysis of
\cref{alg:TT_update_parallel} in~\cref{app:alg3_analysis} relies on two
structural properties that the sequential algorithm deliberately gives up in
exchange for its broader oversampling choices. First, nestedness of the row sets
makes $\mathcal{I}_{\le j}\cup\mathcal{K}_{\le j}$ a subset of the rows of the
restricted unfolding accessed at step $j$, so that the computed approximation
satisfies the exact recursion \cref{eq:recursion_app}; once $\mathcal{K}_{\le j}$
is freed, the selected rows must instead be reconstructed through the matrix
$F_{\le j}$ assembled from the previously computed cores (lines 9--12 of
\cref{alg:TT_update_sequential}), and the recursion is replaced by a perturbed
one in which the pseudo-inverse $F_{\le j}^\dagger$ acts on rows that are
themselves only approximately represented by the current cores. Bounding the
resulting error requires controlling products of previous-level errors together
with the conditioning of $F_{\le j}$, i.e., a nonlinear error propagation in
place of the linear recursion behind \cref{eq:master_app}. Second, in \cref{alg:TT_update_parallel} the column draws $\mathcal{L}_{>j}$
are independent across dimensions and each row draw $\mathcal{K}_{\le j}$ is
uniform over a candidate set fixed by the previous level, so the matrices being
sampled at level $j$ are determined by the index sets alone, the failure
probabilities in the proof of \cref{thm:alg3_main} simply add up, and the
incoherence hypotheses can be imposed on the data. In
\cref{alg:TT_update_sequential}, by contrast, the draws at step $j$ also
determine the compressed matrices $F_{\le j}^\dagger X_j(\cdots)$ from which
later steps sample: the sampled data are themselves random, the events at
different levels are coupled, and the incoherence of \cref{ass:incoh_app}
needs to be shown to propagate through these compressed matrices rather
than being assumed of the data. Neither obstacle appears conceptually
fundamental---in particular, deterministic, a posteriori certifiable bounds in
terms of quantities such as $\sigma_{\min}(F_{\le j})$ remain available at every
step---but a sharp probabilistic analysis requires new ingredients, and we leave
it to future work. Meanwhile, the experiments in \cref{sec:numerical_tests} show
that \cref{alg:TT_update_sequential} attains accuracy comparable to, and in some
of the more challenging high-dimensional tests better than,
\cref{alg:TT_update_parallel}, consistent with the intuition that broader
oversampling can only enrich the sampled information.

\subsection{Variants and extensions of proposed algorithms}
\label{sec:algorithm_extensions}
\Cref{alg:TT_update_parallel} and \cref{alg:TT_update_sequential} share the similarity that the TT cores are constructed mainly from ON bases of the column spaces of tensor unfoldings. This is reflected in the construction of the TT cores in the order of $\mathcal{T}_1,\dots,\mathcal{T}_d$. In the matrix regime, this means that we initiate the projection-based post-processing with ID. Therefore, one straightforward way to extend both of our proposed algorithms for accuracy and efficiency improvement is to build TT cores with ON bases of both column and row spaces of tensor unfoldings, which is similar to using CUR instead of ID in the matrix setting.

Since finding ON bases $U_j$ in line 6 of the parallel algorithm~\cref{alg:TT_update_parallel} is independent and simultaneous, we can alternatively perform this step only for $1 \le j \le \floor{d/2}$. For $\floor{d/2}+1 \le j \le d-1$, we construct $V_j$ with ON columns for sampled submatrices of $X_j^T$. This can ``embarrassingly" double the parallelization efficiency. In addition, this also alleviates the aggressiveness of using $(\mathcal{I}_{\le j-1} \cup \mathcal{K}_{\le j-1}) \otimes \mathbb{I}_j$ to select only $(r+p)n$ from a total of $n^j$ rows of $X_j$ for large $j$, and allows us to build basis $U_j$ that can better characterize $X_j$. In this way, the TT cores can be more accurate. We summarize this variant of~\cref{alg:TT_update_parallel} in~\cref{alg:TT_update_parallel_2sided}, and for simplicity, we assume the dimension $d$ to be odd. In summary, in~\cref{alg:TT_update_parallel_2sided}, in order to construct cores from both ends, we need to construct 4 oversampling sets instead of 2, labeled as $\mathcal{K}_{\le j}^L, \mathcal{K}_{> j}^R, \mathcal{K}_{\le j}^L, \mathcal{K}_{> j}^R$, and construct the last middle TT core with both $U_{(d-1)/2}$ and $V_{(d-1)/2}$ (see lines 13-16).

\begin{algorithm}
\caption{TT-PEID-Par2: Two-sided version of TT-PEID-Par~\cref{alg:TT_update_parallel}.}
\begin{algorithmic}[1]
\label{alg:TT_update_parallel_2sided}
\Require {Tensor $\mathcal{X}$, pivot sets $(\mathcal{I}_{\le j},\mathcal{J}_{>j})$, and oversampling parameter $p$.}
\Ensure {The TT format of $\mathcal{X}$ in the form of TT cores $\mathcal{T}_1,\dots,\mathcal{T}_d$.}
\State Set empty index set $\mathcal{K}_{\le 0}^L$ and $\mathcal{K}_{>d}^R$.
\For {$1 \le j \le (d-1)/2$}
\State Construct index sets $\mathcal{K}_{\le j}^L$ and $\mathcal{K}_{>d-j}^R$ of size $p$ sampled uniformly from $\mathcal{K}_{\le j-1}^L \otimes \mathbb{I}_{j} \ \backslash \ \mathcal{I}_{\le j}$ and $\mathcal{K}_{> d-j+1}^R \otimes \mathbb{I}_{d-j} \ \backslash \ \mathcal{J}_{> d-j}$.
\State Construct index sets $\mathcal{L}_{>j}^L$ and $\mathcal{L}_{\le d-j}^R$ of size $p$ sampled uniformly from $\mathbb{I}_{j+1} \otimes \cdots \otimes \mathbb{I}_d \ \backslash \ \mathcal{J}_{> j}$ and $\mathbb{I}_{1} \otimes \cdots \otimes \mathbb{I}_{d-j} \ \backslash \ \mathcal{I}_{\le d-j}$.
\State Construct $Y = X_j\left(\left(\mathcal{I}_{\le j-1} \cup \mathcal{K}_{\le j-1}^L\right) \otimes \mathbb{I}_j, \mathcal{J}_{>j} \cup \mathcal{L}_{>j}^L\right)$.
\State Construct $Z = X_{d-j}^T\left(\mathcal{I}_{\le d-j} \cup \mathcal{L}_{\le d-j}^R, \mathbb{I}_{d-j+1} \otimes \left(\mathcal{J}_{>d-j+1} \cup \mathcal{K}_{>d-j+1}^R\right)\right)$.
\State Find $U_j$ and $V_j$ with $r$ ON columns from $Y$ and $Z$.
\State Set $\tilde{U}_j = {\rm reshape}(U_j, r+p, nr)$ and $\tilde{V}_j = \rm{reshape}(V_j, r+p, nr)$ for $j > 1$.
\EndFor
\State Set $\mathcal{T}_1 = U_1$ and $\mathcal{T}_{d} = V_1^T$.
\For {$2 \le j \le (d-1)/2$}
\State Construct $\mathcal{T}_j = {\rm reshape}\left(\left[U_{j-1}\left(\mathcal{I}_{\le j-1}\cup\mathcal{K}_{\le j-1}^L,:\right)\right]^\dagger \tilde{U}_j, r, n, r\right)$.
\State Construct $\mathcal{T}_{d-j+1} = {\rm reshape}\left(\tilde{V}_j^T\left[V_{j-1}^T\left(:,\mathcal{J}_{>d-j+2}\cup\mathcal{K}_{>d-j+2}^R\right)\right]^\dagger, r, n, r\right)$.
\EndFor
\State Set $\mathcal{T}_{(d+1)/2} = X_{(d-1)/2}(\mathcal{I}_{\le(d-1)/2}\cup \mathcal{K}^L_{\le (d-1)/2},\mathbb{I}_{(d+1)/2}\otimes (\mathcal{J}_{>(d+1)/2}\cup \mathcal{K}^R_{>(d+1)/2}))$
\State Update $\mathcal{T}_{(d+1)/2} = \left[U_{(d-1)/2}\left(\mathcal{I}_{\le (d-1)/2}\cup\mathcal{K}_{\le (d-1)/2}^L,:\right)\right]^\dagger \mathcal{T}_{(d+1)/2}$.
\State Update $\mathcal{T}_{(d+1)/2} = {\rm reshape}\left(\mathcal{T}_{(d+1)/2},rn,r\right)$
\State Update $\mathcal{T}_{(d+1)/2} =\mathcal{T}_{(d+1)/2} \left[V_{(d-1)/2}^T\left(:,\mathcal{J}_{>(d+1)/2}\cup\mathcal{K}_{>(d+1)/2}^R\right)\right]^\dagger$ and $\mathcal{T}_{(d+1)/2} = {\rm reshape}\left(\mathcal{T}_{(d+1)/2},r,n,r\right)$.
\end{algorithmic}
\end{algorithm}

Similarly, we can develop a 2-sided variant of the sequential algorithm~\cref{alg:TT_update_sequential} (see~\cref{alg:TT_update_sequential_2sided} for the pseudo-code when $d$ is odd) by constructing 4 oversampling sets, and treating the middle TT core $\mathcal{T}_{(d+1)/2}$ as the last core to build. In 1-sided sequential TT approximation algorithms, error propagates and accumulates when proceeding to the next dimension, i.e., error for constructing TT core $\mathcal{T}_j$ is carried over to the rest of the cores $\mathcal{T}_k$ with $k > j$, and it can be especially problematic for very high $d$. The 2-sided algorithm partially fixes the problem. However,~\cref{alg:TT_update_sequential_2sided} requires roughly double memory footprints than~\cref{alg:TT_update_sequential}, as constructions from both ends cannot be combined, leading to storage of both $Y$ and $Z$ matrices in the algorithm.

\begin{algorithm}
\caption{TT-PEID-Seq2: Two-sided version of TT-PEID-Seq~\cref{alg:TT_update_sequential}.}
\begin{algorithmic}[1]
\label{alg:TT_update_sequential_2sided}
\Require {Tensor $\mathcal{X}$, pivot sets $(\mathcal{I}_{\le j},\mathcal{J}_{>j})$, and oversampling parameter $p$.}
\Ensure {The TT format of $\mathcal{X}$ in the form of TT cores $\mathcal{T}_1,\dots,\mathcal{T}_d$.}
\State Construct index sets $\mathcal{L}_{>1}^L$ and $\mathcal{L}_{\le d-1}^R$ of size $p$ sampled uniformly from $\mathbb{I}_{2} \otimes \cdots \otimes \mathbb{I}_d \ \backslash \ \mathcal{J}_{>1}$ and $\mathbb{I}_{1} \otimes \cdots \otimes \mathbb{I}_{d-1} \ \backslash \ \mathcal{I}_{\le d-1}$.
\State Set $Y = X_1\left(:, \mathcal{J}_{>1} \cup \mathcal{L}_{>1}^L\right)$ and $Z = X_{d-1}^T\left( \mathcal{I}_{\le d-1} \cup \mathcal{L}_{\le d-1}^R,:\right)$.
\For {$1 \le j \le (d-1)/2$}
\State Find $U_j$ and $V_j$ with $r$ ON columns from $Y$ and $Z$.
\State Set $\mathcal{T}_1 = U_1$ or $\mathcal{T}_j = {\rm reshape}(U_j,r,n,r)$ if $j > 1$.
\State Set $\mathcal{T}_d = V_1^T$ or $\mathcal{T}_{d-j+1} = {\rm reshape}(V_j^T,r,n,r)$ if $j > 1$.
\State Construct index sets $\mathcal{L}_{>j+1}^L$ and $\mathcal{L}_{\le d-j-1}^R$ of size $p$ sampled uniformly from $\mathbb{I}_{j+2} \otimes \cdots \otimes \mathbb{I}_d  \ \backslash \ \mathcal{J}_{> j+1}$ and $\mathbb{I}_{1} \otimes \cdots \otimes \mathbb{I}_{d-j-1}  \ \backslash \ \mathcal{I}_{\le d-j-1}$.
\State Construct index sets $\mathcal{K}_{\le j}^L$ and $\mathcal{K}_{>d-j}^R$ of size $p$ sampled uniformly from $\mathbb{I}_{1} \otimes \cdots \otimes \mathbb{I}_{j} \ \backslash \ \mathcal{I}_{\le j}$ and $\mathbb{I}_{d-j+1} \otimes \cdots \otimes \mathbb{I}_{d} \ \backslash \ \mathcal{J}_{> d-j}$.
\State Set empty matrices $F_{\le j}$ and $G_{>d-j}$ of size $(r+p)$-by-$r$.
\For {$1 \le \ell \le r+p$}
\State Set $m$ as the $\ell$th element of $\mathcal{I}_{\le j} \cup \mathcal{K}_{\le j}^L$.
\State Find $m_k$ with $1 \le k \le j$ such that $m = \sum_{k=1}^j m_k n^{j-k}$.
\State Construct the $\ell$th row by $F_{\le j}(\ell,:) = \mathcal{T}_1(m_1,:)\prod_{k=2}^j \mathcal{T}_k(:,m_k,:)$.
\State Set $p$ as the $\ell$th element of $\mathcal{J}_{> d-j} \cup \mathcal{K}_{>d-j}^R$.
\State Find $p_k$ with $1 \le k \le j$ such that $p = \sum_{k=1}^j p_k n^{j-k}$.
\State Construct the $\ell$th row by $G_{>d-j}(\ell,:) = \mathcal{T}_d^T(:,p_1)\prod_{k=2}^j \mathcal{T}_{d+1-k}^T(:,p_k,:)$.
\EndFor
\If {$j < (d-1)/2$}
\State Set $Y = F_{\le j}^\dagger X_j\left(\mathcal{I}_{\le j}\cup\mathcal{K}_{\le j}^L,\mathbb{I}_{j+1} \otimes \left(\mathcal{J}_{> j+1}\cup\mathcal{L}_{> j+1}^L\right)\right)$ and $Y = {\rm reshape}\left(Y, rn, r+p \right)$.
\State Set $Z = G_{>d-j}^\dagger X_{d-j-1}^T\left(\mathcal{I}_{\le d-j-1}\cup\mathcal{L}_{\le d-j-1}^R,\mathbb{I}_{d-j} \otimes \left(\mathcal{J}_{> d-j}\cup\mathcal{K}_{> d-j}^R\right)\right)$ and $Z = {\rm reshape}\left(Z, rn, r+p \right)$.
\Else
\State Collect right index set $\mathcal{J}_R :=\mathbb{I}_{(d+1)/2}\otimes \left(\mathcal{J}_{>(d+1)/2}\cup \mathcal{K}^R_{>(d+1)/2}\right)$.
\State Set $\mathcal{T}_{(d+1)/2} = X_{(d-1)/2}\left(\mathcal{I}_{\le (d-1)/2}\cup\mathcal{K}_{\le (d-1)/2}^L,\mathcal{J}_R\right)$.
\State Update $\mathcal{T}_{(d+1)/2} = {\rm reshape}\left(F_{\le (d-1)/2}^\dagger \mathcal{T}_{(d+1)/2},rn,r\right)$.
\State Set $\mathcal{T}_{(d+1)/2} = \mathcal{T}_{(d+1)/2}\left(:,\mathcal{J}_{> (d+1)/2}\cup\mathcal{K}_{> (d+1)/2}^R\right)\left(G_{>(d+1)/2}^T\right)^\dagger$ and $\mathcal{T}_{(d+1)/2} = {\rm reshape}\left(\mathcal{T}_{(d+1)/2}, r, n, r\right)$.
\EndIf
\EndFor
\end{algorithmic}
\end{algorithm}

If computational resources are sufficient, we can perform multiple times of 1-sided algorithms for a better approximation. Since the order of generating the TT cores matter, which is reflected in whether to use column spaces or row spaces of the tensor unfoldings, the simplest way to improve accuracy is to construct two sets of TT cores: $\mathcal{R}_1, \dots,\mathcal{R}_d$ from dimension 1 to $d$ using column spaces of unfoldings, and $\mathcal{S}_1, \dots,\mathcal{S}_d$ from dimension $d$ to 1 using row spaces of unfoldings. Both sets of cores provide approximations to the target tensor, so we can concatenate the cores together to form a new set of TT cores $\mathcal{T}_1 \in \R^{n \times 2r}$, $\mathcal{T}_d \in \R^{2r \times n}$, and $\mathcal{T}_j \in \R^{2r \times n \times 2r}$ for $2 \le j \le d-1$. We can represent these cores mathematically as:
\begin{align}
\mathcal{T}_1 &= \begin{bmatrix} \mathcal{R}_1 & \mathcal{S}_d^T \end{bmatrix}, \quad \mathcal{T}_d = \begin{bmatrix} \mathcal{R}_d \\ \mathcal{S}_1^T \end{bmatrix}, \nonumber \\
\mathcal{T}_j(1:r,i_j,1:r) &= \mathcal{R}_j(:,i_j,:), \ \mathcal{T}_j(r+1:2r,i_j,r+1:2r) = \mathcal{S}_j^T(:,i_j,:), \ 1 \le i_j \le n. \label{eq:concatenate_tensor}
\end{align}
Finally, we conduct TT-rounding for a better and storage-efficient approximation. In practice, for a target tensor $\mathcal{X}$, we can denote a new tensor $\mathcal{Y}$ to be the reverse of $\mathcal{X}$, i.e.
\begin{equation} \label{eq:reverse_tensor}
    \mathcal{X}_{i_1,\cdots,i_d} = \mathcal{Y}_{i_d,\cdots,i_1}, \quad 1 \le i_j \le n, \quad 1 \le j \le d,
\end{equation}
then the pivot sets of $\mathcal{Y}$ are $(\mathcal{J}_{>d-j},\mathcal{I}_{\le d-j})$, and we can find $\mathcal{S}_1, \dots,\mathcal{S}_d$ by running 1-sided algorithms on $\mathcal{Y}$. We summarize this procedure in~\cref{alg:TT_update_rounding}.

\begin{algorithm}
\caption{TT-PEID-R: Rounding on results of two runs of one-sided~\cref{alg:TT_update_parallel} or~\cref{alg:TT_update_sequential} for better accuracy.}
\begin{algorithmic}[1]
\label{alg:TT_update_rounding}
\Require {Tensor $\mathcal{X}$, pivot sets $(\mathcal{I}_{\le j},\mathcal{J}_{>j})$, and oversampling parameter $p$.}
\Ensure {The TT format of $\mathcal{X}$ in the form of TT cores $\mathcal{T}_1,\dots,\mathcal{T}_d$.}
\State Use~\cref{alg:TT_update_parallel} or~\cref{alg:TT_update_sequential} on $\mathcal{X}$ and $(\mathcal{I}_{\le j},\mathcal{J}_{>j})$ to get TT cores $\mathcal{R}_1, \dots,\mathcal{R}_d$.
\State Use~\cref{alg:TT_update_parallel} or~\cref{alg:TT_update_sequential} on $\mathcal{Y}$ from~\cref{eq:reverse_tensor} and $(\mathcal{J}_{>d-j},\mathcal{I}_{\le d-j})$ to get TT cores $\mathcal{S}_1, \dots,\mathcal{S}_d$.
\State Set new TT cores $\mathcal{T}_1,\dots,\mathcal{T}_d$ as~\cref{eq:concatenate_tensor} and perform TT rounding.
\end{algorithmic}
\end{algorithm}

\section{Numerical Tests}
\label{sec:numerical_tests}
In this section we show the performance of the proposed~\cref{alg:TT_update_parallel,alg:TT_update_sequential,alg:TT_update_parallel_2sided,alg:TT_update_sequential_2sided,alg:TT_update_rounding} attached to a standard TTACA implementation~\cite{shi2026distributed} on four different tensors\footnote{The codes for the algorithm implementation can be found at \href{https://github.com/dhayes95/PEID}{https://github.com/dhayes95/PEID}.}, (1) the Hilbert tensor, (2) two tensors corresponding to the evaluations of Mat\'ern and Thin Plate Spline kernel functions, and (3) the Maxwellian tensor encountered in kinetic simulations of differential equations~\cite{einkemmer2025review}. We start by introducing these tensors in~\cref{sec:TensorDef}. Then, in~\cref{sec:TTSVD_Comparisons}, we perform a comparison of accuracy for a given target rank of our proposed algorithms, with standard TTSVD from Tensorly~\cite{tensorly}, MATLAB DMRG~\cite{OseledetsTTToolbox}, TnTorch's cross~\cite{tntorch}, and TT-Sketching~\cite{shi2023parallel}. As different implementations use multiple programming languages and have limitations on tensor size, the baseline comparison is performed on a local machine running a 12th Gen Intel Core i9-12900HX CPU (16 physical cores and 24 logical processors). Next, we show the accuracy and timing with respect to oversample sizes for various dimensions of the PEID algorithms in~\cref{sec:PEID_Performance}. Lastly, we demonstrate the robustness of the proposed algorithms with noisy input tensors in~\cref{sec:Noisy_test}. For the latter two groups of tests, results are computed on the CPU nodes of Perlmutter, a Cray EX supercomputer hosted by NERSC at Lawrence Berkeley National Laboratory with 2 AMD EPYC 7763 CPUs per node (i.e. $2\times64=128$ cores per node).

For error measurements, we compute an approximate relative error by sampling a set $R$ of indices uniformly in the tensor, and computing the error at these locations, i.e.
$$
\frac{||\mathcal{X} - \tilde{\mathcal{X}}||_{F,R}}{||\mathcal{X}||_{F,R}} = \left(\frac{\sum_{i \in R}\left[\mathcal{X}(i) - \tilde{\mathcal{X}}(i)\right]^2}{\sum_{i \in R} \left[\mathcal{X}(i)\right]^2}\right)^{1/2}.
$$
Here, $\|\cdot\|_{F,R}$ denotes the Frobenius norm restricted to the sampled entries, and the ratio is a sampled estimate of the relative error. For a fair comparison, we set a seed to fix the indices used for accuracy computations. In all tests, $|R| = 100,000$.

In addition, we report the reduction factor for each test. This is a measurement on the reduction of the sampled error, between our proposed ACA initiated PEID, and the original ACA approximation. The reduction factor is defined as 
$$
\text{RF} = \frac{||\mathcal{X} - \tilde{\mathcal{X}}_{\rm{ACA}}||_{F,R}}{||\mathcal{X} - \tilde{\mathcal{X}}_{\rm{PEID}}||_{F,R}},
$$
where $\tilde{\mathcal{X}}_{\rm{ACA}}$ is the approximation using only the information from TTACA, and $\tilde{\mathcal{X}}_{\rm{PEID}}$ is the resulting approximation of~\cref{alg:TT_update_parallel,alg:TT_update_sequential,alg:TT_update_parallel_2sided,alg:TT_update_sequential_2sided,alg:TT_update_rounding}. For all tests, the core ranks used are adaptively determined by specifying a tolerance for the magnitudes of the pivot selected in the ACA process.~\Cref{sec:PEID_Performance,sec:TTSVD_Comparisons} use the same TT ranks that TTACA produces for all PEID algorithms for a fair comparison.~\Cref{sec:Noisy_test} allows PEID to use an adaptive rank selection. This is done by specifying a PEID tolerance level that controls a SVD truncation added in the orthornormalization phase of~\cref{alg:TT_update_parallel,alg:TT_update_sequential,alg:TT_update_parallel_2sided,alg:TT_update_sequential_2sided,alg:TT_update_rounding}. For example, in~\cref{alg:TT_update_parallel}, at line 6 to determine an ON basis, we perform QR-SVD with a tolerance for SVD truncation, which may alter the size of $U_j$ and therefore the final TT ranks. The remaining algorithms follow the similar logic for the adaptive PEID procedure. All reported ranks are the averaged TT ranks across runs rounded to the nearest integer.

For the oversampling parameter, we use a constant sampling parameter $p$ for all dimensions. In the event that the number of free indices to be sampled in a specific dimension is less than the sampling parameter, we simply take the minimum between the oversampling parameter and the number of free indices. For example, in a $10$-$d$ example with mode size $50$, if the core ranks used are $(1,15,...)$ then the first dimension only has $35$ free indices. Thus, even if we set $p = 40$, $\mathcal{K}_{\le 1}$ can still only consist of all $35$ free indices. Finally, all plots that use sample size as a parameter report the average sample sizes across dimensions.

\subsection{Tensor Definitions} \label{sec:TensorDef}
In this section we provide a brief summary of the baseline test tensors for our three sets of tests.
\subsubsection{Hilbert Tensor}\label{sec:HilbertTests}
The first tensor we test on is the Hilbert tensor. It is a common benchmark test, as it is known to have low TT ranks~\cite{shi2021compressibility}. The Hilbert tensor is defined entry-wise as
\begin{equation}\label{eq:Hilbert_tensor}
    \mathcal{X}(i_1,i_2,\dots,i_d) = \frac{1}{i_1+i_2+\dots+i_d +1-d },\quad 1\leq i_j\leq n,\quad 1\leq j\leq d.
\end{equation}
If indexing starts at $0$, we can remove the $-d$ from the denominator. We consider two Hilbert tensors of dimensions $4$ and $10$, with uniform mode size $1600$ and $200$, respectively.

\subsubsection{Kernel-Based Tensors}
The next sets of experiments are on tensors corresponding to kernel function evaluations. Both cases are also studied in~\cite{khan2025parametric}. To define such tensors, we assume that the dimension $d$ is even, and $\pmb{x},\pmb{y}$ are vectors in $ [a_0,a_1]^{d/2}$ and $[b_0,b_1]^{d/2}$ respectively. Furthermore, let $\zeta_i$ be the $i$th Chebyshev node on $[-1,1]$ with $1 \le i \le n$, $\xi_i = \frac{a_0-a_1}{2}\zeta_i + \frac{a_0+a_1}{2}$, and  
 $\eta_j = \frac{b_0-b_1}{2}\zeta_j + \frac{b_0+b_1}{2}$. Using these components, we can form vectors $\pmb{x}_{i_1\cdots i_{d/2}} = \left[\xi_{i_1},\cdots,\xi_{i_{d/2}}\right]$ and $\pmb{y}_{j_1\cdots j_{d/2}} = \left[\eta_{j_1},\cdots,\eta_{j_{d/2}}\right]$ for $1 \le i_1,\cdots,i_{d/2},j_1,\cdots,j_{d/2}\le n$, and compute a scalar $\gamma_{i_1,\cdots,i_{d/2},j_1,\cdots,j_{d/2}} = \big|\big|\pmb{x}_{i_1\cdots i_{d/2}}-\pmb{y}_{j_1\cdots j_{d/2}}\big|\big|$. Here, $\big|\big|\cdot\big|\big|$ denotes the vector Euclidean norm.

Our first studied kernel function is the Mat\'ern kernel, which uses the order-$(5/2)$ modified Bessel function of the second kind $K_{5/2}(\cdot)$, and the resulting tensor is
\begin{equation}\label{eq:Matern}
    \mathcal{X}\left(i_1,\cdots,i_{d/2},j_1,\cdots,j_{d/2}\right)
    = \frac{2^{-3/2}}{\Gamma(5/2)}
    \left(\sqrt{5}\gamma_{i_1,\cdots,i_{d/2},j_1,\cdots,j_{d/2}}\right)^{5/2}
    K_{5/2}\left(\sqrt{5}\gamma_{i_1,\cdots,i_{d/2},j_1,\cdots,j_{d/2}}\right).
\end{equation}
The second kernel is referred to as the Thin Plate Spline kernel, denoted as TPS in test results, and it involves a logarithm to build the tensor
\begin{equation}\label{eq:ThinPlateSpline}
    \mathcal{X}\left(i_1,\cdots,i_{d/2},j_1,\cdots,j_{d/2}\right) = \gamma_{i_1,\cdots,i_{d/2},j_1,\cdots,j_{d/2}}^2\log\left(\gamma_{i_1,\cdots,i_{d/2},j_1,\cdots,j_{d/2}}^2\right).
\end{equation}
In our tests, we use $[a_0,a_1] = [2,3]$, and $[b_0,b_1] = [0,1]$. These are chosen to be disjoint to avoid any singularities in tensor entries.

\subsubsection{Maxwellian distribution from gas dynamics}
For our last set of experiments, we test on tensors built from Maxwellian distributions, which naturally arise in gas dynamics~\cite{einkemmer2025review}. These tensors are defined entry-wise by function evaluations on a prescribed discretization grid of displacement and velocity domains. Here, we define the $4$-$d$ tensor, and the $6$-$d$ version is a natural extension. For $4$-$d$, we assume that $x,y\in [-1/2,1/2]$, $v_x,v_y\in[-3,3]$, and define
\begin{equation} \label{eq:2d2vFunction}
f(x,y,v_x,v_y) = \varrho(x,y)\left[\exp{\left(-b_x^- - b_y^-\right)}+\exp{\left(-b_x^+ - b_y^+\right)}\right],
\end{equation}
where
$$
\varrho(x,y) = \left(\frac{\rho(x)}{2\sqrt{2\pi T(x)}} + \frac{\rho(y)}{2\sqrt{2\pi T(y)}}\right),\quad b_x^{\pm} = \frac{|v_x \pm 0.75|^2}{2T(x)},\quad b_y^{\pm} = \frac{|v_y \pm 0.75|^2}{2T(y)},
$$
with
$$
\rho(w)=1+0.875\sin(2\pi w),\quad T(w)= 0.5 + 0.4\sin(2\pi w), \quad w = x,y.
$$
By taking $n$ equi-distanced points on each domain, we can create a $4$-$d$ grid of size $n\times n\times n\times n$, and we can calculate the entries of the Maxwellian tensor $\mathcal{X}$ by
$$
\mathcal{X}(i_1,i_2,i_3,i_4) = f(x_{i_1},y_{i_3},(v_x)_{i_2},(v_y)_{i_4}).
$$

\subsection{Comparison with other methods}\label{sec:TTSVD_Comparisons}
This section focuses on 4--$d$ Hilbert, Maxwellian, Mat\'ern, and Thin Plate Spline tensors with mode size $100$. We compare the performance of the proposed PEID algorithms against TTSVD through Tensorly, the TnTorch cross function, TT-Sketching, and DMRG through Matlab. TnTorch and DMRG are interpolation-type methods, while TT-Sketching is a randomized method based on orthogonal projections. For our tests, we first run TTACA, and use the output TT ranks for TTSVD, TnTorch, and TT-Sketching. Matlab DMRG is a tolerance based implementation, so its results are produced by manually selecting a tolerance to obtain nearly identical TT ranks as TTACA. In this way, these experiments compare approximation accuracy at matched or approximately matched TT ranks. We do not report cross-method timing comparisons here because the implementations use different software environments and computational backends, which reflect implementation choices as well as algorithmic cost.

\begin{figure}[!htbp]
    \centering
    \includegraphics[width=0.8\linewidth]{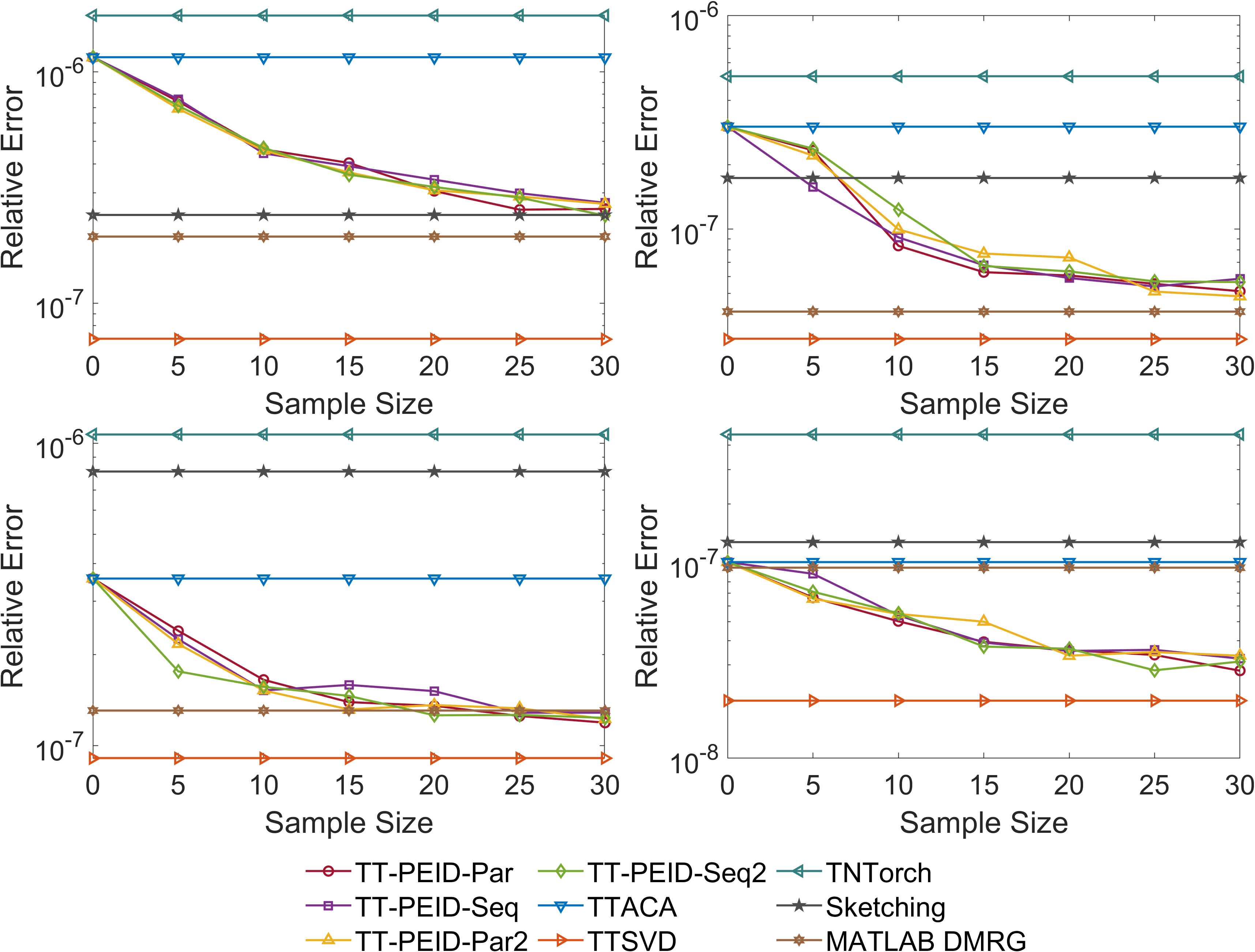}
    \caption{Comparison of PEID against common tensor compression techniques. Error versus Sample size for 4-$d$ tensors with mode size $100$. Top left: Hilbert - rank $(11,11,11)$, top right: 2d2v Maxwellian - rank $(12,5,22)$, bottom left: Mat\'ern - rank $(7,12,6)$, and bottom right: Thin Plate Spline - rank $(6,12,6)$. }
    \label{fig:TTSVD_comparison}
\end{figure}

In these tests, PEID algorithms are run 100 times and the median error values are plotted.~\Cref{fig:TTSVD_comparison} plots the error versus sample size for all four cases. For these plots, TTACA, TTSVD, TnTorch, Sketching, and DMRG are run only once and the errors are extended across the plot as horizontal lines for reference.~\Cref{fig:TTSVD_comparison} shows that in all cases, as sample size increases, the observed errors of PEID algorithms decrease. TTSVD achieves the smallest errors in all tests, but at 30 oversampling points, the PEID algorithms attain errors that are generally comparable to or smaller than those of the other methods.

\subsection{Performance of PEID}\label{sec:PEID_Performance}
In this section, we illustrate the accuracy and computational efficiency of the proposed PEID algorithms. Similar to~\cref{sec:TTSVD_Comparisons}, all tests in this section are run 100 times and the reported results are the median over all runs. The sample sizes used for testing are $5,10,\ldots,50$. TTACA ranks are determined separately for each trial and may vary between trials. The tolerance values used for TTACA and PEID along with the average ranks observed are reported with the results.

\subsubsection{Hilbert Tensor}
\begin{figure}[!htbp]
    \centering
    \includegraphics[width=0.8\linewidth]{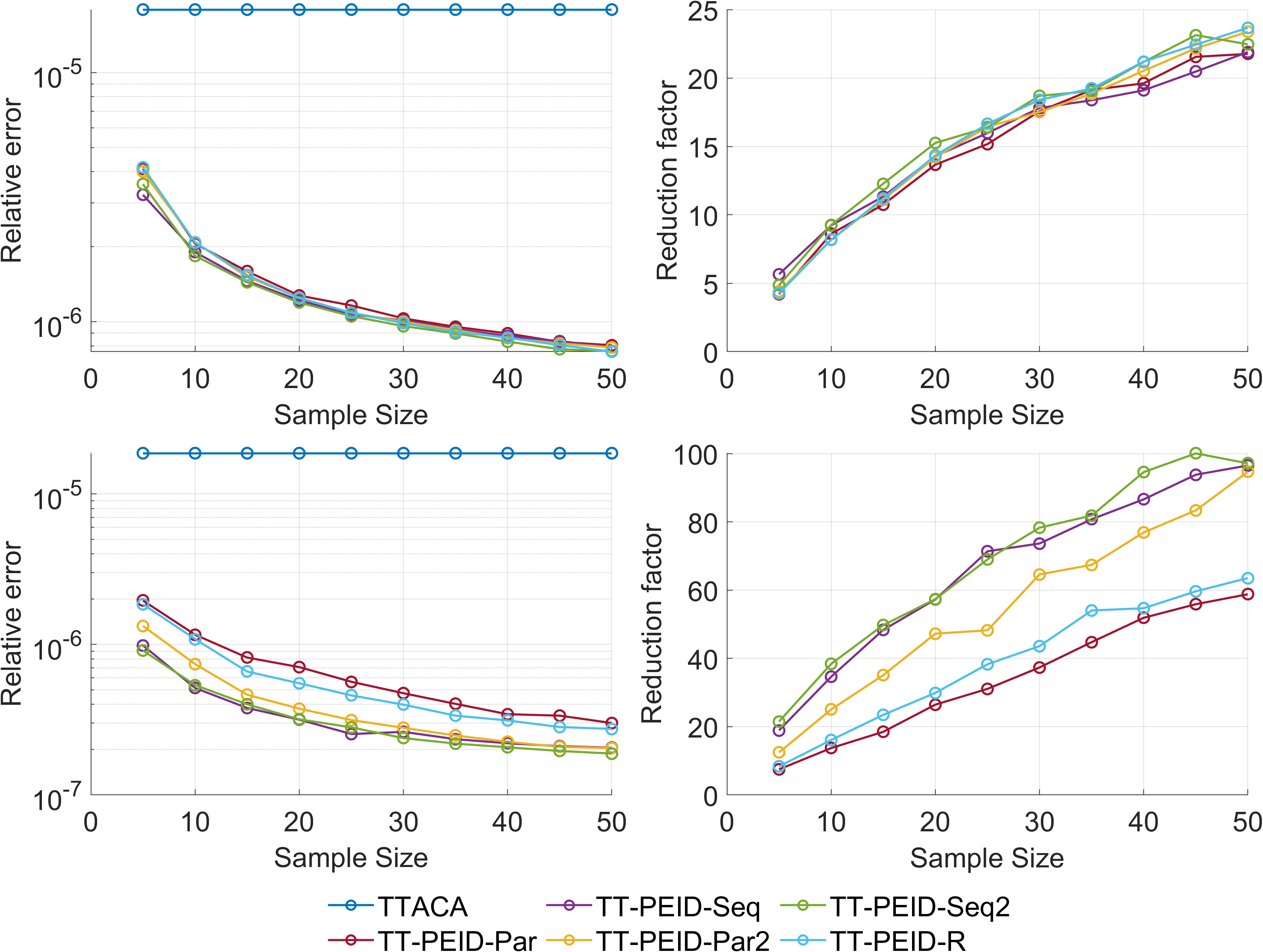}
    \caption{Relative error and reduction factor of Hilbert tensor. Top row: 4-$d$ Hilbert tensor with mode size 1600 and average observed ranks $(15,15,15)$. Bottom row: 10-$d$ Hilbert tensor with mode size 200 and average observed ranks $(13, 14, 14, 14, 14, 14, 14, 14, 13)$.}
    \label{fig:HilbertErrorRFvsSampleSize}
\end{figure}

The first test measures the reduction in error for increasing sample sizes $p$ for a fixed TTACA input. In this case, the tolerance value used for TTACA was $10^{-7}$ for both 4-$d$ and 10-$d$.~\Cref{fig:HilbertErrorRFvsSampleSize} shows the results for both 4-$d$ and 10-$d$ tests. On the top row, for the 4-$d$ results, we see that as sample size increases, the observed reduction in error increases, up to approximately 24 times reduction for 50 oversample points per dimension. For this test case,~\cref{alg:TT_update_parallel,alg:TT_update_sequential,alg:TT_update_parallel_2sided,alg:TT_update_sequential_2sided,alg:TT_update_rounding} perform similarly. The second row plots the 10-$d$ test case. We observe again that with increasing sample sizes, the reduction increases, up to approximately 100 times. In this test, there is a preference towards~\cref{alg:TT_update_sequential,alg:TT_update_sequential_2sided}. However, all other algorithms obtain good reductions in error. For example, with 50 oversampling points, the reduction factor is at least 55.

\begin{figure}[!htbp]
    \centering
    \includegraphics[width=0.8\linewidth]{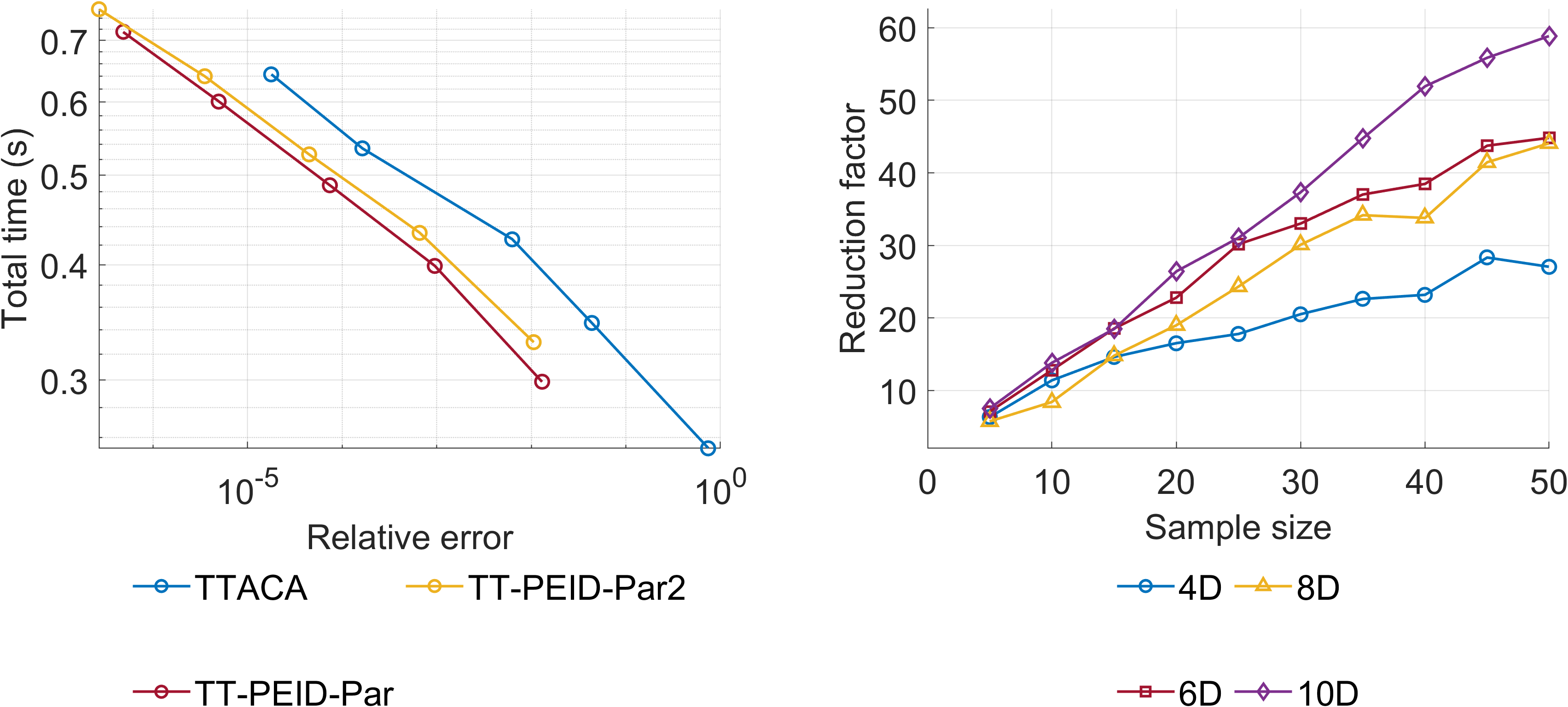}
    \caption{Left: Time in seconds taken to run vs relative error achieved for TTACA plus~\cref{alg:TT_update_parallel,alg:TT_update_parallel_2sided} on a 10-$d$ Hilbert tensor with mode size 200. Right: Reduction factor of \cref{alg:TT_update_parallel} for $4,6,8,10$ dimensional Hilbert tensors with mode sizes fixed at 200 vs sample size.}
    \label{fig:HilbertTimeandDim}
\end{figure}

Next, we measure the time taken in seconds to run just TTACA, and TTACA followed by~\cref{alg:TT_update_parallel,alg:TT_update_parallel_2sided} for a range of TTACA tolerance values.
The purpose of this test is to show that in terms of computation time, it is advantageous to run TTACA with a larger tolerance, followed by oversampling to reduce approximation error. The left plot in~\cref{fig:HilbertTimeandDim} demonstrates this, as the blue curve for TTACA lies above the curves for~\cref{alg:TT_update_parallel,alg:TT_update_parallel_2sided}. For the error levels considered in this experiment, the results indicate that using a less restrictive TTACA tolerance followed by PEID can achieve a desired approximation accuracy at a lower computational cost than TTACA alone. For example, in this plot, for a computation time of roughly 0.4 seconds, TTACA with~\cref{alg:TT_update_parallel} can achieve an error of $10^{-3}$ while TTACA alone achieves an error of $10^{-2}$. The right plot of~\cref{fig:HilbertTimeandDim} contains results when mode size is fixed at 200 and dimension is varied. We observe that for Hilbert tensors, there are generally more advantages in oversampling as dimension increases. 

\subsubsection{Kernel-Based Tensors}

\begin{figure}[!htbp]
    \centering
    \includegraphics[width=0.8\linewidth]{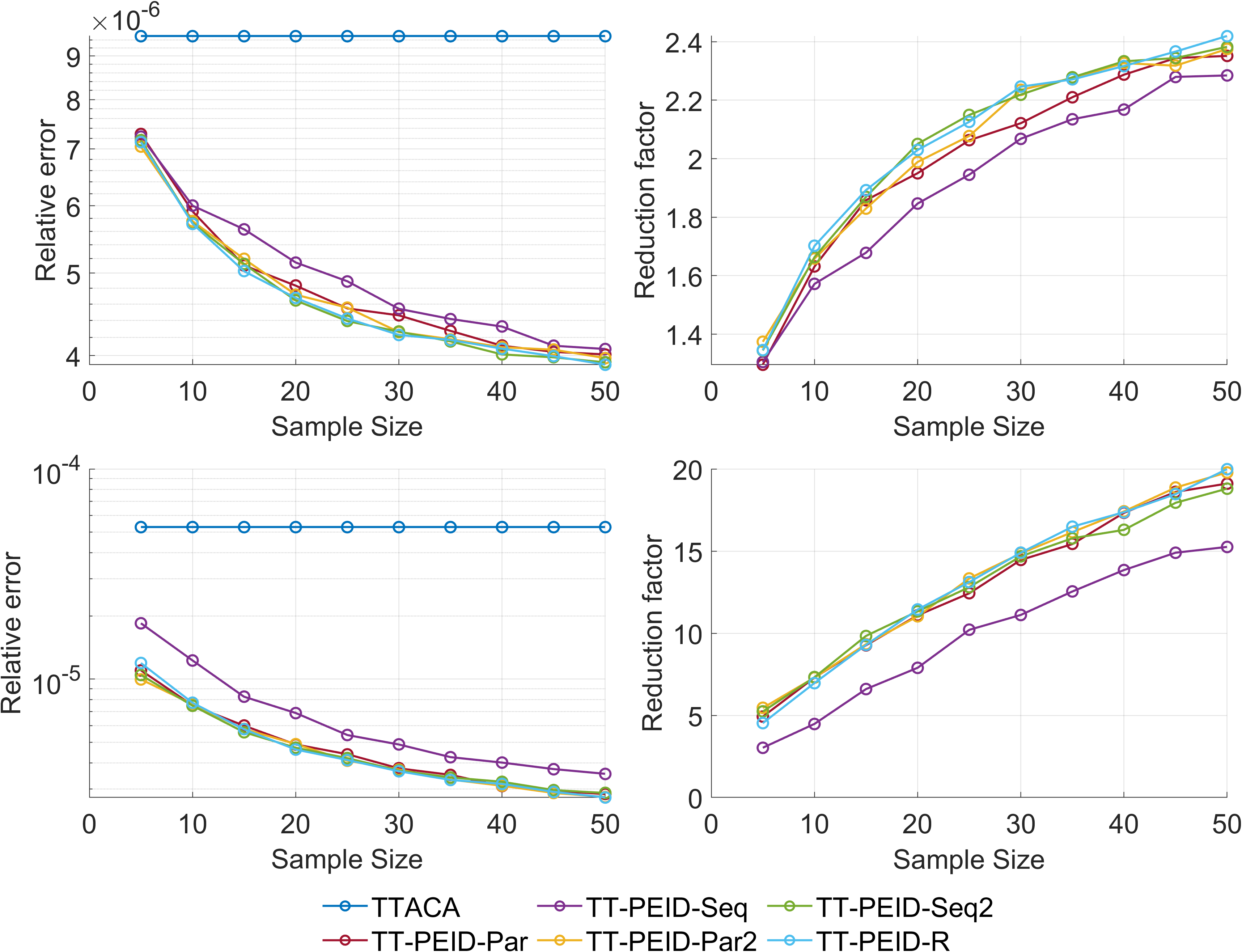}
    \caption{Relative error and reduction factor of Mat\'ern tensor. Top row: 4-$d$ Mat\'ern tensor with mode size 1600 and average observed ranks $(5, 8, 5)$. Bottom row: 10-$d$ Mat\'ern tensor with mode size 200 and average observed ranks $(5, 15, 35, 77, 118, 79, 36, 16, 5)$.}
    \label{fig:MaternErrorRFvsSampleSize}
\end{figure}

As with the Hilbert tensor, our first set of tests measure the error reduction for 4-$d$ and 10-$d$ Mat\'ern tensors with mode sizes 1600, 200 and tolerance values $10^{-4}$ and $10^{-6}$ respectively using different oversampling size. The results for these tests are shown in~\cref{fig:MaternErrorRFvsSampleSize}. We see that the error reduction increases as sample size increases. For 4-$d$ we reduce the error by approximately 2.4 times, and for 10-$d$ we reduce the error by up to approximately 20 times. Although this is not as drastic as the Hilbert tensor, it is noteworthy that the Mat\'ern tensor is a more difficult test case, since the TT ranks observed are significantly larger than those of the Hilbert tensor. 

\begin{figure}[!htbp]
    \centering
    \includegraphics[width=0.8\linewidth]{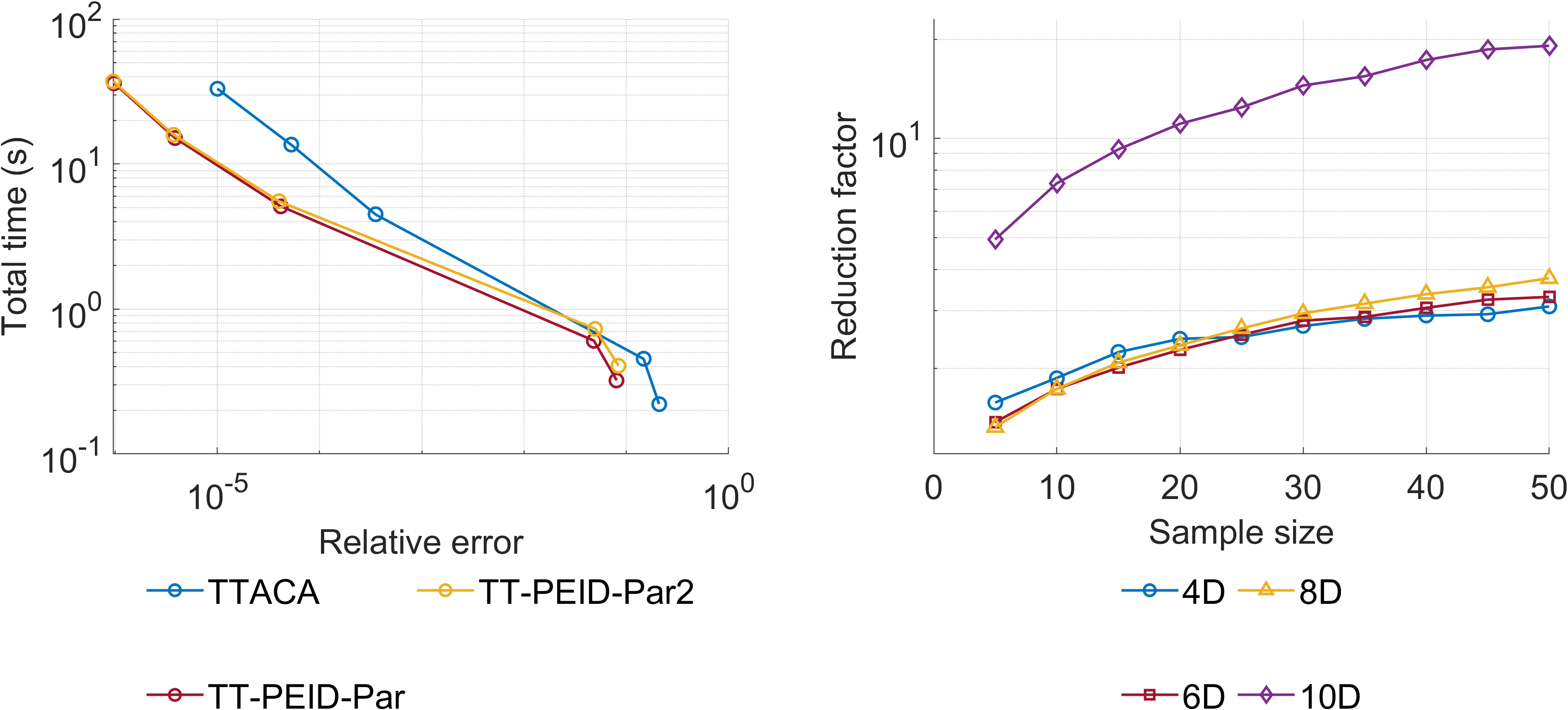}
    \caption{Left: Time in seconds taken to run vs relative error achieved for TTACA plus~\cref{alg:TT_update_parallel,alg:TT_update_parallel_2sided} on a 10-$d$ mode size 200 Mat\'ern tensor. Right: Reduction factor of \cref{alg:TT_update_parallel} for $4,6,8,10$ dimensional Mat\'ern tensors with mode sizes fixed at 200 vs sample size.}
    \label{fig:MaternTimeandDim}
\end{figure}

Due to the larger ranks of the Mat\'ern tensors, run time is larger and is reflected by the scale of the vertical axis in the left plot of~\cref{fig:MaternTimeandDim}. This plot shows the timing results for the 10-$d$ case over a range of TTACA tolerance levels. In this test, we see a similar benefit to the Hilbert tensor case where TTACA with~\cref{alg:TT_update_parallel} can reach smaller errors for the same runtime as TTACA alone. For example, at a runtime of roughly 35 seconds, TTACA with~\cref{alg:TT_update_parallel} achieves an error of $10^{-6}$ while TTACA alone can only achieve error $10^{-5}$. In the right plot of~\cref{fig:MaternTimeandDim}, which contains results where mode size is fixed at 200 and dimension is varied, we can see that for dimensions 4, 6, and 8, we see a reduction of error of about 1.5--2 times while a significantly larger benefit is observed in the 10-$d$ case. The error bounds offer an explanation of this behavior. By~\cref{thm:alg3_main}, the error of the projection enhanced approximation at each separation is governed by the smaller of two quantities: the skeleton error $\beta_k$ inherited from the TTACA pivots, and the oversampling term $(1+c_1\delta^{-1/2})\Delta_k$, which is a constant multiple of the optimal rank-$r$ error $\Delta_k$ and does not depend on the pivot quality at all. The TTACA approximation, in contrast, is controlled by the skeleton error alone (see \cref{thm:ttaca_app}). Therefore, the reduction factor is large when the greedy skeleton is far from optimal, i.e., when $\beta_k$ exceeds $\Delta_k$ by a large factor. For the $4$-dimensional Mat\'ern tensor in particular, the observed ranks are $(5,8,5)$, which indicates that the singular values of the unfoldings decay rapidly, and the greedy pivots are already nearly optimal, so both approximations are close to $\Delta_k$ and the reduction factor stays near $2$. For the $10$-dimensional case, the middle ranks exceed $100$, so the singular values decay slowly, the greedy pivoting loses accuracy relative to the optimal error, and the oversampling branch of the bound takes over. What is more, the two bounds also differ in how errors propagate across separations. In~\cref{thm:ttaca_app}, the error at separation $k$ is amplified by the factor $(r2^{r-1})^{k-1}$, whereas in~\cref{thm:PEID_format_main,thm:alg3_main} it is amplified by $\prod_{j<k}\gamma_j$ and $\prod_{j\le k}\eta_j$ respectively. Whenever the coefficient norms $\gamma_j$, or their sampled surrogates $\eta_j$, are smaller than $r2^{r-1}$, generic for well-chosen pivots, the ratio of the two amplification factors grows geometrically with $k$. As a result, the bounds predict a reduction factor that increases with the dimensionality.

\begin{figure}[!htbp]
    \centering
    \includegraphics[width=0.8\linewidth]{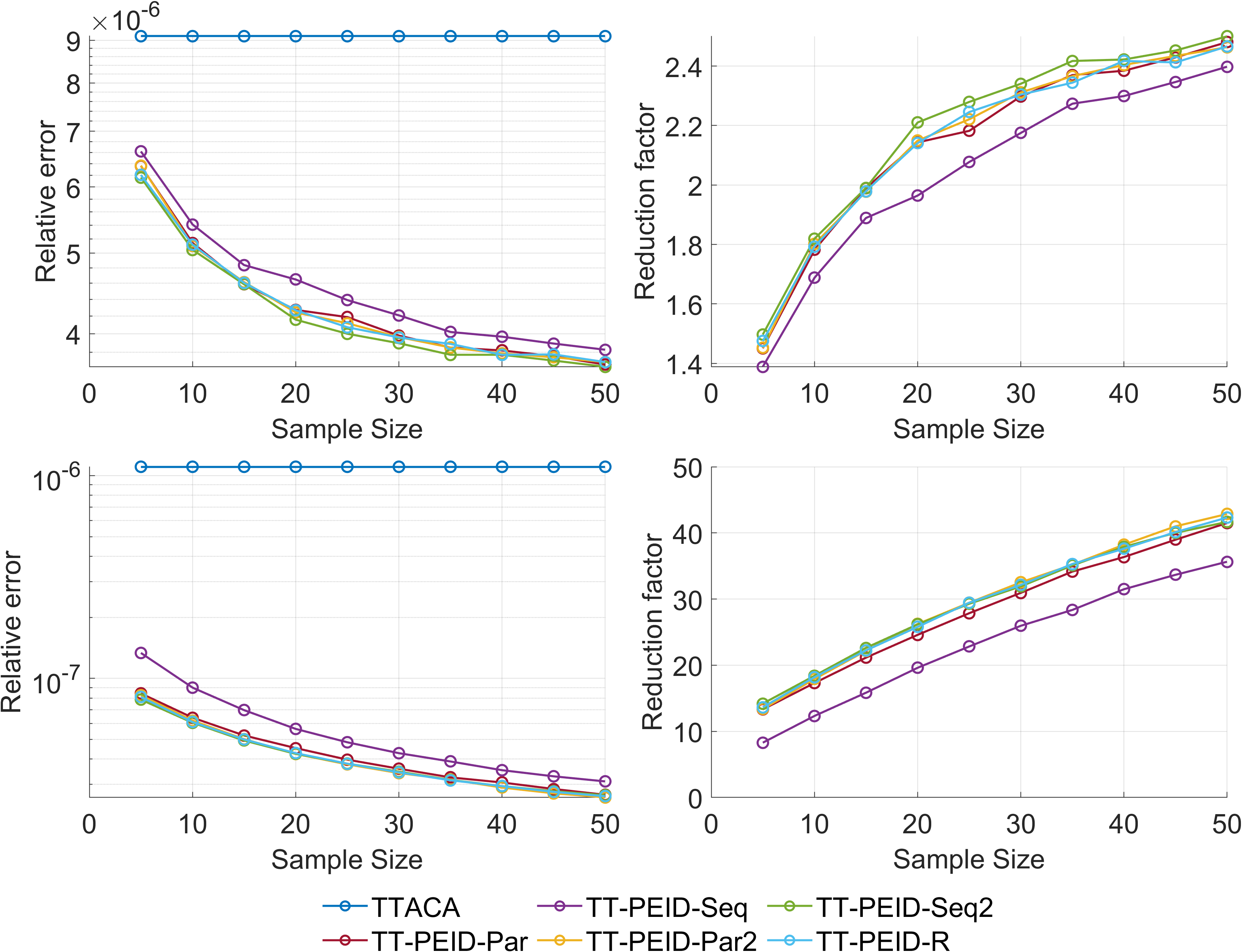}
    \caption{Relative error and reduction factor of TPS tensor. Top row: 4-$d$ TPS tensor with mode size 1600 and average observed ranks $(5,7,5)$. Bottom row: 10-$d$ TPS tensor with mode size 200 and average observed ranks $(5, 16, 32, 65, 100, 65, 32, 15, 6)$.}
    \label{fig:TPSErrorRFvsSampleSize}
\end{figure}

\begin{figure}[!htbp]
    \centering
    \includegraphics[width=0.8\linewidth]{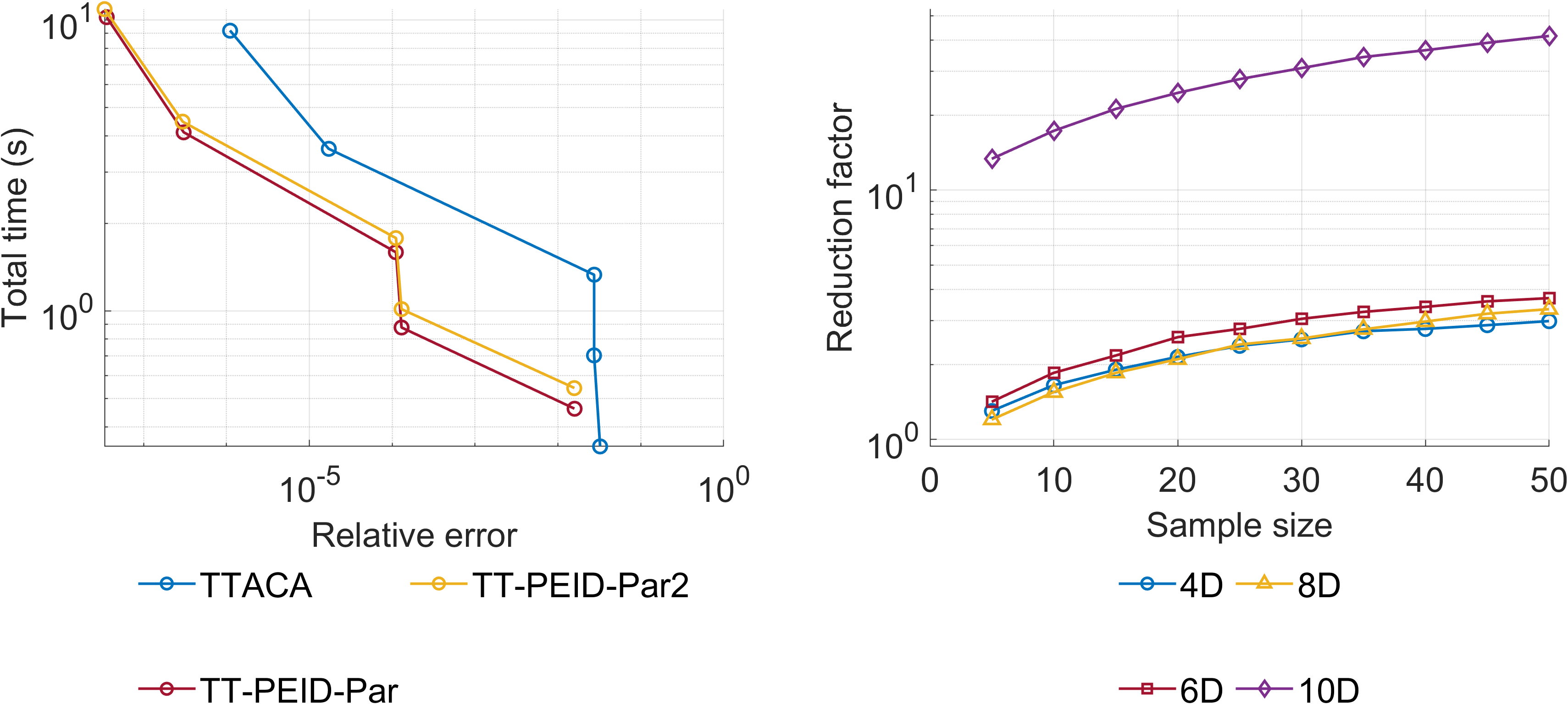}
    \caption{Left: Time in seconds taken to run vs relative error achieved for TTACA plus~\cref{alg:TT_update_parallel,alg:TT_update_parallel_2sided} on a 10-$d$ TPS tensor with mode size 200. Right: Reduction factor of \cref{alg:TT_update_parallel} for $4,6,8,10$ dimensional TPS tensors with mode size fixed at 200 vs sample size.}
    \label{fig:TPSTimeandDim}
\end{figure}

~\Cref{fig:TPSErrorRFvsSampleSize,fig:TPSTimeandDim} provides the results of the Thin Plate Spline for the same test cases as Hilbert and Mat\'ern with tolerance values $10^{-4}$ and $10^{-7}$ for 4-$d$ and 10-$d$. From~\cref{fig:TPSErrorRFvsSampleSize}, we see that performing oversampling reduces the error by roughly 2.5 times in the 4-$d$ case. However, in the 10-$d$ case, we are able to reduce the error by up to approximately 43 times using~\cref{alg:TT_update_parallel_2sided}. For the 10-$d$ Thin Plate Spline tensor, we notice that pure TTACA tends to stall roughly around relative error $10^{-2}$, until the TT ranks reach a high enough threshold, when the error quickly decreases. This stalling behavior can be observed in the left plot of~\cref{fig:TPSTimeandDim}, where TTACA error remains relatively unchanged with lower tolerances until the threshold is reached. In fact, this stalling is a pivot-selection failure rather than a loss of low rank structure. The optimal errors $\Delta_k$ keep decreasing as the ranks grow, but the greedy pivots fail to capture the dominant directions in this scenario, so the skeleton error $\beta_k$ stays flat. The oversampling branch of \cref{thm:alg3_main}, however, is insensitive to this failure since uniformly sampled fibers see the dominant subspaces whenever those subspaces are incoherent. Furthermore, when looking at the time scale of the left plot in~\cref{fig:TPSTimeandDim}, we find that reducing the error by two orders of magnitude can provide significant time savings compared to simply running TTACA with a strict tolerance. In the right plot of~\cref{fig:TPSTimeandDim}, we can see significant advantages when using~\cref{alg:TT_update_parallel} as the dimension increases for Thin Plate Spline tensors. This behavior agrees with the explanation of TTACA stalling in 10-$d$.

\subsubsection{Maxwellian Distribution}
\begin{figure}[!htbp]
    \centering
    \includegraphics[width=0.8\linewidth]{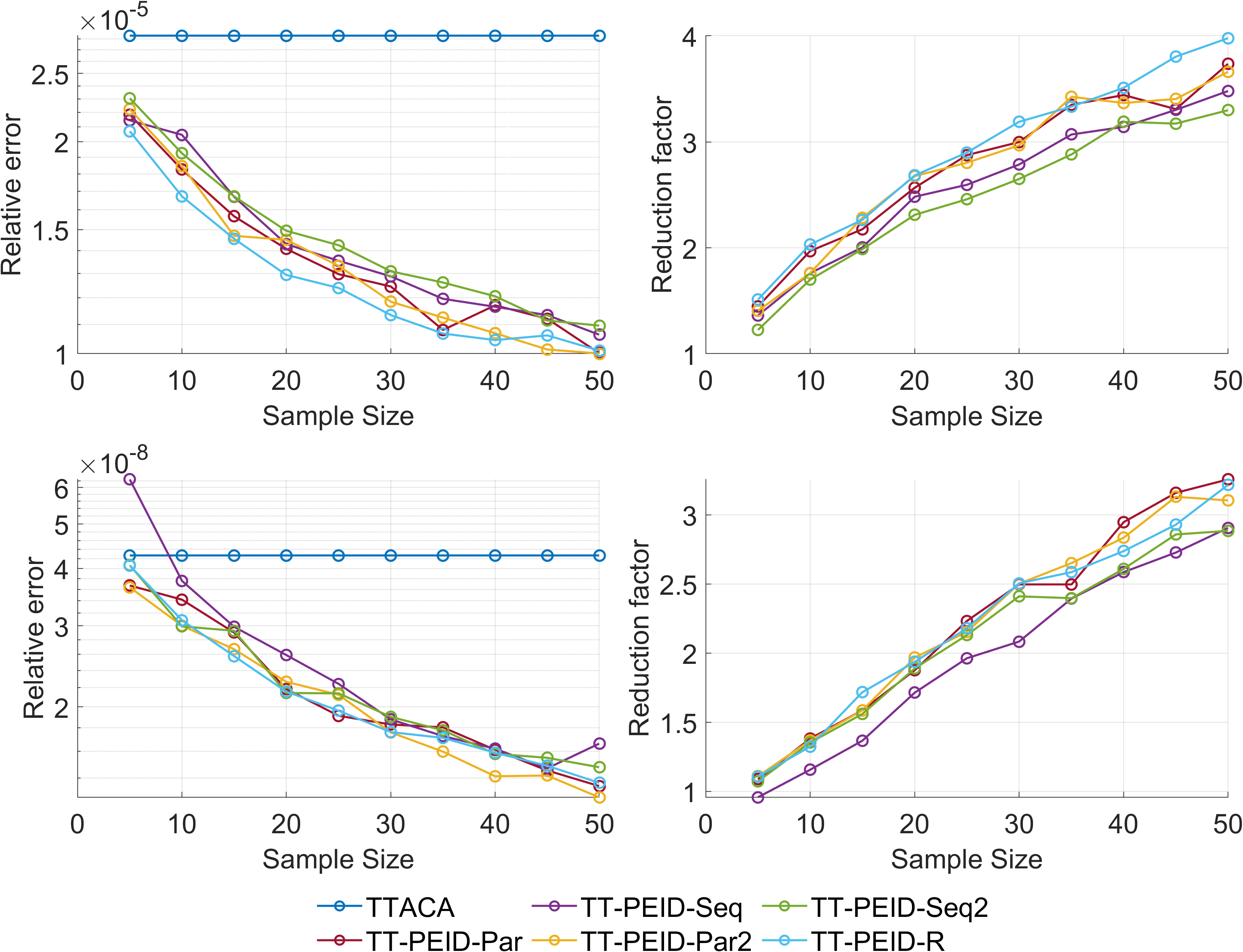}
    \caption{Relative error and reduction factor of Maxwellian tensor. Top row: 4-$d$ Maxwellian tensor with mode size 1600 and average observed ranks $(9, 5, 15)$. Bottom row: 6-$d$ Maxwellian tensor with mode size 800 and average observed ranks $(12, 6, 37, 5, 22)$.}
    \label{fig:maxwellianErrorRFvsSampleSize}
\end{figure}

For the Maxwellian test case, we present the results for 4-$d$ mode size 1600 and 6-$d$ mode size 800 tensors with tolerances $10^{-4}$ and $10^{-6}$ in~\cref{fig:maxwellianErrorRFvsSampleSize} as these are the only physically relevant dimensions of interest. The top row shows 4-$d$ results, which exhibit similar behavior to all previous tests, and we observe a reduction in error of approximately 4 times. The bottom two plots display similar behavior for 6-$d$ with reduction in error of roughly 3.3 times.

\subsection{Performance with noisy data}\label{sec:Noisy_test}
\begin{figure}[!htbp]
    \centering
    \includegraphics[width=0.8\linewidth]{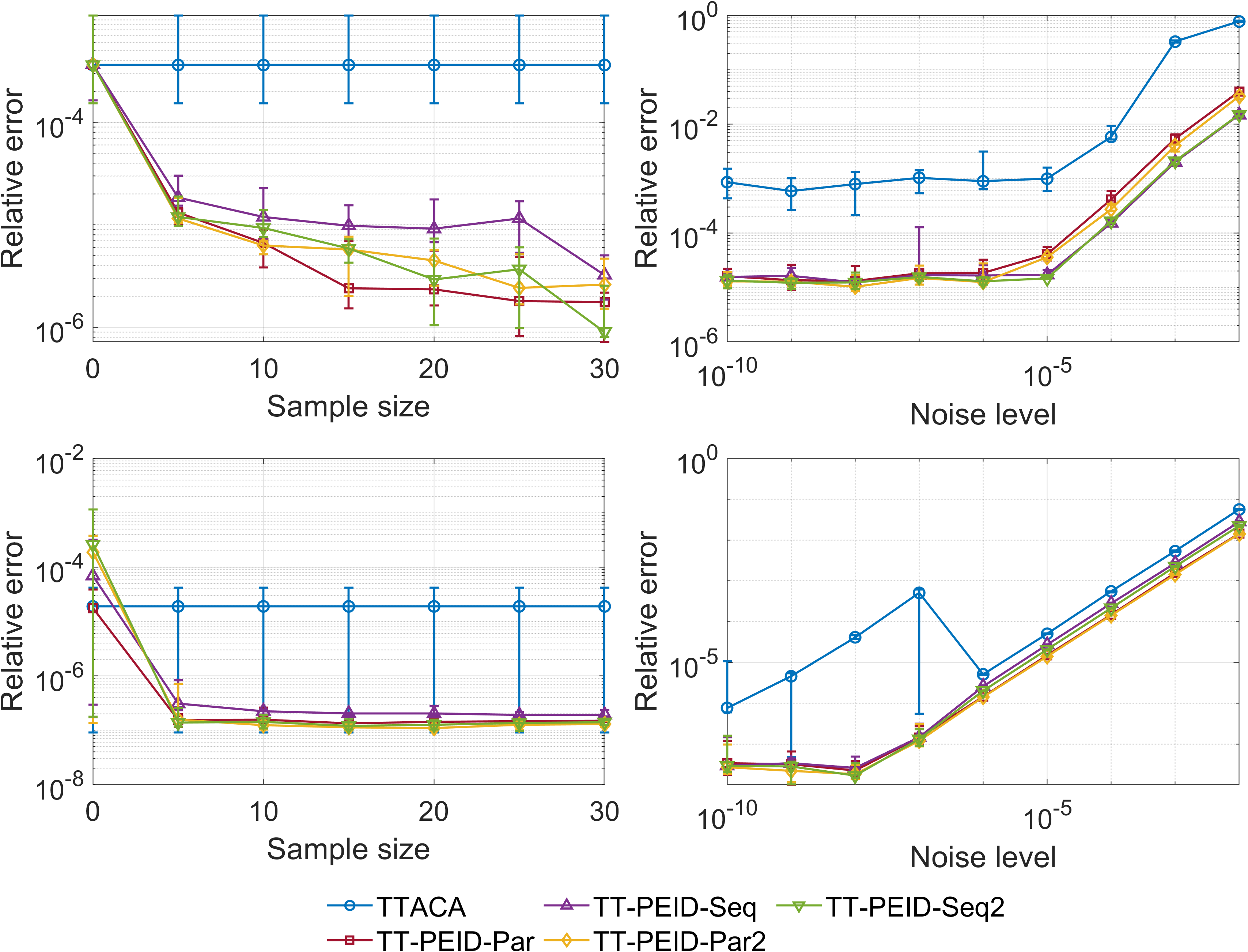}
    \caption{Relative error for 6D Hilbert (top) and Maxwellian (bottom) tensors with mode size $800$ with noise added. Left plots have a fixed noise size of $\mu = 10^{-8}$ and TTACA tolerance of $10^{-6}$. For Hilbert, average rank for TTACA and TT-PEID-Par with 30 oversample points are $(12,13,13,13,12)$ and $(13,14,15,15,14)$. Average Maxwellian rank for TTACA and TT-PEID-Par with 30 oversample points are $(13, 7,37,6,22)$ and $(10,5,30,5,18)$. Right plots have a fixed sample size of 30 and TTACA tolerance of $10^{-6}$.}
    \label{fig:Noisy_results}
\end{figure}

This section tests the robustness of PEID in the presence of noisy data. For these tests, our target tensor takes the form $(1 + \mu\mathcal{Y})\mathcal{X}$ entrywise, where $\mu$ is a noise level scaling factor, and $\mathcal{Y}$ is a tensor with random Gaussian entries with mean 0 and standard deviation 1. The reported errors are measured with the clean data, i.e. $||\mathcal{X} - \tilde{\mathcal{X}}||_{F,R}/||\mathcal{X}||_{F,R}$ where $\tilde{\mathcal{X}}$ is the TT approximation of the noisy tensor $(1+\mu\mathcal{Y})\mathcal{X}$. In addition, we run TTACA with tolerance $10^{-6}$. In test cases where $\mu>10^{-6}$, single runs experience large run times due to large observed ranks. Therefore, each experiment is repeated only 10 times, and the reported results are the median relative errors with error bars representing the 25th and 75th percentiles.~\Cref{fig:Noisy_results} shows the results for 6-$d$ Hilbert (top row) and 6-$d$ Maxwellian (bottom row) tensors on two types of experiments: (1) we fix the noise level at $10^{-8}$ and vary the sample size, and (2) we fix the sample size at 30 and vary the noise level. For the fixed noise level tests, PEID uses an adaptive truncation tolerance of $10^{-7}$, while no adaptive truncation is applied in the varying noise level tests.

The first column of plots corresponds to fixed noise level. We see that in the presence of noisy data, PEID reduces error as well as the variation, as demonstrated by the overall shrinking sizes of the error bars. A similar conclusion can be drawn in the second column of plots where sample size is fixed at 30 and $\mu$ is varied. In both columns, for almost all plot points, the error bars of PEID are either comparable to or smaller than those of TTACA, aside from a noise level of $10^{-7}$ for the Hilbert test case in the top right plot. We also see that for both TTACA and PEID, for $\mu$ values less than the TTACA tolerance of $10^{-6}$, the methods appear to achieve comparable error measurements regardless of the presence of noise in the data. For $\mu$ values larger than the TTACA tolerance, we see the same trend in TTACA and PEID errors, with PEID recovering the approximation up to the noise level of the data.

\section{Acknowledgement} The authors would like to thank Dr. Lexing Ying in suggesting the references \cite{engquist2007fast, cortinovis2024sublinear} and enlightening discussions. Large language models are used to assist in the development of theoretical results and implementations. All mathematical statements, derivations, and proofs are independently verified and finalized by the authors.

\appendix

\section{Existing TT decomposition algorithms}
In this appendix, we outline TT-sketching algorithms in sequential and dimension parallel setting in~\Cref{alg:TTsketching,alg:parallel_TTsketching}. This helps readers to understand our proposed PEID algorithms in~\cref{sec:algorithm}.
\begin{algorithm}
\caption{TT-Sketching: Given a tensor, compute an approximant tensor in TT format using randomized sketching. }
\begin{algorithmic}[1]
\label{alg:TTsketching}
\Require {A tensor $\mathcal{X} \in \R^{n_1 \times \dots \times n_d}$ and desired TT ranks $r_1,\dots,r_{d-1}$.}
\Ensure {TT cores $\mathcal{G}_1, \dots, \mathcal{G}_d$ of an approximant $\tilde{\mathcal{X}}$.}
\State Set $Y = X_1$ the first unfolding of $\mathcal{X}$.
\State Set $r_0 = 1$.
\For {$1 \le j \le d-1$}
\State Use randomized range finder to find $U_j$ with $r_j$ ON columns of $Y$.
\State Set $\mathcal{G}_j = {\rm reshape}(U_j, r_{j-1}, n_j, r_j)$.
\State Set $Y = {\rm reshape}\left(U_j^TY,r_jn_{j+1},\prod_{k=j+2}^d n_k\right)$.
\EndFor
\State Set $\mathcal{G}_d = {\rm reshape}(Y, r_{d-1}, n_d, 1)$.
\end{algorithmic}
\end{algorithm}

\begin{algorithm}
\caption{Parallel-TT-Sketching: Given a tensor, compute an approximant tensor in TT format using randomized sketching with parallelization in dimensionality. }
\begin{algorithmic}[1]
\label{alg:parallel_TTsketching}
\Require {A tensor $\mathcal{X} \in \R^{n_1 \times \dots \times n_d}$ and desired TT ranks $r_1,\dots,r_{d-1}$.}
\Ensure {TT cores $\mathcal{G}_1, \dots, \mathcal{G}_d$ of an approximant $\tilde{\mathcal{X}}$.}
\For {$1 \le j \le d-1$}
\State Use randomized range finder to find $U_j$ with $r_j$ ON columns of $X_j$.
\EndFor
\For {$1 \le k \le d-2$}
\State Calculate $W_{k+1} = U_k^T \ {\rm reshape}(U_{k+1}, \prod_{i=1}^k n_i ,n_{k+1}r_{k+1})$.
\State Set $\mathcal{G}_{k+1} = {\rm reshape}(W_{k+1}, r_k, n_{k+1}, r_{k+1})$.
\EndFor
\State Set $\mathcal{G}_1 = U_1$ and $\mathcal{G}_d = U_{d-1}^TX_{d-1}$.
\end{algorithmic}
\end{algorithm}

\section{Error analysis of the projection enhanced TT format}\label{app:format_analysis}

Recall that we write
\[
C_k = X_k(:,\mathcal{J}_{>k}),\qquad R_k = X_k(\mathcal{I}_{\le k},:),\qquad 1\le k\le d-1,
\]
for the $k$th unfolding of a tensor $\mathcal{X}$. Also, let
\begin{equation}\label{eq:gamma_def}
\gamma_k := \left\|X_kR_k^\dagger\right\|_2\;\ge\;1
\end{equation}
denote the coefficient norm of the row interpolative decomposition induced by the pivot rows $\mathcal{I}_{\le k}$, a classical measure of row skeleton quality~\cite{cheng2005compression,gu1996efficient}. Then, the following result presents a complete version of \cref{thm:PEID_format_simplified}, and this appendix proves it. 
\begin{theorem}\label{thm:PEID_format_main}
For $1\le k\le d-1$, assume that $R_k$ has full row rank, and that the pivot pair $(\mathcal{I}_{\le k},\mathcal{J}_{>k})$ is produced by $r$ steps of ACA with complete pivoting applied to $X_k$, where $\rank(X_k)>r$. Then the TT approximation $\tilde{\mathcal{X}}$ defined by \cref{eq:ttcross_cores,eq:ttcross_PEID_build} satisfies
\[
\left\|\mathcal{X}-\tilde{\mathcal{X}}\right\|_F\;\le\;2\,n^{d/2}\,4^{r}\rho_r\sum_{k=1}^{d-1}\Bigl(\prod_{j=1}^{k-1}\gamma_j\Bigr)\,\sigma_{r+1}(X_k),
\]
where $\gamma_j$ is defined in~\cref{eq:gamma_def}, $\sigma_{r+1}(X_k)$ is the $(r+1)$st singular value of the $k$th unfolding, and $\rho_r\le2\sqrt{r+1}\,(r+1)^{\ln(r+1)/4}$ is the growth factor of Gaussian elimination with complete pivoting~\cite{wilkinson1961error,cortinovis2020maximum}.
\end{theorem}

Throughout, $\|\cdot\|_{\max}$ denotes the Chebyshev (entrywise maximum) norm of a matrix or tensor, which is unfolding-independent, as is the Frobenius norm. For a matrix $M$ we write $P_M:=MM^\dagger$ for the orthogonal projection onto its column space and $M_r$ for a best rank-$r$ approximation. In addition, we repeatedly use the singular value interlacing property
\begin{equation}\label{eq:interlacing_app}
\sigma_i\bigl(A(\mathcal{S},\mathcal{T})\bigr)\le\sigma_i(A)\qquad\text{for all }i\text{ and all row/column selections }\mathcal{S},\mathcal{T},
\end{equation}
see, e.g., \cite[Cor.~7.3.6]{horn2013matrix}, and its consequence $\|A(\mathcal{S},\mathcal{T})-A(\mathcal{S},\mathcal{T})_r\|_F\le\|A-A_r\|_F$.

\subsection{Factorization of the format into projections}
If we reshape $M_k:=X_{k-1}(\mathcal{I}_{\le k-1},\mathbb{I}_k\otimes\mathcal{J}_{>k})$ as an $r\times n\times r$ array as in \cref{eq:ttcross_cores}, then
\begin{equation}\label{eq:slice_identity}
M_k(:,i_k,:)\;=\;C_k\bigl(\mathcal{I}_{\le k-1}\otimes\{i_k\},\,:\bigr),
\end{equation}
which is valid entrywise as both sides equal $\mathcal{X}(\alpha,i_k,\mathcal{J}_{>k})$ for $\alpha\in\mathcal{I}_{\le k-1}$. Note that \cref{eq:slice_identity} uses no nestedness. 
Additionally, for a tensor $\mathcal{Y}$ of the same size as $\mathcal{X}$, we define the following operators acting through the $k$th unfolding
\begin{equation}
\mathcal{Q}_k\mathcal{Y}:\ Y_k\mapsto P_{C_k}\,Y_k,
\qquad
\Psi_k\mathcal{Y}:\ Y_k\mapsto X_kR_k^\dagger\;Y_k(\mathcal{I}_{\le k},:),
\qquad
\mathcal{Q}_{d-1}'\mathcal{Y}:\ Y_{d-1}\mapsto Y_{d-1}\,R_{d-1}^\dagger R_{d-1}.
\end{equation}
Each $\mathcal{Q}_k$ (and $\mathcal{Q}_{d-1}'$) is an orthogonal projection applied to one side of one unfolding. If $R_k$ has full row rank, then $\Psi_k$ is an (oblique) projection, since $\bigl(X_kR_k^\dagger Y_k(\mathcal{I}_{\le k},:)\bigr)(\mathcal{I}_{\le k},:)=R_kR_k^\dagger\,Y_k(\mathcal{I}_{\le k},:)=Y_k(\mathcal{I}_{\le k},:)$, with
\begin{equation}\label{eq:gamma_norm}
\|\Psi_k\|_{2}\;=\;\|X_kR_k^\dagger\|_2\;=\;\gamma_k\;\ge\;1 ,
\end{equation}
where $\gamma_k\ge1$ because $X_kR_k^\dagger$ contains $I_r$ in the rows $\mathcal{I}_{\le k}$.

\begin{lemma}[Factorization of the PEID format]\label{lem:factorization_app}
Assume each $R_k$, $1\le k\le d-1$, has full row rank. Then the TT tensor $\tilde{\mathcal{X}}$ defined by \cref{eq:ttcross_cores,eq:ttcross_PEID_build} satisfies
\begin{equation}\label{eq:factorization_id}
\tilde{\mathcal{X}} \;=\; \mathcal{Q}_1\Psi_1\;\mathcal{Q}_2\Psi_2\cdots\mathcal{Q}_{d-2}\Psi_{d-2}\;\mathcal{Q}_{d-1}\mathcal{Q}_{d-1}'\;\mathcal{X} .
\end{equation}
No nestedness of the pivot sets is required.
\end{lemma}

\begin{proof}
Since $R_k$ has full row rank, $\Psi_k$ is an (oblique) projection. We use induction on $d$. For $d=2$, $\tilde{\mathcal{X}}=C_1Y_1R_1=P_{C_1}X_1R_1^\dagger R_1=\mathcal{Q}_1\mathcal{Q}_1'\mathcal{X}$.

Let $d\ge3$. Contracting the cores \cref{eq:ttcross_cores} and grouping, we get
\[
\tilde X_1 \;=\; C_1\,Y_1\,B,\qquad
B\bigl(:,(i_2,\dots,i_d)\bigr)= M_2(:,i_2,:)\,Y_2\,M_3(:,i_3,:)\cdots Y_{d-1}\,R_{d-1}(:,i_d)\in\R^{r}.
\]
Consider the tensor obtained by ${\rm reshape}\left(\mathcal{X},n^2,n,\dots,n\right)$, so its unfoldings coincide with $X_2,\dots,X_{d-1}$. Therefore, the format \cref{eq:ttcross_cores,eq:ttcross_PEID_build} built on separations $2,\dots,d-1$ is well defined for it, and by the induction hypothesis (applied to this order-$(d-1)$ tensor), we have
\[
\mathcal{V} \;=\; \mathcal{Q}_2\Psi_2\cdots\mathcal{Q}_{d-2}\Psi_{d-2}\;\mathcal{Q}_{d-1}\mathcal{Q}_{d-1}'\;\mathcal{X} ,
\]
with first core $C_2$. More explicitly, $\mathcal{V}\bigl(\beta,i_3,\dots,i_d\bigr)=C_2(\beta,:)\,Y_2\,M_3(:,i_3,:)\cdots Y_{d-1}R_{d-1}(:,i_d)$ for $\beta\in[n^2]$. By the slice identity \cref{eq:slice_identity} with $k=2$, we can represent $B$ from above as
\[
B\bigl(:,(i_2,\dots,i_d)\bigr)
= C_2\bigl(\mathcal{I}_{\le1}\otimes\{i_2\},:\bigr)\,Y_2\cdots
= \mathcal{V}\bigl(\mathcal{I}_{\le1}\otimes\{i_2\},\,i_3,\dots,i_d\bigr),
\]
i.e., $B=V_1(\mathcal{I}_{\le1},:)$, the pivot rows of the first (original, unmerged) unfolding of $\mathcal{V}$. Hence
\[
\tilde X_1 = C_1\,C_1^\dagger X_1R_1^\dagger\;V_1(\mathcal{I}_{\le1},:) = \bigl(\mathcal{Q}_1\Psi_1\,\mathcal{V}\bigr)_1 ,
\]
which is \cref{eq:factorization_id}.
\end{proof}

\begin{remark}[Where the obliqueness lives]\label{rem:obliqueness_app}
\Cref{lem:factorization_app} makes precise a structural fact: in the TT chaining, the row factors $R_k^\dagger$ of the intermediate cores $Y_k$ do \emph{not} act as orthogonal row projections; they are consumed as the oblique gluing operators $\Psi_k$ that route information into the next core through the $r$ pivot rows. Only the outermost row factor survives as a genuine orthogonal projection $\mathcal{Q}_{d-1}'$. Thus projection enhancement removes obliqueness \emph{within} each separation (the $\mathcal{Q}_k$ have unit norm), but obliqueness \emph{between} separations is intrinsic to gluing rank-$r$ cores through skeleton rows. This is also why the algorithms of \cref{sec:algorithm} are naturally ``one-sided'': the intermediate right-projections would not buy orthogonality anyway.
\end{remark}

\begin{theorem}[Master bound for the PEID format]\label{thm:format_master}
Assume each $R_k$ has full row rank, and set
\[
\epsilon_k^{\mathrm{col}}:=\|(I-P_{C_k})\,X_k\|_F,\qquad
\epsilon_k^{\mathrm{row}}:=\|X_k\,(I-R_k^\dagger R_k)\|_F,\qquad 1\le k\le d-1 .
\]
Then the format \cref{eq:ttcross_cores,eq:ttcross_PEID_build} satisfies
\begin{equation}\label{eq:format_master_bound}
\|\mathcal{X}-\tilde{\mathcal{X}}\|_F\;\le\;\sum_{k=1}^{d-1}\Bigl(\prod_{j=1}^{k-1}\gamma_j\Bigr)\bigl(\epsilon_k^{\mathrm{col}}+\epsilon_k^{\mathrm{row}}\bigr),
\end{equation}
where $\gamma_j$ is computed as in~\cref{eq:gamma_norm}.
\end{theorem}

\begin{proof}
Write \cref{eq:factorization_id} as $\tilde{\mathcal{X}}=O_1O_2\cdots O_{2d-2}\,\mathcal{X}$ with the ordered operator list
\[
(O_1,\dots,O_{2d-2})=(\mathcal{Q}_1,\Psi_1,\mathcal{Q}_2,\Psi_2,\dots,\mathcal{Q}_{d-2},\Psi_{d-2},\mathcal{Q}_{d-1},\mathcal{Q}_{d-1}') .
\]
The telescoping identity $I-O_1\cdots O_m=\sum_{j=1}^m O_1\cdots O_{j-1}(I-O_j)$ and submultiplicativity give
\[
\|\mathcal{X}-\tilde{\mathcal{X}}\|_F\le\sum_{j=1}^{2d-2}\Bigl(\prod_{i<j}\|O_i\|_{2}\Bigr)\,\|(I-O_j)\mathcal{X}\|_F.
\]
Orthogonal projections applied to one side of one unfolding are nonexpansive in the Frobenius norm, so $\|\mathcal{Q}_k\|_2=\|\mathcal{Q}_{d-1}'\|_2=1$, while $\|\Psi_j\|_2=\gamma_j$ by \cref{eq:gamma_norm}; hence the prefix preceding either operator at separation $k$ is at most $\gamma_1\cdots\gamma_{k-1}$. Finally,
\[
(I-\mathcal{Q}_k)\mathcal{X}:\ X_k-P_{C_k}X_k,\qquad
(I-\Psi_k)\mathcal{X}:\ X_k-X_kR_k^\dagger R_k,\qquad
(I-\mathcal{Q}_{d-1}')\mathcal{X}:\ X_{d-1}-X_{d-1}R_{d-1}^\dagger R_{d-1},
\]
whose Frobenius norms are $\epsilon_k^{\mathrm{col}}$, $\epsilon_k^{\mathrm{row}}$, and $\epsilon_{d-1}^{\mathrm{row}}$ respectively; collecting terms yields \cref{eq:format_master_bound}.
\end{proof}

\begin{remark}[Reading the amplification factors]\label{rem:gamma_reading}
Write $\epsilon_k:=\epsilon_k^{\mathrm{col}}+\epsilon_k^{\mathrm{row}}$ for the local error at separation $k$. The bound~\cref{eq:format_master_bound} weights $\epsilon_k$ by the product of the factors $\gamma_j$ with $j<k$, and three features of this weighting are worth making explicit. First, the product is empty for $k=1$ and the factor $\gamma_{d-1}$ never occurs. Both facts reflect the factorization~\cref{eq:factorization_id}: the error committed at separation $k$ is transported to the full tensor through the gluing operators $\Psi_1,\dots,\Psi_{k-1}$ of the preceding separations, and the last separation is closed by the orthogonal projection $\mathcal{Q}_{d-1}'$ rather than by an operator $\Psi_{d-1}$. Second, for $d=2$ there is a single separation and the bound reads $\|\mathcal{X}-\tilde{\mathcal{X}}\|_F\le\epsilon_1$, which is the matrix estimate for $P_CAP_R$. This does not mean that $\gamma_1=1$ for matrices. The factor $\gamma_1=\|X_1R_1^\dagger\|_2$ is in general larger than one, and it is absent only because no further core is glued behind the single separation. For $d=3$ the bound reads $\epsilon_1+\gamma_1\epsilon_2$, so $\gamma_1$ enters as soon as a second separation has to be reached through the pivot rows of the first. Third, $\epsilon_k$ and $\gamma_k$ measure different properties of the skeleton. The local error $\epsilon_k$ depends only on whether the selected fibers can span the dominant subspaces of $X_k$. The factor $\gamma_k$ depends on how well-conditioned the pivot rows are as a basis for extending information to all rows of $X_k$. A skeleton can therefore be accurate and still poorly conditioned. If $X_k$ has exact rank $r$, then $\epsilon_k=0$ for every choice of pivots with nonsingular intersection, whereas $\gamma_k=\|[U_k(\mathcal{I}_{\le k},:)]^{-1}\|_2$ can be arbitrarily large (see~\cref{rem:gamma_eta}). In the bound, a large $\gamma_k$ is harmless as long as the subsequent local errors vanish, and it becomes visible exactly when they do not.
\end{remark}

The per-separation quantities in \cref{thm:format_master} are pure column/row subset selection errors of the full unfoldings, with three consequences. First, they are \emph{optimal} for the chosen skeletons: for any $W$,
\begin{equation}\label{eq:optimality_app}
\epsilon_k^{\mathrm{col}}\le\|X_k-C_kW\|_F,\qquad
\epsilon_k^{\mathrm{row}}\le\|X_k-WR_k\|_F,
\end{equation}
so in particular each is bounded by the Frobenius error of the cross approximation built from the same pivots---a format-level ``projection enhancement never hurts'' statement. Second, they are bounded below by $\Delta_k=\|X_k-(X_k)_r\|_F$, since the projections have rank at most $r$. Third, they are monotone under enlarging the fiber sets, so replacing $\mathcal{J}_{>k}$ by $\mathcal{J}_{>k}\cup\mathcal{L}_{>k}$ and $\mathcal{I}_{\le k}$ by $\mathcal{I}_{\le k}\cup\mathcal{K}_{\le k}$, as in \cref{eq:matrix_cur_over}, can only decrease them.

\subsection{Imported matrix results, extended to rectangular matrices}\label{app:imported}
We use the error bound of~\cite{cortinovis2020maximum} for cross approximation with complete pivoting (equivalently, LU with complete pivoting) \cite[Thm.~6]{cortinovis2020maximum}: if $A\in\R^{N\times N}$ has rank at least $r$, then the index sets $(\mathcal{I},\mathcal{J})$ produced by $r$ steps of complete pivoting satisfy
\begin{equation}\label{eq:CKM_bound}
\|A-A(:,\mathcal{J})A(\mathcal{I},\mathcal{J})^{-1}A(\mathcal{I},:)\|_{\max}\;\le\;4^{r}\,\rho_r\,\sigma_{r+1}(A),
\end{equation}
where $\rho_r$ is the growth factor of Gaussian elimination with complete pivoting, for which Wilkinson's classical bound gives $\rho_r\le 2\sqrt{r+1}\,(r+1)^{\ln(r+1)/4}$ \cite{wilkinson1961error,cortinovis2020maximum}. Here and below, the \emph{remainder} $Z_k$ after $k$ steps of ACA is defined by the recursion of \cite[Alg.~1]{cortinovis2020maximum}: $Z_0:=A$ and $Z_{k+1}:=Z_k-p_{k+1}^{-1}Z_k(:,j_{k+1})Z_k(i_{k+1},:)$, where $(i_{k+1},j_{k+1})$ is the entry of $Z_k$ of maximum modulus and $p_{k+1}:=Z_k(i_{k+1},j_{k+1})$. After $r$ steps the remainder is exactly the cross-approximation residual, $Z_r=A-A(:,\mathcal{J})A(\mathcal{I},\mathcal{J})^{-1}A(\mathcal{I},:)$, a full-size matrix that vanishes on the pivot rows and columns; see \cite[Eq.~(5) and Rem.~5]{cortinovis2020maximum}. Tensor unfoldings are often rectangular, so we record that \cref{eq:CKM_bound} still holds.

\begin{lemma}[Rectangular version of \cref{eq:CKM_bound}]\label{lem:rectCKM_app}
Let $A\in\R^{M\times N}$ have rank at least $r<\min\{M,N\}$. Then the index sets $(\mathcal{I},\mathcal{J})$ produced by $r$ steps of ACA with complete pivoting applied to $A$ satisfy \cref{eq:CKM_bound}.
\end{lemma}

\begin{proof}
Without loss of generality, assume $M\le N$. Otherwise, we work on $A^T$. Pad $\bar A:=\begin{bmatrix}A\\ 0\end{bmatrix}\in\R^{N\times N}$. We claim that complete pivoting applied to $\bar A$ selects the same pivots as for $A$ and that the remainders satisfy $\bar Z_k=\begin{bmatrix}Z_k\\ 0\end{bmatrix}$ for all $k\le r$. This holds for $k=0$; inductively, the maximizing entry of $\bar Z_k$ lies in the top block (the bottom block is zero, and $Z_k\ne0$ because $\rank(A)\ge r$), so the selected pivot coincides with that of $A$, and the rank-one update leaves the zero block intact. Consequently the remainder after $r$ steps satisfies $\|\bar Z_r\|_{\max}=\|Z_r\|_{\max}$, while $\sigma_{r+1}(\bar A)=\sigma_{r+1}(A)$. Applying \cite[Thm.~6]{cortinovis2020maximum} to $\bar A$ gives the claim, and the growth factor of the padded instance can be bounded by the same $\rho_r$.
\end{proof}

\subsection{Proof of \texorpdfstring{\cref{thm:PEID_format_main}}{Theorem}}
\begin{proof}[Proof of {\cref{thm:PEID_format_main}}]
By \cref{eq:optimality_app} with $W=\hat A_k^{-1}R_k$ (resp.\ $W=C_k\hat A_k^{-1}$), where $\hat A_k=X_k(\mathcal{I}_{\le k},\mathcal{J}_{>k})$, both $\epsilon_k^{\mathrm{col}}$ and $\epsilon_k^{\mathrm{row}}$ are bounded by the Frobenius norm of the ACA remainder of $X_k$, hence by $\sqrt{n^k\cdot n^{d-k}}=n^{d/2}$ times its Chebyshev norm, and \cref{lem:rectCKM_app} bounds the latter by $4^r\rho_r\sigma_{r+1}(X_k)$. Insert into \cref{eq:format_master_bound}.
\end{proof}

\begin{remark}\label{rem:comlete_replacement}
For the factor $4^r\rho_r n^{d/2}\sigma_{r+1}(X_k)$, a probabilistic refinement replacing it by a constant multiple of the best rank-$r$ unfolding error is given in~\cref{app:alg3_analysis}.
\end{remark}

\begin{remark}[Controlling $\gamma_k$; consistency with~\cref{app:alg3_analysis}]\label{rem:gamma_eta}
The amplification factors $\gamma_k=\|X_kR_k^\dagger\|_2$ are the row-ID coefficient norms of the pivot rows, the classical measure of skeleton quality in the ID/CUR literature. Three regimes: (i) if the rows are selected by strong RRQR or maximum-volume criteria, the coefficient matrix has entries bounded by a modest constant $f$, giving bounds of the form $\sqrt{1+f^2r(n^k-r)}$ \cite{gu1996efficient,cheng2005compression,goreinov2010find}; (ii) $\gamma_k$ is small precisely when the pivots are good row skeletons, which is what the stability experiments of \cref{fig:examples24} probe; (iii) if $X_k$ has exact rank $r$ with orthonormal column basis $U_k$, then $X_kR_k^\dagger=U_k[U_k(\mathcal{I}_{\le k},:)]^{-1}$, so $\gamma_k=\|[U_k(\mathcal{I}_{\le k},:)]^{-1}\|_2$, which is precisely the quantity $\eta_k$ of~\cref{app:alg3_analysis} with row set $\mathcal{I}_{\le k}$ and no row oversampling. Consequently,~\cref{app:alg3_analysis} is the analysis of the computable, sampled surrogate of the format \cref{eq:ttcross_cores,eq:ttcross_PEID_build}, in which full unfoldings are replaced by sampled ones and $\gamma_k$ by the certifiable $\eta_k$; and row oversampling in the format can only decrease each $\gamma_k$ (\cref{lem:rowmono_app}).
\end{remark}

\section{Error analysis of TT-PEID-Par (\texorpdfstring{\cref{alg:TT_update_parallel}}{Algorithm 3})}\label{app:alg3_analysis}
This appendix proves \cref{thm:alg3_main}, via a deterministic master bound (\cref{thm:master_app}) in which each algorithmic ingredient contributes a separate term, followed by probabilistic bounds for the oversampling contributions.

\subsection{Setup}\label{app:alg3_setup}
Let $(\mathcal{I}_{\le k},\mathcal{J}_{>k})$ be nested skeleton pivots of cardinality $r$ as in \cref{eq:tt_nestedness}, and let $\mathcal{K}_{\le k}$ and $\mathcal{L}_{>k}$ be the oversampling sets of cardinality $p$ produced by \cref{alg:TT_update_parallel}. For notational simplicity, we define
\[
\mathcal{S}_k:=\mathcal{I}_{\le k}\cup\mathcal{K}_{\le k},\qquad \mathcal{T}^{\mathrm{c}}_k:=\mathcal{J}_{>k}\cup\mathcal{L}_{>k},\qquad |\mathcal{S}_k|=|\mathcal{T}^{\mathrm{c}}_k|=r+p,
\]
with the convention $\mathcal{S}_0:=\{1\}$, and the superscript on $\mathcal{T}^{\mathrm{c}}_k$ distinguishes this column index set from the TT cores $\mathcal{T}_j$. Since both $\mathcal{I}$- and $\mathcal{K}$-sets are nested, we have
\begin{equation}\label{eq:S_nested}
\mathcal{S}_k\subseteq\mathcal{S}_{k-1}\otimes\mathbb{I}_k,\qquad 1\le k\le d-1 .
\end{equation}
In addition, we set the restricted unfoldings and sampled column blocks used in line~6 of \cref{alg:TT_update_parallel} as
\[
A_k:=X_k(\mathcal{S}_{k-1}\otimes\mathbb{I}_k,\,:)\in\R^{m_k\times N_k},\qquad
C_k^{\mathrm{s}}:=A_k(:,\mathcal{T}^{\mathrm{c}}_k)\in\R^{m_k\times(r+p)},
\]
with $m_1=n$, $m_k=(r+p)n$ for $k\ge2$, and $N_k=n^{d-k}$. For the analysis we fix the following construction of the orthonormal factors and cores (lines~6 and 9--10 of \cref{alg:TT_update_parallel}): $U_k\in\R^{m_k\times r}$ consists of the $r$ dominant left singular vectors of $C_k^{\mathrm{s}}$; the cores are $\mathcal{T}_1=U_1$, $\mathcal{T}_k=\mathrm{reshape}\bigl([U_{k-1}(\mathcal{S}_{k-1},:)]^\dagger\tilde U_k,\,r,n,r\bigr)$ for $2\le k\le d-1$ with $\tilde U_k=\mathrm{reshape}(U_k,\,r+p,\,nr)$, and $\mathcal{T}_d=[U_{d-1}(\mathcal{S}_{d-1},:)]^\dagger X_{d-1}(\mathcal{S}_{d-1},:)$. We assume each $U_k(\mathcal{S}_k,:)\in\R^{(r+p)\times r}$ has full column rank, and define
\begin{equation}\label{eq:eta_def}
\eta_k:=\bigl\|[U_k(\mathcal{S}_k,:)]^\dagger\bigr\|_2=\frac{1}{\sigma_r\bigl(U_k(\mathcal{S}_k,:)\bigr)}\;\ge\;1,
\qquad
\theta_k:=\bigl\|(C_k^{\mathrm{s}})^\dagger A_k\bigr\|_2 ,
\end{equation}
where $\eta_k\ge1$ by \cref{eq:interlacing_app} since $U_k$ has orthonormal columns. Otherwise, if $U_k(\mathcal{S}_k,:)$ does not have full column rank, $\eta_k=\infty$.

The output of \cref{alg:TT_update_parallel} admits the following recursive description. Write $X^{(1)}:=X_1=A_1$, and for $2\le k\le d$, let $X^{(k)}:={\rm reshape}\left(X_{k-1}(\mathcal{S}_{k-1},:),r+p,n,\dots,n\right)$ with $d-k+1$ modes. Then by \cref{eq:S_nested}, the second unfolding of $X^{(k)}$ is exactly $A_k$, and the rows $\mathcal{S}_k$ of $A_k$ reproduce $X^{(k+1)}$. Contracting the cores shows that the computed approximation satisfies
\begin{equation}\label{eq:recursion_app}
\hat X^{(k)} = U_k\,[U_k(\mathcal{S}_k,:)]^\dagger\,\hat X^{(k+1)},\qquad 1\le k\le d-1,
\qquad \hat X^{(d)}:=X^{(d)},
\end{equation}
where each $\hat X^{(k)}$ is understood as a matrix of the same size as $A_k$, and $\tilde{\mathcal{X}}:=\hat X^{(1)}$ is the TT tensor returned by \cref{alg:TT_update_parallel}. The last level is exact because $\mathcal{T}_d$ uses the rows $X_{d-1}(\mathcal{S}_{d-1},:)$ of the data itself. Nestedness \cref{eq:S_nested} is used exactly once, and structurally: it makes $\mathcal{S}_k$ a valid row set of $A_k$, so that \cref{eq:recursion_app} closes. This is also a theoretical explanation of why \cref{alg:TT_update_parallel} must restrict its oversampling sets to be nested, while \cref{alg:TT_update_sequential}, which frees $\mathcal{K}_{\le k}$, must pay for it by rebuilding the row map $F_{\le k}$ from the computed cores.

\subsection{Deterministic master bound}

\begin{lemma}[Oblique projection step]\label{lem:oblique_app}
Let $U\in\R^{m\times r}$ have orthonormal columns and let $\mathcal{S}$ be a row set with $U(\mathcal{S},:)$ of full column rank. Then $\Pi:=U\,[U(\mathcal{S},:)]^\dagger E_\mathcal{S}^T$, where $E_\mathcal{S}^TA=A(\mathcal{S},:)$, is a projection onto $\mathrm{range}(U)$ with $\|\Pi\|_2=\|[U(\mathcal{S},:)]^\dagger\|_2=:\eta$, and for every matrix $A$ with $m$ rows,
\[
\|(I-\Pi)A\|_F\;\le\;\eta\,\|(I-UU^T)A\|_F.
\]
\end{lemma}

\begin{proof}
$\Pi U=U[U(\mathcal{S},:)]^\dagger U(\mathcal{S},:)=U$, so $\Pi^2=\Pi$ and $(I-\Pi)U=0$; hence $(I-\Pi)A=(I-\Pi)(I-UU^T)A$ and $\|(I-\Pi)A\|_F\le\|I-\Pi\|_2\|(I-UU^T)A\|_F$. Since $U$ and $E_\mathcal{S}$ have orthonormal columns, $\|\Pi\|_2=\|[U(\mathcal{S},:)]^\dagger\|_2=\eta$, and for any projection $\Pi\notin\{0,I\}$ one has $\|I-\Pi\|_2=\|\Pi\|_2$ \cite{szyld2006many}.
\end{proof}

\begin{theorem}[Master bound for \cref{alg:TT_update_parallel}]\label{thm:master_app}
With the conventions of~\cref{app:alg3_setup} and $e_k:=\|(I-U_kU_k^T)A_k\|_F$, $1\le k\le d-1$, the output $\tilde{\mathcal{X}}$ of \cref{alg:TT_update_parallel} satisfies
\begin{equation}\label{eq:master_app}
\|\mathcal{X}-\tilde{\mathcal{X}}\|_F\;\le\;\sum_{k=1}^{d-1}\Bigl(\prod_{j=1}^{k}\eta_j\Bigr)\,e_k .
\end{equation}
\end{theorem}

\begin{proof}
Let $\mathrm{Err}_k:=\|X^{(k)}-\hat X^{(k)}\|_F$, so $\mathrm{Err}_1=\|\mathcal{X}-\tilde{\mathcal{X}}\|_F$ and $\mathrm{Err}_d=0$. Fix $k\le d-1$ and write $\Pi_k:=U_k[U_k(\mathcal{S}_k,:)]^\dagger E_{\mathcal{S}_k}^T$. Using \cref{eq:recursion_app} and $A_k(\mathcal{S}_k,:)=X^{(k+1)}$,
\[
X^{(k)}-\hat X^{(k)}
=\bigl(A_k-\Pi_kA_k\bigr)
\;+\;U_k[U_k(\mathcal{S}_k,:)]^\dagger\bigl(X^{(k+1)}-\hat X^{(k+1)}\bigr).
\]
\Cref{lem:oblique_app} bounds the first term by $\eta_k e_k$, and the second by $\eta_k\,\mathrm{Err}_{k+1}$ because $\|U_k[U_k(\mathcal{S}_k,:)]^\dagger\|_2=\eta_k$. Hence $\mathrm{Err}_k\le\eta_k e_k+\eta_k\,\mathrm{Err}_{k+1}$, and unrolling the recursion from $k=1$ with $\mathrm{Err}_d=0$ gives \cref{eq:master_app}.
\end{proof}

All suboptimality caused by the oblique projections is thus confined to the computable factors $\eta_j$, while $e_k$ is a pure column subset selection error for the sampled columns $\mathcal{T}^{\mathrm{c}}_k$ of the restricted unfolding, to which the skeleton and the oversampling contribute. We now make $e_k$ explicit.

\subsection{Bounding \texorpdfstring{$e_k$}{ek}: truncation, skeleton, and oversampling}

\begin{lemma}[Truncation to rank $r$]\label{lem:trunc_app}
With $U_k$ the $r$ dominant left singular vectors of $C_k^{\mathrm{s}}$,
\[
e_k\;\le\;\|(I-P_{C_k^{\mathrm{s}}})A_k\|_F\;+\;\theta_k\,\|C_k^{\mathrm{s}}-(C_k^{\mathrm{s}})_r\|_F
\;\le\;\|(I-P_{C_k^{\mathrm{s}}})A_k\|_F\;+\;\theta_k\,\Delta_k .
\]
\end{lemma}

\begin{proof}
Split $A_k=(I-P_{C_k^{\mathrm{s}}})A_k+P_{C_k^{\mathrm{s}}}A_k$. Since $\|I-U_kU_k^T\|_2\le1$, the first part contributes at most $\|(I-P_{C_k^{\mathrm{s}}})A_k\|_F$. For the second part, $U_kU_k^TC_k^{\mathrm{s}}=(C_k^{\mathrm{s}})_r$, so
$(I-U_kU_k^T)P_{C_k^{\mathrm{s}}}A_k=\bigl(C_k^{\mathrm{s}}-(C_k^{\mathrm{s}})_r\bigr)(C_k^{\mathrm{s}})^\dagger A_k$,
whose Frobenius norm is at most $\|C_k^{\mathrm{s}}-(C_k^{\mathrm{s}})_r\|_F\,\|(C_k^{\mathrm{s}})^\dagger A_k\|_2$, using $\|MB\|_F\le\|M\|_F\|B\|_2$. Finally, $C_k^{\mathrm{s}}$ is a submatrix of $X_k$, so \cref{eq:interlacing_app} gives $\|C_k^{\mathrm{s}}-(C_k^{\mathrm{s}})_r\|_F\le\Delta_k$.
\end{proof}

For this bound, if TT core sizes $r+p$ are acceptable, one may take $U_k$ to be an orthonormal basis of $\mathrm{range}(C_k^{\mathrm{s}})$, in which case $e_k=\|(I-P_{C_k^{\mathrm{s}}})A_k\|_F$ exactly and the $\theta_k\Delta_k$ term disappears. In addition, if $U_k$ is computed by a randomized range finder instead of a truncated SVD, the standard additional factors of \cite[Thm.~10.5]{halko2011finding} appear.

The remaining quantity $\|(I-P_{C_k^{\mathrm{s}}})A_k\|_F$ is monotone in the column set: enlarging $\mathcal{T}^{\mathrm{c}}_k$ can only decrease it. Since $\mathcal{T}^{\mathrm{c}}_k\supseteq\mathcal{J}_{>k}$ and $\mathcal{T}^{\mathrm{c}}_k\supseteq\mathcal{L}_{>k}$, we obtain two independent upper bounds and may take their minimum: at the level of the bound, oversampling can only help the skeleton, and the skeleton can only help oversampling.

\begin{lemma}[Skeleton branch]\label{lem:skel_app}
Suppose $(\mathcal{I}_{\le k},\mathcal{J}_{>k})$ is produced by $r$ steps of ACA with complete pivoting on $X_k$ with $\rank(X_k)>r$. Then
\[
\|(I-P_{C_k^{\mathrm{s}}})A_k\|_F\;\le\;\sqrt{m_kN_k}\;4^{r}\rho_r\,\sigma_{r+1}(X_k),
\qquad m_kN_k\le (r+p)\,n^{\,d-k+1}.
\]
More generally, the left-hand side is bounded by the Frobenius norm of the ACA remainder of $X_k$ restricted to the rows $\mathcal{S}_{k-1}\otimes\mathbb{I}_k$.
\end{lemma}

\begin{proof}
By nestedness, $\mathcal{I}_{\le k}\subseteq\mathcal{I}_{\le k-1}\otimes\mathbb{I}_k\subseteq\mathcal{S}_{k-1}\otimes\mathbb{I}_k$, so the pivot rows of $X_k$ survive in $A_k$, and the blocks $A_k(:,\mathcal{J}_{>k})$, $\hat A_k=X_k(\mathcal{I}_{\le k},\mathcal{J}_{>k})$ are submatrices of $X_k$'s cross data. Since orthogonal projection onto the columns $\mathcal{T}^{\mathrm{c}}_k\supseteq\mathcal{J}_{>k}$ is at least as accurate as any approximation of the form $A_k(:,\mathcal{J}_{>k})W$,
\[
\|(I-P_{C_k^{\mathrm{s}}})A_k\|_F\le\|A_k-A_k(:,\mathcal{J}_{>k})\,\hat A_k^{-1}X_k(\mathcal{I}_{\le k},:)\|_F.
\]
Every entry of the matrix on the right is an entry of the full ACA remainder $X_k-X_k(:,\mathcal{J}_{>k})\hat A_k^{-1}X_k(\mathcal{I}_{\le k},:)$ (restricting rows commutes with forming the cross), so its Chebyshev norm is bounded by $4^r\rho_r\sigma_{r+1}(X_k)$ by \cref{lem:rectCKM_app}, and the Frobenius norm by $\sqrt{m_kN_k}$ times that.
\end{proof}

\subsection{Oversampling branches: columns and rows}\label{app:oversampling}
Both oversampling sets of \cref{alg:TT_update_parallel} are drawn uniformly at random---$\mathcal{L}_{>k}$ among the columns of $A_k$ and $\mathcal{K}_{\le k}$ among its rows---and uniform sampling can capture a subspace only if no single row of its orthonormal basis carries a disproportionate share of the mass. The quantity that measures this is the \emph{coherence}: for an $r$-dimensional subspace of $\R^{N}$ with orthonormal basis $V\in\R^{N\times r}$,
\[
\mu(V):=\frac{N}{r}\max_{1\le i\le N}\|V(i,:)\|_2^2\in\Bigl[1,\frac{N}{r}\Bigr],
\]
which is a property of the subspace, not of the chosen basis. The value $\mu(V)=1$ means the mass is spread evenly over all coordinates, while $\mu(V)=N/r$ means the subspace contains a canonical basis vector, so that any uniform sample missing the corresponding coordinate loses that direction entirely. Coherence is the standard hypothesis under which uniform sampling is analysed in the skeleton/CUR literature \cite{chiu2013sublinear,cortinovis2024sublinear}. Since columns are sampled from $\R^{N_k}$, the column branch needs incoherence of the dominant \emph{right} singular subspace of $X_k$. In addition, since rows are sampled from $\R^{m_k}$, the row branch needs incoherence of the dominant \emph{left} singular subspace of $A_k$, which is spanned by the computed basis $U_k$. We state the two hypotheses separately, because they enter the analysis differently and the second one is checkable at run time.
 
\begin{assumption}\label{ass:incoh_app}
For each $1\le k\le d-1$, write $X_k=B_k+E_k$ with $B_k=(X_k)_r$ (so $\|E_k\|_F=\Delta_k$), let $V_k\in\R^{N_k\times r}$ be an orthonormal basis of the row space of $B_k$, and assume: (i) (incoherence) $\mu(V_k)\le\mu$ and $\mu r^2\le N_k/2$; (ii) (restriction regularity) $\rank\bigl(B_k(\mathcal{S}_{k-1}\otimes\mathbb{I}_k,:)\bigr)=r$, i.e., the retained rows do not lose the dominant row space. (Condition (ii) holds whenever the retained rows contain $r$ linearly independent rows of $B_k$; since $\mathcal{I}_{\le k}\subseteq\mathcal{S}_{k-1}\otimes\mathbb{I}_k$ and the skeleton block $\hat A_k$ is nonsingular, it is implied by the skeleton pivots capturing the dominant subspace, and it is generic.)
\end{assumption}
 
Under \cref{ass:incoh_app}(ii) the restricted unfolding inherits the split $A_k=B_k'+E_k'$ with $B_k':=B_k(\mathcal{S}_{k-1}\otimes\mathbb{I}_k,:)$ of rank $r$, $\|E_k'\|_F\le\Delta_k$, and---crucially---the row space of $B_k'$ equals that of $B_k$, so its right singular subspace has the same coherence.
 
\begin{lemma}[Column oversampling branch]\label{lem:oversampling_app}
Let \cref{ass:incoh_app} hold for index $k$, and let $\mathcal{L}_{>k}$ of size $p$ be drawn uniformly without replacement from $[N_k]\setminus\mathcal{J}_{>k}$. There is an absolute constant $c_0$ such that if $p\ge c_0\,\mu\,r\,\log(r/\delta)$, then with probability at least $1-2\delta$,
\[
\|(I-P_{C_k^{\mathrm{s}}})A_k\|_F\;\le\;\bigl(1+c_1\,\delta^{-1/2}\bigr)\,\Delta_k,
\]
with an absolute constant $c_1\le2\sqrt2$.
\end{lemma}
 
\begin{proof}
By monotonicity it suffices to bound $\|(I-P)A_k\|_F$ for $P$ the orthogonal projection onto $\mathrm{range}\bigl(A_k(:,\mathcal{L}_{>k})\bigr)$. Write $A_k=B_k'+E_k'$ as above and let $V\in\R^{N_k\times r}$ be an orthonormal basis of the row space of $B_k'$, with coherence at most $\mu$.
 
\emph{Step 1 (well-conditioned sampled rows of $V$).} The sampling is uniform over $\Omega:=[N_k]\setminus\mathcal{J}_{>k}$. Since $\|V(\mathcal{J}_{>k},:)^TV(\mathcal{J}_{>k},:)\|_2\le r\cdot\mu r/N_k\le1/2$ by \cref{ass:incoh_app}(i), we have $\sum_{i\in\Omega}V(i,:)^TV(i,:)\succeq\tfrac12 I_r$. A matrix Chernoff bound for uniform sampling without replacement \cite{tropp2012user,gross2010note} then yields an absolute constant $c_0$ such that, for $p\ge c_0\mu r\log(r/\delta)$, with probability at least $1-\delta$,
\begin{equation}\label{eq:chernoff_app}
\sigma_r\bigl(V(\mathcal{L}_{>k},:)\bigr)\;\ge\;\tfrac12\sqrt{p/N_k},
\qquad\text{hence}\qquad
\bigl\|\bigl(V(\mathcal{L}_{>k},:)^T\bigr)^\dagger\bigr\|_2\le 2\sqrt{N_k/p}.
\end{equation}
 
\emph{Step 2 (structural identity).} On the event \cref{eq:chernoff_app}, $V(\mathcal{L}_{>k},:)^T\in\R^{r\times p}$ has full row rank, so with $B_k'=W\Sigma V^T$,
\[
B_k'(:,\mathcal{L}_{>k})\,\bigl(V(\mathcal{L}_{>k},:)^T\bigr)^\dagger V^T
= W\Sigma\, V(\mathcal{L}_{>k},:)^T\bigl(V(\mathcal{L}_{>k},:)^T\bigr)^\dagger V^T
= W\Sigma V^T=B_k'.
\]
Therefore the surrogate $\hat A:=A_k(:,\mathcal{L}_{>k})\bigl(V(\mathcal{L}_{>k},:)^T\bigr)^\dagger V^T$, which lies in the column space of $A_k(:,\mathcal{L}_{>k})$, satisfies
$A_k-\hat A = E_k' - E_k'(:,\mathcal{L}_{>k})\bigl(V(\mathcal{L}_{>k},:)^T\bigr)^\dagger V^T$,
and consequently
\[
\|(I-P)A_k\|_F\le\|A_k-\hat A\|_F
\le\|E_k'\|_F+2\sqrt{N_k/p}\;\|E_k'(:,\mathcal{L}_{>k})\|_F.
\]
 
\emph{Step 3 (size of the sampled residual).} Uniform sampling of $p$ out of $|\Omega|\ge N_k/2$ columns gives $\mathbb{E}\,\|E_k'(:,\mathcal{L}_{>k})\|_F^2=\frac{p}{|\Omega|}\|E_k'(:,\Omega)\|_F^2\le\frac{2p}{N_k}\Delta_k^2$, so by Markov's inequality, with probability at least $1-\delta$, $\|E_k'(:,\mathcal{L}_{>k})\|_F\le\sqrt{2p/(N_k\delta)}\,\Delta_k$. Intersecting the two events and combining the displays proves the claim with $c_1=2\sqrt2$.
\end{proof}
 
The row sets enter \cref{eq:master_app} only through the oblique factors $\eta_k=1/\sigma_r(U_k(\mathcal{S}_k,:))$, which are the tensor analogue of the (Q)DEIM/GappyPOD selection constant $\|(P^TV)^{-1}\|_2$ \cite{drmac2016new,peherstorfer2020stability} and of the oblique-projection factors in the oversampled CUR analysis of \cite{park2025accuracy}. They are monotone in the row set, in the same way that the column error is monotone in the column set.
 
\begin{lemma}[Row monotonicity]\label{lem:rowmono_app}
Let $U\in\R^{m\times r}$ have orthonormal columns and let $\mathcal{S}\subseteq\mathcal{S}'$ be row sets. Then $\sigma_r\bigl(U(\mathcal{S}',:)\bigr)\ge\sigma_r\bigl(U(\mathcal{S},:)\bigr)$, hence $\|[U(\mathcal{S}',:)]^\dagger\|_2\le\|[U(\mathcal{S},:)]^\dagger\|_2$. In particular, every row added to $\mathcal{S}_k$---whether a skeleton pivot or an oversample---can only decrease $\eta_k$ and hence improve \cref{eq:master_app}, on every realization.
\end{lemma}
 
\begin{proof}
$U(\mathcal{S}',:)^TU(\mathcal{S}',:)=U(\mathcal{S},:)^TU(\mathcal{S},:)+\sum_{i\in\mathcal{S}'\setminus\mathcal{S}}U(i,:)^TU(i,:)\succeq U(\mathcal{S},:)^TU(\mathcal{S},:)$, and $\lambda_{\min}$ is monotone with respect to the Loewner order; take square roots.
\end{proof}
 
By \cref{lem:rowmono_app}, a bound on $\eta_k$ follows from a bound on $\sigma_r(U_k(\mathcal{K}_{\le k},:))$ alone, i.e., from the uniformly sampled rows, exactly as the column branch used only $\mathcal{L}_{>k}$. The uniform row sampling of an orthonormal matrix succeeds, however, only when that matrix is incoherent: if, for instance, the first $r$ rows of $U_k$ form the identity and all other rows vanish ($\mu(U_k)=m_k/r$), then every sample $\mathcal{K}_{\le k}$ avoiding those $r$ rows gives $\sigma_r(U_k(\mathcal{K}_{\le k},:))=0$, and the randomly sampled rows contribute nothing to \cref{eq:master_app}. This is why the following hypothesis is required: it is precisely the left-sided counterpart of \cref{ass:incoh_app}(i).
 
\begin{assumption}\label{ass:incohU_app}
For each $1\le k\le d-1$, the computed basis $U_k$ satisfies $\mu(U_k)\le\mu_U$ and $\mu_U r^2\le m_k/2$.
\end{assumption}
 
Two features of \cref{ass:incohU_app} deserve comment. First, it is stated on the \emph{computed} basis $U_k$ rather than on the dominant left singular subspace of $A_k$, because $U_k$ is the matrix whose rows the algorithm actually samples: when $\rank(A_k)=r$ the two subspaces coincide and the distinction disappears, but in general a truncated or randomized range finder returns a basis whose coherence may differ from that of the exact dominant subspace, and it is the former that governs \cref{lem:rowrandom_app}. Second, and as a consequence, the hypothesis is \emph{verifiable at run time}: $\mu(U_k)$ is determined by the row norms of $U_k$ and costs $\mathcal{O}(m_kr)$ to evaluate, negligible compared with the $\mathcal{O}(m_kr^2)$ already spent on $U_k$. Checking it turns the probabilistic guarantee of \cref{thm:alg3_main} into a certifiable one, which is not possible for \cref{ass:incoh_app}, a statement about the unknown data. Additionally, the side condition $\mu_U r^2\le m_k/2$ is mild: it only ensures that excluding the $r$ skeleton rows from the sampling pool (recall $\mathcal{K}_{\le k}\cap\mathcal{I}_{\le k}=\emptyset$) removes at most half of the mass of every direction in $\mathrm{range}(U_k)$. When the assumption fails, or when one prefers to avoid it, \cref{rem:alg3_remarks}(iii) describes deterministic selections that bound $\eta_k$ without any incoherence hypothesis.
 
\begin{lemma}[Row oversampling branch]\label{lem:rowrandom_app}
In the local row indexing of $A_k$, the set $\mathcal{K}_{\le k}$ drawn in line~3 of \cref{alg:TT_update_parallel} is uniform without replacement over $[m_k]\setminus\mathcal{I}_{\le k}$ (nestedness \cref{eq:S_nested} makes the candidate set $\mathcal{S}_{k-1}\otimes\mathbb{I}_k\setminus\mathcal{I}_{\le k}$ coincide with the complement of the skeleton rows of $A_k$), and it is independent of $U_k$ conditionally on $\mathcal{S}_{k-1}$, since $U_k$ is a function of $\mathcal{S}_{k-1}$ and the column draws $\mathcal{L}_{>k}$ only. Let \cref{ass:incohU_app} hold and let $p\ge c_0\,\mu_U\,r\,\log(r/\delta)$, with the absolute constant $c_0$ of \cref{lem:oversampling_app}. Then, with conditional probability at least $1-\delta$,
\[
\sigma_r\bigl(U_k(\mathcal{K}_{\le k},:)\bigr)\;\ge\;\tfrac12\sqrt{p/m_k},
\qquad\text{hence by \cref{lem:rowmono_app}}\qquad
\eta_k\;\le\;2\sqrt{m_k/p}\;\le\;2\sqrt{(r+p)\,n/p}\,.
\]
\end{lemma}
 
\begin{proof}
This is Step~1 of the proof of \cref{lem:oversampling_app} with $V$ replaced by $U_k$, which has orthonormal columns, so $\sum_{i=1}^{m_k}U_k(i,:)^TU_k(i,:)=I_r$. The excluded skeleton rows are harmless because $\|U_k(\mathcal{I}_{\le k},:)^TU_k(\mathcal{I}_{\le k},:)\|_2\le r\cdot\mu_U r/m_k\le\tfrac12$ by \cref{ass:incohU_app}, whence $\sum_{i\notin\mathcal{I}_{\le k}}U_k(i,:)^TU_k(i,:)\succeq\tfrac12 I_r$; the matrix Chernoff bound for uniform sampling without replacement \cite{tropp2012user,gross2010note} then gives the stated singular value bound with probability at least $1-\delta$.
\end{proof}
 
\subsection{Proof of \texorpdfstring{\cref{thm:alg3_main}}{the main theorem}}
 
\begin{proof}[Proof of {\cref{thm:alg3_main}}]
Start from the master bound \cref{eq:master_app}. By \cref{lem:trunc_app}, $e_k\le\|(I-P_{C_k^{\mathrm{s}}})A_k\|_F+\theta_k\Delta_k$, and the first term is bounded both by \cref{lem:skel_app} (deterministically, by $\beta_k$) and by \cref{lem:oversampling_app} (by $(1+c_1\delta^{-1/2})\Delta_k$, on an event $\mathcal{E}_k^{\mathrm{col}}$ of probability at least $1-2\delta$); since both bound the same quantity, their minimum does too. Hence, on $\bigcap_{k}\mathcal{E}_k^{\mathrm{col}}$,
\begin{equation}\label{eq:alg3_app_bound}
\|\mathcal{X}-\tilde{\mathcal{X}}\|_F\;\le\;
\sum_{k=1}^{d-1}\Bigl(\prod_{j=1}^{k}\eta_j\Bigr)
\Bigl[\,\min\bigl\{\beta_k,\;\bigl(1+c_1\delta^{-1/2}\bigr)\Delta_k\bigr\}\;+\;\theta_k\,\Delta_k\Bigr].
\end{equation}
This intermediate inequality still contains the oblique factors $\eta_j$; since they are computable once the $U_j$ are available, \cref{eq:alg3_app_bound} is an a posteriori certificate for the output of \cref{alg:TT_update_parallel}, valid with probability at least $1-2(d-1)\delta$ by a union bound over the $d-1$ column events. To obtain the explicit statement of \cref{thm:alg3_main}, apply \cref{lem:rowrandom_app}: on an event $\mathcal{E}_j^{\mathrm{row}}$ of conditional probability at least $1-\delta$ (conditioning on $\mathcal{S}_{j-1}$ and $\mathcal{L}_{>j}$), $\eta_j\le2\sqrt{m_j/p}\le2\sqrt{(r+p)n/p}$, so that $\prod_{j\le k}\eta_j\le\bigl(2\sqrt{(r+p)n/p}\bigr)^k$. Controlling the row events by sequential conditioning on $j=1,\dots,d-1$ contributes a further $(d-1)\delta$ to the failure probability, for a total of $3(d-1)\delta$. Substituting into \cref{eq:alg3_app_bound} gives \cref{eq:alg3_main_bound}; the simplified form follows by bounding each factor $\bigl(2\sqrt{(r+p)n/p}\bigr)^k$ by its value at $k=d-1$, $\min\{\cdot,\cdot\}\le(1+c_1\delta^{-1/2})\Delta_k$, $\theta_k\le\theta$, and $\Delta_k\le\max_k\Delta_k$.
\end{proof}
 
\begin{remark}[Interpretation and scope]\label{rem:alg3_remarks}
(i) Applying the Cauchy--Schwarz inequality to \cref{eq:alg3_app_bound} with the oversampling branch shows that, up to the oblique amplification and a constant depending on $\delta$, \cref{alg:TT_update_parallel} matches the TTSVD-type benchmark $\bigl(\sum_k\Delta_k^2\bigr)^{1/2}$ \cite{oseledets2011tensor}, which itself is within $\sqrt{d-1}$ of the best TT approximation error at ranks $r$. (ii) \Cref{lem:oversampling_app,lem:rowrandom_app} require $p\gtrsim\max\{\mu,\mu_U\}r\log r$ uniform samples, whereas the experiments in \cref{sec:numerical_tests} use much smaller $p$; the min in \cref{eq:alg3_main_bound} means the guarantee never falls below the deterministic skeleton branch, and the gap mirrors the state of the art for uniform sampling in the matrix literature \cite{chiu2013sublinear,cortinovis2024sublinear}. (iii) Alternatively to \cref{lem:rowrandom_app}, the row set for the oblique projections can be (re-)selected deterministically by maximum-volume or strong RRQR procedures applied to the small matrix $U_k$, with classical guarantees $\eta_k\le\sqrt{1+r(m_k-r)}$ for $r$ rows \cite{goreinov2010find} and improvements of the form $\sqrt{1+r(m_k-r)/(p+1)}$ for $r+p$ rows \cite{osinsky2018pseudo,osinsky2018rectangular,cortinovis2024sublinear}. This requires no incoherence hypothesis, and by \cref{lem:rowmono_app}, appending such rows to $\mathcal{S}_k$ inherits the same guarantee.
\end{remark}

\section{A baseline bound for the raw interpolatory TT cross}\label{app:ttaca_baseline}
For comparison, this appendix records a Chebyshev-norm bound for the \emph{raw} nested TT cross interpolant with cores built from $\hat A_k^{-1}=X_k(\mathcal{I}_{\le k},\mathcal{J}_{>k})^{-1}$, i.e., the interpolation-based format that our algorithms post-process~\cite{shi2026distributed}. Together with \cref{thm:PEID_format_main} and \cref{thm:alg3_main}, it quantifies at the format level what the projection enhancement and the oversampling buy, and it provides theoretical context for the reduction factors reported in \cref{sec:numerical_tests}.

Assume the pivots are nested as in \cref{eq:tt_nestedness} and each $\hat A_k$ is nonsingular. For $1\le k\le d-1$ define the interpolation operator $\Pi_k^{\mathrm{cr}}$ acting on tensors $\mathcal{Y}$ through the $k$th unfolding,
\begin{equation}\label{eq:Pi_cr_def}
\bigl(\Pi_k^{\mathrm{cr}}\mathcal{Y}\bigr)_k \;:=\; C_k\,\hat A_k^{-1}\, Y_k(\mathcal{I}_{\le k},:) .
\end{equation}
Since $C_k(\mathcal{I}_{\le k},:)=\hat A_k$, each $\Pi_k^{\mathrm{cr}}$ is an oblique projection. The nested TT cross interpolant is the TT tensor $\tilde{\mathcal{X}}^{\mathrm{cr}}$ with entries
\begin{equation}\label{eq:cross_product_formula}
\tilde{\mathcal{X}}^{\mathrm{cr}}(i_1,\dots,i_d)
= X_1(i_1,\mathcal{J}_{>1})\,\hat A_1^{-1}\,
X_2\bigl(\mathcal{I}_{\le1}\otimes\{i_2\},\mathcal{J}_{>2}\bigr)\,\hat A_2^{-1}\cdots
\hat A_{d-1}^{-1}\, X_{d-1}\bigl(\mathcal{I}_{\le d-1},\{i_d\}\bigr),
\end{equation}
the standard TT-cross formula \cite{oseledets2010tt,savostyanov2014quasioptimality,dolgov2020parallel}, which under \cref{eq:tt_nestedness} interpolates $\mathcal{X}$ exactly on the pivot crosses.

\begin{lemma}[Product representation]\label{lem:cross_product_app}
$\tilde{\mathcal{X}}^{\mathrm{cr}} = \Pi_1^{\mathrm{cr}}\,\Pi_2^{\mathrm{cr}}\cdots\Pi_{d-1}^{\mathrm{cr}}\,\mathcal{X}$.
\end{lemma}

\begin{proof}
Set $W^{(d)}:=\mathcal{X}$ and $W^{(k)}:=\Pi_k^{\mathrm{cr}}W^{(k+1)}$ for $k=d-1,\dots,1$. By \cref{eq:Pi_cr_def}, $W^{(1)}_1 = C_1\hat A_1^{-1}\,W^{(2)}_1(\mathcal{I}_{\le1},:)$. The entries of $W^{(2)}=\Pi_2^{\mathrm{cr}}W^{(3)}$ with first index $i_1\in\mathcal{I}_{\le1}$ form the rows $\mathcal{I}_{\le1}\otimes\mathbb{I}_2$ of the second unfolding $W^{(2)}_2=C_2\hat A_2^{-1}W^{(3)}_2(\mathcal{I}_{\le2},:)$; restricting those rows on both sides and using \cref{eq:slice_identity} gives
$W^{(2)}_2(\mathcal{I}_{\le1}\otimes\mathbb{I}_2,:) = M_2\,\hat A_2^{-1}\,W^{(3)}_2(\mathcal{I}_{\le2},:)$. Iterating for $k=3,\dots,d-1$ and terminating with $W^{(d)}_{d-1}(\mathcal{I}_{\le d-1},:)=X_{d-1}(\mathcal{I}_{\le d-1},:)$ yields exactly \cref{eq:cross_product_formula}.
\end{proof}

\begin{lemma}[Interpolation operator norm under complete pivoting]\label{lem:kappa_app}
Let $(\mathcal{I},\mathcal{J})$ be produced by $r$ steps of complete pivoting applied to $A\in\R^{M\times N}$, with $\hat A=A(\mathcal{I},\mathcal{J})$ nonsingular and $C=A(:,\mathcal{J})$. Then every entry of $C\hat A^{-1}$ is bounded in modulus by $2^{r-1}$ and $\kappa:=\|C\hat A^{-1}\|_\infty\le r\,2^{r-1}$, where $\|\cdot\|_\infty$ is the maximum absolute row sum.
\end{lemma}

\begin{proof}
After permuting rows and columns into pivot order, $r$ steps of LU with complete pivoting give $A(:,\mathcal{J})=\begin{bmatrix}L_{11}\\ L_{21}\end{bmatrix}U_{11}$ and $\hat A=L_{11}U_{11}$ (the remainder vanishes on pivot columns), so
$C\hat A^{-1}=\begin{bmatrix}I_r\\ L_{21}L_{11}^{-1}\end{bmatrix}$ in that ordering. Let $D$ collect the pivots (the diagonal of $L_{11}$) and set $\tilde L_{11}:=L_{11}D^{-1}$, $\tilde L_{21}:=L_{21}D^{-1}$; then $L_{21}L_{11}^{-1}=\tilde L_{21}\tilde L_{11}^{-1}$. Complete pivoting bounds all multipliers by one, i.e., $\tilde L_{11}$ is unit lower triangular and all entries of $\tilde L_{11},\tilde L_{21}$ have modulus at most one. For such $\tilde L_{11}$, the classical bound $|(\tilde L_{11}^{-1})_{ij}|\le 2^{i-j-1}$ for $i>j$ (see, e.g., \cite[Chap.~8]{higham2002accuracy}) yields $\|\tilde L_{11}^{-1}\|_\infty\le 2^{r-1}$, each entry of $\tilde L_{21}\tilde L_{11}^{-1}$ is bounded by $2^{r-1}$, and $\|\tilde L_{21}\tilde L_{11}^{-1}\|_\infty\le r\,2^{r-1}$. Row and column permutations do not change $\|\cdot\|_\infty$ up to reordering rows.
\end{proof}

\begin{theorem}[Raw TT cross with complete pivoting]\label{thm:ttaca_app}
Suppose the pivots are nested as in \cref{eq:tt_nestedness} and each pivot pair $(\mathcal{I}_{\le k},\mathcal{J}_{>k})$ is produced by $r$ steps of ACA with complete pivoting applied to $X_k$, where $\rank(X_k)>r$. Then
\[
\|\mathcal{X}-\tilde{\mathcal{X}}^{\mathrm{cr}}\|_{\max}\;\le\;\sum_{k=1}^{d-1}\bigl(r\,2^{r-1}\bigr)^{k-1}\,4^{r}\rho_r\,\sigma_{r+1}(X_k)
\;\le\;(d-1)\,\bigl(r\,2^{r-1}\bigr)^{d-2}\,4^{r}\rho_r\;\max_{1\le k\le d-1}\sigma_{r+1}(X_k).
\]
Since $\sigma_{r+1}(X_k)\le\|\mathcal{X}-\mathcal{Y}\|_F$ for every TT tensor $\mathcal{Y}$ with TT ranks at most $r$ \cite{oseledets2011tensor}, the interpolant is quasi-optimal with constant $(d-1)(r2^{r-1})^{d-2}4^r\rho_r$ relative to the best TT approximation in the Frobenius norm.
\end{theorem}

\begin{proof}
By \cref{lem:cross_product_app} and the telescoping identity
$I-\Pi_1^{\mathrm{cr}}\cdots\Pi_{d-1}^{\mathrm{cr}}
=\sum_{k=1}^{d-1}\Pi_1^{\mathrm{cr}}\cdots\Pi_{k-1}^{\mathrm{cr}}(I-\Pi_k^{\mathrm{cr}})$,
\[
\|\mathcal{X}-\tilde{\mathcal{X}}^{\mathrm{cr}}\|_{\max}\le\sum_{k=1}^{d-1}\Bigl(\prod_{j=1}^{k-1}\|\Pi_j^{\mathrm{cr}}\|_{\max\to\max}\Bigr)\,\|(I-\Pi_k^{\mathrm{cr}})\mathcal{X}\|_{\max}.
\]
For any tensor $\mathcal{Y}$, entry $(a,b)$ of $(\Pi_j^{\mathrm{cr}}\mathcal{Y})_j$ equals $\sum_{s}(C_j\hat A_j^{-1})_{as}\,Y_j(i_s,b)$, whence $\|\Pi_j^{\mathrm{cr}}\|_{\max\to\max}\le\|C_j\hat A_j^{-1}\|_\infty\le r2^{r-1}$ by \cref{lem:kappa_app}. Furthermore, $\bigl((I-\Pi_k^{\mathrm{cr}})\mathcal{X}\bigr)_k$ is exactly the remainder of $r$ steps of complete-pivoting ACA on $X_k$, so \cref{lem:rectCKM_app} bounds its Chebyshev norm by $4^r\rho_r\sigma_{r+1}(X_k)$.
\end{proof}

\begin{remark}[Sharpness and scope]\label{rem:sharpness_app}
(i) The factor $4^r\rho_r$ is inherited from the matrix bound and cannot be substantially improved at this level of generality: already for symmetric positive semidefinite matrices the factor $4^r$ is asymptotically attained \cite{cortinovis2020maximum,harbrecht2012low}, with Kahan-type matrices \cite{kahan1966numerical} providing explicit certificates on which complete-pivoting ACA overestimates $\sigma_{r+1}$ by factors growing exponentially in $r$. Such failures are pivot-selection failures rather than subspace failures: the trailing singular subspaces of these matrices are well spread out, so they are precisely the instances where the oversampling branch of \cref{eq:alg3_main_bound} beats the skeleton branch exponentially. (ii) With maximum-volume pivots instead, the dominance property gives $\kappa_k\le r$ and per-separation errors at most $(r+1)\sigma_{r+1}(X_k)$ \cite{goreinov2001maximal}, recovering, up to constants, the Chebyshev-norm quasi-optimality results of \cite{savostyanov2014quasioptimality,osinsky2019tensor}. (iii) The complete-pivoting assumption is an idealization: practical linearly-scaled TTACA uses partial or heuristic pivoting, for which no general a priori bound can exist, and complete pivoting per unfolding does not automatically produce nested sets. Comparing \cref{thm:ttaca_app} with \cref{thm:PEID_format_main,thm:alg3_main} shows what post-processing buys at the level of the bounds: projection enhancement replaces the interpolation amplification $r2^{r-1}$ per separation by the coefficient norms $\gamma_k$ (or their computable surrogates $\eta_k$), with per-separation errors that are optimal for the chosen fibers, and oversampling further replaces the exponential per-separation constant, in probability, by an absolute constant times the optimal error $\Delta_k$.
\end{remark}

\bibliographystyle{elsarticle-num}
\bibliography{references}

\end{document}